\documentclass{scrartcl}
\usepackage{mathtools, amsmath, amssymb, amsthm, float}
\usepackage{dsfont}
\usepackage{bbm,bm}
\usepackage[scr=boondoxo]{mathalpha}
\usepackage{stmaryrd}
\usepackage[style=alphabetic, giveninits=true, maxbibnames=99]{biblatex}
\usepackage{algorithm}
\usepackage{algpseudocodex}
\algrenewcommand{\algorithmiccomment}[1]{\textcolor{gray}{\itshape\(\triangleright\) #1}}
\usepackage{mathtools}

 \usepackage[font={small}]{caption}
\RequirePackage[OT1]{fontenc}
 \usepackage{lscape}

\usepackage[textsize=footnotesize]{todonotes}
\usepackage{tikz}
\usepackage{pgfplots}
\usetikzlibrary{intersections, pgfplots.fillbetween}
\usetikzlibrary{arrows.meta,backgrounds}
\usepackage[textsize=footnotesize]{todonotes}
\usepackage[left=3cm, right=2.5cm, top=2.5cm, bottom=2.5cm]{geometry}

\newcommand{\bfxi}{\boldsymbol{\xi}}
\newcommand{\bfzeta}{\boldsymbol{\zeta}}
\newcommand{\bfX}{\mathbf{X}}

\newcommand{\bfm}{\mathbf{m}}
\newcommand{\bfk}{\mathbf{k}}

\newcommand{\cA}{\mathcal {A}}
\newcommand{\cB}{\mathcal {B}}

\newcommand{\cE}{\mathcal {E}}

\newcommand{\cG}{\mathcal {G}}
\newcommand{\cH}{\mathcal {H}}
\newcommand{\cI}{\mathcal {I}}
\newcommand{\cK}{\mathcal {K}}
\newcommand{\cL}{\mathcal {L}}

\newcommand{\cN}{\mathcal {N}}

\newcommand{\cP}{\mathcal {P}}

\newcommand{\cR}{\mathcal {R}}
\newcommand{\cS}{\mathcal {S}}
\newcommand{\cT}{\mathcal {T}}
\newcommand{\cU}{\mathcal {U}}
\newcommand{\cV}{\mathcal {V}}

\newcommand{\cZ}{\mathcal {Z}}

\newcommand{\IE}{\mathbb{E}}

\newcommand{\IK}{\mathbb{K}}

\newcommand{\IM}{\mathbb{M}}
\newcommand{\IN}{\mathbb{N}}

\newcommand{\IP}{\mathbb{P}}

\newcommand{\IR}{\mathbb{R}}

\newcommand{\mem}{\mathscr{m}}
\newcommand{\comm}{\mathscr{c}}
\newcommand{\offmem}{M}
\newcommand{\offcomm}{K}
\newcommand{\compgraph}{\cK}
\newcommand{\stargraph}{\cS}
\newcommand{\compline}{\cL_K}

\newcommand{\BC}{\mathrm{BC}}
\newcommand{\RIG}{\mathrm{RIG}}

\newcommand{\CT}{\mathrm{CT}}

\newcommand{\at}{\mathrm{aT}}
\newcommand{\one}{\mathbf{1}}
\newcommand{\1}{\mathbbm{1}}

\renewcommand{\root}{o}

\newcommand*{\R}{\mathbb R}
\newcommand{\Z}{\mathbb{Z}}
\newcommand{\E}{\mathbb{E}}
\renewcommand{\P}{\mathbb{P}}

\newcommand{\bitp}{\boldsymbol{p}}
\newcommand{\bitq}{\boldsymbol{q}}

\newcommand{\leftpar}{\cV^{\mem}}
\newcommand{\rightpar}{\cV^{\comm}}

\newtheorem{theorem}{Theorem}[section]
\newtheorem{metatheorem}[theorem]{Meta-theorem}
\newtheorem{lemma}[theorem]{Lemma}
\newtheorem{corollary}[theorem]{Corollary}
\newtheorem{definition}[theorem]{Definition}

\newtheorem{remark}[theorem]{Remark}
\newtheorem{assumption}[theorem]{Assumption}
\newtheorem{proposition}[theorem]{Proposition}
\newtheorem{claim}[theorem]{Claim}

\newtheorem{example}[theorem]{Example}

\makeatletter
\renewenvironment{proof}[1][\proofname]{%
	\par\pushQED{\qed}\normalfont%
	\topsep6\p@\@plus6\p@\relax
	\trivlist\item[\hskip\labelsep\bfseries#1\@addpunct{.}]%
	\ignorespaces
}{%
	\popQED\endtrivlist\@endpefalse
}
\makeatother
\numberwithin{equation}{section}
\title{How communities help epidemics persist}
\author{J\'ulia Komj\'athy\footnote{Delft University of Technology; Mekelweg 4, 2628 CD Delft, Netherlands, \texttt{j.komjathy@tudelft.nl}}, \,Marco Seiler\footnote{Goethe University Frankfurt; Robert-Mayer-Straße 10, 60486 Frankfurt am Main, Germany, \texttt{seiler@math.uni-frankfurt.de}}  \,and Daniel Valesin\footnote{University of Warwick; Coventry, CV4 7AL, United Kingdom, \texttt{daniel.valesin@warwick.ac.uk}}}
\date{\today}

\begin{document}
\maketitle

\begin{abstract}
   Empirical studies show that communities such as households and workplaces strongly influence how an infection spreads and how long it persists in a population. To understand this effect, we study the contact process (SIS epidemic model) on random intersection graphs and their local weak limit. In the contact process, infected vertices recover with rate~$1$ and transmit the infection to each of their neighbours with rate~$\lambda$. The random intersection graph is defined on a finite population; each individual belongs to a random number of communities (captured by a membership distribution), and each community has a random size (captured by a community size distribution). Communities are all taken to be complete graphs. The local weak limit of the random intersection graph, called the branching community tree, has complete graphs arranged in a tree structure. 

   Letting~$M$ have the membership distribution and~$K$ have the community size distribution, we prove that if either~$\mathbb E[e^{cM}]=\infty$ for all~$c > 0$ or~$\mathbb E[e^{cK\log K}] = \infty$ for all~$c > 0$, then for all~$\lambda > 0$, the infection on the finite graph survives for exponentially long with high probability, and on the local limit graph, it survives forever with positive probability. On the other hand, if for some~$c > 0$ both~$\mathbb E[e^{cM}]<\infty$ and~$\mathbb E[e^{K\log(cK)}]<\infty$ hold, then the process with small~$\lambda$ dies out in the local limit graph. We compare these findings with  known results on the behaviour of the contact process on a ``community-free'' random graph
with similar degree distributions, to highlight the impact of communities.
\end{abstract}
{\footnotesize\textit{2020 Mathematics Subject Classification --} Primary 60K37; Secondary 05C82, 92D30  
	\\
\textit{Keywords --}  Contact Process; Metastability; Random Graphs; Random Intersection Graph}

\section{Introduction}
The communities we are part of determine our life to a large extent. Our extended families, cultural-, sports-clubs, our friends, our colleagues all form our social circles. It has been debated whether epidemics spread faster or slower due to the presence of these modular structures, or whether the presence of such structures increases or decreases pandemic thresholds \cite{hui2007modular, nadini2018epidemic, hebert2010propagation}. In this paper we rigourously study the contact process (also called the Susceptible-Infected-Susceptible or SIS model) on a random network model incorporating the presence of communities. The contact process is a continuous time Markov process that we informally describe now. On any network $G$ with vertices representing individuals, each vertex can be either infected or susceptible  at any given time $t\ge 0$. When infected, vertices heal independently of everything else at rate $1$, and while infected, they send the infection to each of their susceptible neighbours  at rate $\lambda>0$, where $\lambda$ is called the infection rate. On a large but finite graph $G$, the infection eventually vanishes as the all-healthy configuration is an absorbing state. Due to sustained infections, this could take as long as exponentially long in the size of the graph. In this case we talk about \emph{long persistence} of the epidemic in the population. We are motivated by the question:
\begin{equation*}
\begin{aligned}
&\textit{How does the presence of communities in the underlying contact network}\\
&\qquad \qquad\textit{affect the critical infection rate for long persistence?}
\end{aligned}
\end{equation*}

This question is nontrivial and the answer has been debated, as we explain now. While communities may support repeated reinfections, rare connections between large-enough communities may cause bottlenecks and thus prohibit spreading. Earlier studies report increased SIS epidemic thresholds in certain community-structured network models: Zhao and Gao \cite{hui2007modular} find a threshold increasing with modular strength in a small-world model, while the mean-field analysis of H\'ebert-Dufresne \emph{et al.} \cite{hebert2010propagation} predicts a higher threshold than in random networks with the same degree distribution, under the assumption of zero degree correlations. We mention that these models do not easily translate to the model we use for our results.

For SIS, it has been observed \cite{cota2018griffiths, onge2021master, onge2021social} that communities may create long persistence, due to the following effect. As reinfections within a single dense community are frequent, the infection may be sustained in communities above a certain $\lambda$-dependent size, while it is absent in most of the network. If such communities form a connected communicating network, these rare regional effects may cause sustained infections for SIS. Onge \emph{et al.\ }call this mesoscopic localisation, and estimate the necessary infection rate for this to happen by replacing the communication between different communities, (i.e., the effect of incoming rare infection to a single community) by a mean field external force field.  This gives a prediction for the pandemic threshold in terms of community sizes and membership distributions for individuals, which shall come closest to our mathematical results. 

To be able to answer the question rigourously, we choose the \emph{random intersection graph (RIG)} to model the the underlying  population structure. This choice allows for easy comparison with existing literature \cite{bhamidi2021survival, onge2021master}. We shall follow the model variant proposed in \cite{newman2003properties, coupechoux2014clustering, van2021random}, a variant of RIG that incorporates both tunable community-distributions and tunable membership distributions; as real networks often show a high variability of their degrees \cite{newman2018networks}. Variations on similar network models incorporating communities are plentiful \cite{coupechoux2014clustering, newman2003properties, singer1996random, ball2014epidemics}, and our results and some of the analysis we provide carries through these network models either directly or with minor modifications. 

In the random intersection graph model we prescribe, for each vertex, how many communities it shall be a member of, stored in a membership list $\mathbf m^{(n)}=(m_v)_{v\in \mathcal V}$, in total $h_n=\sum_{v\in V} m_v$ memberships. Further, we also prescribe for each community how many members it shall contain, stored in a community size list $\mathbf k^{(n)}=(k_1,k_2, \dots k_{a_n})$, where $a_n$ is the total number of communities $\mathcal C$ in the model, and in total we must then have again $h_n=\sum_{w\in \mathcal C}k_w$ community memberships. We then match the two different types of memberships uniformly at random using a random permutation. We finally create the RIG graph by superimposing each community $w$ as a complete graph on its members in the final graph. For a precise definition, see Definition \ref{def:random-intersection-graph} below and for a detailed explanation and properties such as degree distribution, clustering, community overlaps, and local weak limit of this model see \cite{van2021random}, while \cite{van2022phase} for the phase transition for the existence of a unique giant (linear sized) component. We shall assume that the empirical distribution of both the individual memberships and the community sizes converge, together with their $\eta$th moment, for some $\eta>1$, to those of some limiting variables $M$ and $K$.

The local weak limit of the random intersection graph captures the limit of the graph as seen from a typical vertex \cite{benjamini2001recurrence, aldous2004objective}. For random intersection graphs, this is an infinite graph that we call a \emph{branching community tree}, which inherits some tree structure from the underlying bipartite model of matching the memberships of individuals and communities, while containing communities chained together in a tree structure. The local weak limit was identified in \cite{van2021random}, and we describe it in Section \ref{sec:Definition-community-tree}, see also Figure \ref{fig:BCLWL}. As the limiting graph is infinite with positive probability, we may speak about local versus global survival of the process on it. We speak about local survival when the infection returns to any vertex infinitely many times after consecutive healings. 

We are ready to state our main results and answer the heuristic question. We chose to do this informally in a meta-theorem. 
\begin{metatheorem}\label{thm:meta} Consider the SIS/contact process on a random intersection graph with limiting individual membership distribution $M$ and community size distribution $K$, respectively. 
\begin{enumerate}
\item If either $\mathbb E[e^{cM}]=\infty$ for all $c>0$  or $\mathbb E[e^{c K\log K}]=\infty$ for all $c>0$, then the contact process shows sustained infection for all $\lambda>0$ for all sufficiently large network size $n$; i.e., the critical infection rate is $0$. The infection is metastable in the finite graph, and is sustained exponentially long with respect to the number of vertices in the graph. Both the local and global critical infection rates are $0$ on the limiting branching community tree. 
\item If for some $c>0$ both $\E[e^{cM}]<\infty$ and $\E[e^{K\log(cK)}]<\infty$ hold, then the contact process has critical infection rate $>0$ in the large network limit. Both the local and global critical infection rates are strictly positive on the limiting branching community tree.
\end{enumerate}
\end{metatheorem}
For the precise statements, see Theorems \ref{thm:Survival-Criterion-Community-Tree} and \ref{thm:Existence-of-Subcritical-Phase} for the results on the local weak limit and Theorems \ref{thm:survival-finite-graph} and \ref{thm:rigc-extinction} for exponentially long persistence on finite graphs. 

\paragraph{Long persistence even via light tailed communities.} Next we contrast our results to existing literature. Metastability is the phenomenon where a process persists for a long time in a set of states that seem stationary, but are not true equilibrium states, until the process finally enters its true equilibrium. 
The long persistence of the contact process in a finite population is a form of metastability. It has been studied on multiple graph structures, for example on finite subsets of $\Z$ \cite{durrett1988finite,durrett1988finiteII}, or on regular trees \cite{stacey2001trees,cranston2014trees}, and more recently on random graphs \cite{chatterjee2009power,mourrat2016phase}.
 For the contact process on random graphs model mimicking real networks specifically, in the mathematical literature, Huang and Durrett \cite{huang2020contact} followed by Bhamidi \emph{et.\ al.}\ \cite{bhamidi2021survival} showed for the \emph{configuration model}, that the asymptotic degree distribution having a finite moment generating function in some open neighborhood of $0$ is both necessary and sufficient for the critical infection rate to be zero. The configuration model \cite{bollobas1980probabilistic, molloy1995critical} is a community-free version of our model: individuals have prescribed degrees $(d_1, \dots, d_n)$, which, upon being represented as half-edges,  are matched uniformly into edges to form a graph. Such a graph has asymptotically vanishing clustering and no well-defined communities under finite average degree assumption.  
With $D$ being the asymptotic degree distribution, $\E[e^{cD}]=\infty$ for all $c>0$ implies that the infection persists for all $\lambda>0$ with positive probability for sufficiently large network size $n$, while $\E[e^{cD}]<\infty$ for some $c>0$ causes the infection to go extinct within polynomial time in the network size \cite{bhamidi2021survival}.  

The asymptotic degree distribution of the random intersection graphs can be written as $D=\sum_{i=1}^M \widetilde K$, where $\widetilde K$ is the size-biased  distribution of $K$ shifted down by $1$. Our result in Meta-theorem \ref{thm:meta} says that the contact process persists in the network when either $\E[e^{cM}]=\infty$ or $\E[e^{cK\log K}]=\infty$ for all $c>0$. Assuming bounded community sizes $K<\Delta$ for some constant $\Delta<\infty$, $D\le \Delta M$ and our theorem recovers the results in \cite{bhamidi2021survival,huang2020contact}. Indeed, then finiteness vs infiniteness of $\E[e^{cD}]$ is equivalent to that of $\E[e^{cM}]$, meaning that in this case communities are not able to help the infection enough to achieve critical rate $0$, and survival must come from high degree vertices captured by the tail of $M$. 
The effect of communities on the critical rate becomes noticeable when communities  have a light-tailed but unbounded distribution, for instance, following a geometric distribution with arbitrary fixed success parameter $p>0$. Assuming bounded distribution for community memberships $M$, in this case, the tail of $D=\sum_{i=1}^{M}\widetilde K_i$ is exponential and $\E[e^{cD}]<\infty$ for some $c>0$, so \cite{bhamidi2021survival} would predict extinction for small infection rates $\lambda>0$. As, however, $\E[e^{c K\log K}]=\infty$ in this case for all $c>0$, our results modify this to say that the effect of these communities is strong enough to achieve long persistence on finite graphs and eternal survival of the infection on the limiting branching community trees, for any fixed $\lambda>0$.

 \paragraph{Extinction for extremely light tailed communities.} If, however, the limiting community size distribution is Poisson with an arbitrary fixed mean $\mu>0$, then $\E[e^{K \log (cK)}]<\infty$ is finite for all $c<1/(e\mu)$, so if additionally also $\E[e^{cM}]<\infty$ then we see a subcritical phase for small $\lambda>0$ for the contact process. In some sense Poisson is exactly on the threshold, as its tail is asymptotically $\P(K\ge k)\asymp \mu^k/k!\asymp \exp(-k \log(k)+(\log(\mu)+1)k -\frac12 \log k)$.
 We can contrast this result with nonrigourous findings from physics. 
Onge \emph{et.\ al.} \cite{onge2021master} predicts that SIS is subcritical if $\lambda$ satisfies
\begin{equation}\label{eq:physics-formula-0}
\lambda\frac{\mathbb E[M(M-1)]}{\mathbb E[M]}\cdot \frac{\mathbb E[K(K-1)A(K,\lambda)]}{\mathbb E[K]}<1,
\end{equation}
where $\lambda(k-1)A(k, \lambda)$ captures the number of infection events generated by a single seed inside a clique of size $k$, with   
\begin{equation}\label{eq:physics-formula-1}
A(k,\lambda)=\int_{0}^{\infty}(1+\lambda u)^{k-2} e^{-u}\mathrm du=\lambda^{k-2}(k-2)!
\sum_{\ell=0}^{k-2}\frac{\lambda^{-\ell}}{\ell!}.
\end{equation}
The authors derive this formula by replacing incoming infections to cliques by a homogeneous external force field.
While the goal of  \cite{onge2021master} is to give numerical predictions on the critical infection rate, we can use this equation to check whether the critical rate is strictly positive in the large network limit as follows. 
For every fixed $\lambda>0$, we can bound
\[
\lambda^{k-2}(k-2)!
\le A(k,\lambda)
\le e^{1/\lambda}\lambda^{k-2}(k-2)!
\]
Using this in \eqref{eq:physics-formula-1}, and assuming bounded $M$ to avoid hub-effects,   
the formula in \eqref{eq:physics-formula-0} predicts that if for some small $\lambda_0>0$ we have $\E[\lambda_0^{K}K!]<\infty$, then for sufficiently small $\lambda>0$  we shall see subcritical contact process behaviour regardless of the network size $n$, implying that the critical infection rate is strictly positive in the large network limit. Using elementary analysis one readily shows that $\E[\lambda_0^{K}K!]<\infty$ for some $\lambda_0>0$ is equivalent of the finiteness of $\E[e^{K \log(cK)}]$ for some $c>0$ in our Meta-theorem \ref{thm:meta}. In other words, our mathematically rigorous method reproduces the latest physics predictions on the non-triviality of the critical infection rate $\lambda_{1}$ caused by communities.
We mention that there is a small gap left open between the two criteria in part 1 and part 2 of Meta-theorem \ref{thm:meta}, which we investigate further in Section \ref{sec:lwl-extinction}.
\paragraph{Methodology.} For the nontriviality of $\lambda_1>0$, i.e., part 2 of Meta-theorem \ref{thm:meta}, we use a probabilistic, -- rather than analytic, as in \cite{onge2021master} -- method, and stochastically dominate the contact process by a varying-rate \emph{reinforced contact process} on the local weak limit branching community tree. The branching community tree is a decorated tree where additional edges between individuals in the same generation are present to form the communities.  The reinforced contact process, invented in \cite{bhamidi2021survival}, is amenable to tree-recursion methods:  infections `upward' in the tree are replaced by keeping the whole ancestry line of every infected individual infected as long as the individual is infected. To compensate for the effect of within generation edges, we modify the original reinforced contact process by increasing the rate on each downward edge by a factor of the number of within generation edges of the receiving vertex. Using this construction, we can show that the obtained inhomogeneous-rate reinforced contact process has expected extinction time bounded from above by $2$ for all small $\lambda>0$, when the infection is started from the root, whenever $\E[e^{K \log (c K)}]<\infty$ for some $c>0$.
We find it remarkable that this upper-bounding method matches the non-rigorous predictions in \cite{onge2021master}.
See Section \ref{sec:lwl-extinction} for more details.  

We prove long persistence (Meta-theorem \ref{thm:meta}, part 1) via two different methods; one for the finite graphs and one for the limit branching community tree,  both using embedding. In the branching community tree, we show that it is a.s. possible 
to embed an infinite chain of communities of a certain infection-rate-dependent size $k_\lambda$, 
at bounded distance $r_\lambda$ from each other, as indicated on Figure \ref{fig:community-line}.
This embedding is possible for arbitrary small infection rates $\lambda>0$ under the condition $\E[e^{cK\log K}]=\infty$. 
There are two competing effects that we explain: we show that the infection in any given community
persist $T_{\mathrm{per}}(\lambda):=\exp(k_\lambda \log (\lambda k_\lambda))$ long with probability $1-o_\lambda(1)$.
Within this time, the community needs to communicate with other similar communities. Transmission is possible whenever $\lambda^{r_\lambda} T_{\mathrm{per}}(\lambda)\gg 1$.
Within distance $r_\lambda$, one finds roughly $\mu^{r_\lambda}\P(K\ge k_\lambda)\gg 1$ many communities of such size; where $\mu$ is the exponential growth rate of the branching community tree; giving a condition on admissible pairs $(k_\lambda, r_\lambda)$. 
For large $\lambda$, it is possible to find admissible pairs under even milder conditions; however, 
$\E[e^{cK\log K}]=\infty$ for all $c>0$ is necessary for a solution to exist for all $\lambda>0$.

The survival of the local weak limit graph does not necessarily guarantee long persistence on finite graphs, as the objects necessarily for such survival may not be present in the finite graphs at the relevant scales. Most of the paper is dedicated to show that this is not the case, and we can push persistence to \emph{exactly exponential time-scale} in the size of the network.  
For showing that the healing of the whole graph happens on \emph{at most} an exponential scale, we develop what we call \emph{the method of non-contributions} per vertex. This method is robust, and it allows us to show extinction of the infection on strictly exponential time-scales even if the graph is not sparse in the strict sense; i.e., the number of edges is superlinear in $n$. 
This naturally happens when the community-size distribution follows a power-law with $\P(K=k)\propto k^{-\tau_c}, \tau_c \in (2,3)$, as already the largest communities in the graph are of size $n^{1/(\tau_c-1)}$, giving rise to $n^{2/(\tau_c-1)}\gg n$ many edges, which is also the order of the total number of edges in the graph. See Section \ref{sec:long-survival-on-rig}.

\begin{figure}[!t]
  \centering
  \includegraphics[width=0.95\textwidth]{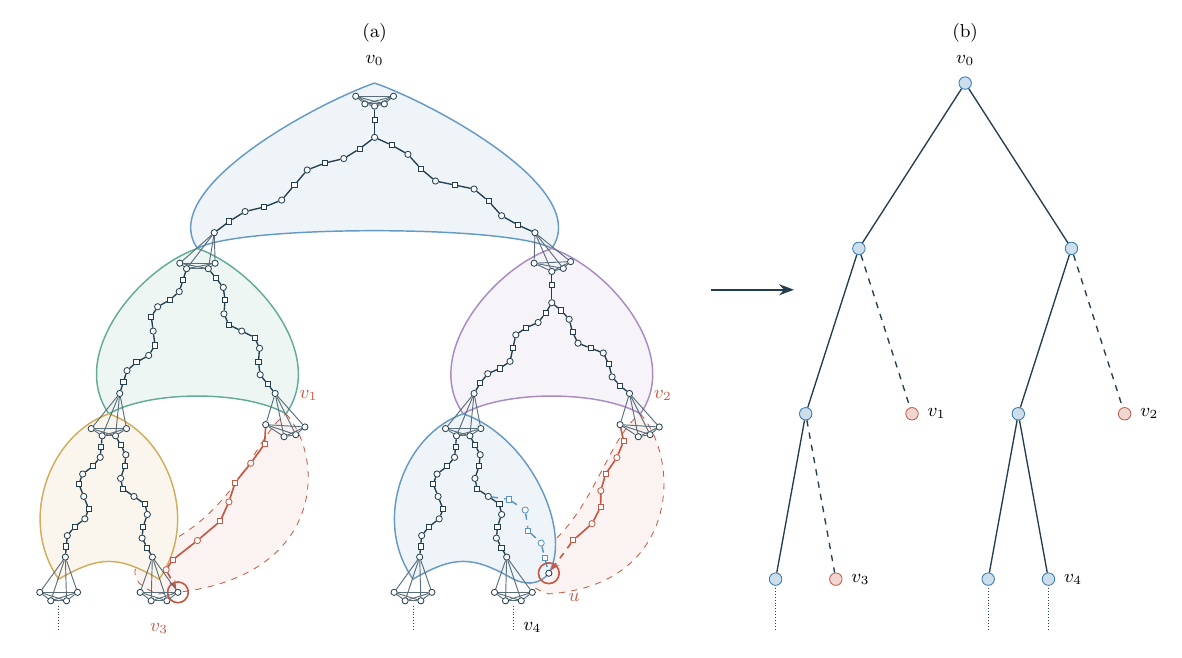}
  \caption{Subfigure (a) shows part of the resulting bead-network, while Subfigure (b) visualises the resulting branching structure. Blue vertices in (b) represent successful beads and red beads signify an unsuccessful exploration round. The core failure criteria, beside some technical ones, is that a loop is found in an exploration round. In the round where we build $v_1$'s bead, one can see in (a) that the breadth-first exploration finds an connecting path between $v_1$ and $v_3$, where $v_3$ was already found in a prior round. This results in a short loop that connects $v_1$ and $v_3$, and thus we take neither of these vertices to be part of the desired network of communicating communities. A loop is found also in the round where we build the bead of $v_2$. There, the loop does not contain the initial community of any other bead. Thus, it suffices to only delete $v_2$ and the connecting path, whereas $v_4$ shall be part of the network later if its bead can be built successfully later.}
  
  \label{fig:bead-embedding}
\end{figure}

Finally, the technically hardest part of the proof is to show that long persistence does occur with positive probability on \emph{at least} exponential time-scale in $n$. 
We identify a linear-sized \emph{network} of communicating communities, each of size $k_\lambda$, where neighboring communities in
the network are within distance $r_\lambda$ at most in the original graph. 
We do this via an embedding algorithm in Section \ref{sec:long-survival-on-rig} that jointly constructs the graph together with the community-network. 
The algorithm is coupling-based, and aims to embed linearly many contant-sized neighborhoods of size $k_\lambda$ communities, or `beads'. Each `bead' represents a successful coupling of the exploration of the bounded-sized neighborhood of a community of size $k_\lambda$ to an independent sample from the local weak limit started at a community of size $k_\lambda$. The beads are concatenated in a way that at the boundary of each bead the next two communities of size $k_\lambda$ are present to explore, so that they can start their own beads. As each exploration is of constant size, the success  probability of the coupling can be made sufficiently high (say $>9/10$) throughout the building of the entire bead-network. The whole bead-network then forms a branching structure of small `beads' tied together at the communities we have found, with dead ends at unsuccessful couplings (due to e.g. exploration loops).  See Figure \ref{fig:bead-embedding}. We show that we can successfully embed linearly many such `beads', forming the desired web. As every bead is connected to at most $3$ more other beads, the coarse structure of the web is a linear sized bounded degree tree, and we show that the contact process persists on this structure exponentially long in $n$. The large scale architecture of our algorithm is somewhat similar in nature to the one in \cite{mountford2016exponetial} in the sense that the embedding algorithm in that paper also embeds small structures similar to our beads. The proof there is not coupling based and the errors that arise as possible causes of failure are different for two reasons, the graph model is different and there is no coupling to a local weak limit.

\paragraph{Further related literature.}
Metastability, equivalently long persistence, of the contact process has been a active area in last two decades, especially on graphs where
the aysmptotic degree distribution $D$ is heavy-tailed. 
In an early series of works \cite{pastor2001epidemic,pastor2001complex,pastor2002finite}, Pastor-Satorras \& Vespignani predicted via non-rigorous mean-field approximations, that if $\P(D=k)\propto k^{-\tau}$, then for small enough infection rates $\lambda>0$ not depending on the network size, the disease persists for an exponentially long time if $\tau>3$, 
while it will die out after polynomial time (in the network size) if $\tau<3$. 
This was later disproved by Chatterjee and Durrett \cite{chatterjee2009power}, who showed that for \emph{all} $\tau<\infty$  for any small fixed $\lambda>0$ and sufficiently large network size $n$, the persistence time is superpolynomial with respect to the graph size. This leads to the limiting infinite graph having $0$ critical infection rate for survival.
This was further improved by Mountford et al.\@~\cite{mountford2016exponetial}, who showed that the persistence time is exponential, for all $\tau\in (2,\infty)$ and by \cite{can2015configuration} who showed the same if $\tau\in (1,2]$. Finally,  \cite{bhamidi2021survival} showed the same for distributions with no exponential moments. Besides the configuration model, metastability of the contact process has also been studied on many other scale free network models, such as preferential attachment models, hyperbolic graphs and scale-free percolation, see \cite{van2017preferential,linker2021hyperbolic,barnier2025scale}. See~\cite{valesin2024contact} for a survey. Recently, Jahnel, L\"{u}chtrath, and M\"{o}nch \cite{jahnel2026sparse} provided a generalisation of \cite{chatterjee2009power} for any sequences of sparse random graphs that converge in the local weak sense.

It is also noteworthy that in spatial models for random graphs, cliques play similar roles as communities in our model. In the random geometric graph model, where vertices form a unit intensity Poisson process in a box of volume $n$ and vertices are connected whenever their distance is at most $r_n$, complete graphs of size $\Theta(r_n^d)$ proliferate the graph. Van Hao~\cite{van2018super} showed that these structures ensure persistence time as long as $\exp(n \log (\lambda r_n^d))$. For sparse graphs ($r_n=r$ fixed) this gives exponential survival above the critical infection rate $\lambda_1$, which is strictly positive in this particular model \cite{menard2016percolation}. While random geometric graphs have a different community overlap structure than random intersection graphs, we observe that the distribution of cliques is Poisson, hence our formulas align in predicting strictly positive critical infection rate.   

\paragraph{Organisation.}
In the next section we define the contact process, the random intersection graph and its local weak limit, the branching community tree. In Section \ref{sec:results} we present our theorems formally followed by more technical discussions. Then, in Section \ref{sec:lwl-survival} we show survival of the contact process on the branching community trees, relating to part 1 of Meta-theorem \ref{thm:meta}. In Section \ref{sec:lwl-extinction} we show extinction for small $\lambda$ on the branching community tree under the conditions in part 2 of Meta-theorem \ref{thm:meta}. In Section \ref{sec:long-survival-on-rig} we prove that a supercritical process persists at most exponentially long in $n$ on the graph, while in Section \ref{sec:persistence-in-rig} we carry out the embedding of the web of communities to show at least exponentially long persistence in $n$ of supercritical infections. 

\section{Models}
In this section we briefly define the contact process, and the underlying random intersection graph that it shall live on; as well as the Benjamini-Schramm local weak limit of the random intersection graphs. 
\subsection{Contact process}
Let $G=(V,E)$ be a locally finite (undirected) graph, where $V$ denotes the set of vertices and $E$ is the set of edges. We denote by $\root$ a designated vertex contained in $V$, which we call the root. Furthermore, we denote for every $v\in V$ by $\cN_v$ the neighbourhood of vertex $v$, i.e.\ all $u\in V$ such that $\{u,v\}\in E$.

The contact process $\bfxi=(\bfxi_t)_{t\geq 0}$ is a time-continuous Markov process with values in $\{0,1\}^{V}$. Let $\bfxi$ be in state $\xi$, then it has for every $v\in V$ transitions
\begin{equation}\label{eq:Infection-Dynamics}
    \begin{aligned}
        \xi\mapsto \xi -\one_{v} & \text{ with rate } \xi(v) \text{ and }\\
        \xi\mapsto \xi +\one_{v} & \text{ with rate } \lambda (1-\xi(v)) \sum_{u\in \cN_v} \xi(u),
    \end{aligned}
\end{equation}
where $\lambda>0$ is called the \emph{infection rate} and $\one_v$ denotes the vector which contains a $1$ at position $v$, i.e.\ $\one_v(v)=1$, and $0$ at any other position, i.e.\ $\one_v(u)=0$ for all $u\neq v$. A vertex $v$ is called \emph{infected} at time $t$ if $\bfxi_t(v)=1$ and healthy otherwise. We will indicate the initial configuration $\xi$ by adding it as a superscript to the process, i.e. $\bfxi^{\xi}$.

Since $\{0,1\}^V$ is isomorphic to $\cP(V)$ we use the common abuse of notation where $\xi$ and $\bfxi_t$ are also identified as (random) subsets of $V$. Thus, we denote by 
\begin{equation}\label{eq:extinction-time-def}
    T_{\text{ext}}^{\xi}:=\inf\{t>0: \bfxi^{\xi}_t\neq \emptyset\}
\end{equation}
the \textit{extinction time} of the infection started from initial state $\xi$. We introduce 
\begin{equation}\label{eq:survival-probabilities-def}
\begin{aligned}
   \theta(\xi,\lambda;G)&:= 
   \P_{G,\lambda}\big(T_{\text{ext}}^{\xi}=\infty \big),\\
   \phi(\xi,\lambda;G)&:=    \P_{G,\lambda}\big(\forall s \geq 0\, \exists \,t \geq s \, : \, \bfxi^{\xi}_{t}(\root)=1 \big),
\end{aligned}
\end{equation}
the probabilities of \textit{global survival} and \textit{local survival}, respectively, of the contact process on a graph $G$ started from initial state $\xi$. If clear by context we omit the graph $G$ in the notation.

\subsubsection{Graphical representation}\label{sec:graphical-rep}
We will now explicitly construct the contact process on a given graph $G=(V,E)$ via a the usual graphical representation. For that we denote by $\vec{E}$ the set of all oriented edges, i.e.\ $\vec{E}$ contains all $(x,y)\in V^2$ with $\{x,y\}\in E$.

Let us denote by $\cI=(\cI_e)_{e\in \vec{E}}$ a Poisson point process on $\vec{E}\times \IR_+$ with intensity measure $\lambda \mathsf{d} t$. Furthermore, let $\cR=(\cR_x)_{x\in V}$ be a standard Poisson point process on $V\times \IR_+$ independent of $\cI$. Now $\cI$ signifies all potential \textit{infection} events and $\cR$ all potential \textit{recovery} events.

\begin{definition}[Contact process infection path]
Let  $(y, s)$ and $(x, t)$ with $0\leq s < t$ be space-time points. We write $(y, s)\longrightarrow (x, t)$ if there exists a sequence of times $s = t_0 < t_1 < \dots < t_n \leq t_{n+1} = t$ and space points $y = x_0,x_1,\dots, x_n = x$ such that $\big((x_{k},x_{k+1}),t_k\big)\in \cI$ and  $\cR\cap\big(\{x_k\}\times[t_k , t_{k+1} )\big)=\emptyset$ for all $k \in \{ 0, \dots , n\}$.
\end{definition}

Now we define the contact process by setting $\bfxi_0(x):=\xi\in \{0,1\}^V$ and
\begin{align}\label{DefinitionCP}
	\bfxi^{\xi}_t:=\begin{cases}
	    1 &\text{if }\exists y\in V \text{ with }\xi(y)=1 \text{ such that } (y,0)\longrightarrow (x,t) \\
        0 & \text{ otherwise}
	\end{cases}
\end{align}
for every $t\geq 0$. 
The process $\bfxi$ defined in this is a Markov process with transition rates as in \eqref{eq:Infection-Dynamics}. For more details on the graphical representation see for example the book of Swart \cite{swart2026interacting}.

\subsection{Random Intersection Graph}\label{sec:rig-def}

We now introduce the \textit{random intersection graph} where the respective communities form complete graphs. We use a simplified version and simplified notation of this model as done in \cite{van2021random}, where it was introduced in greater generality. 

We assume that there are $n\in \IN$ many \textit{individuals} and $a_n\in \IN$ many \textit{communities}, where we assume that $a_n\to\infty$ as $n\to\infty$. Thus, we denote the set of all individuals by $\leftpar_n = [n]$ and the set of communities by $\rightpar_n = [a_n]$. We shall define the model using a two-step procedure: in the first step, we assign to each individual the communities it will be part of, then in the second step we perform a `community projection' where we obtain a graph only on the set of individuals.

We start with the first step. To assign the community memberships of individuals, we construct a bipartite graph where individuals are placed on one side and the communities on the other. The number of community memberships of an individual $x\in\leftpar$ is called its membership degree and denoted by $m_x=m^{(n)}_x$.
For a community $w\in\rightpar$, we denote the community degree, i.e.\ the number of individuals which can be part of community $w$,  by $k_a=k_a^{(n)}$. We express the membership and community degrees as vectors 
\begin{equation*}
    \bfm=\bfm^{(n)} := \big(m_x^{(n)}\big)_{x\in\leftpar_n}\quad \text{ and  } \quad \bfk=\bfk^{(n)} := \big(k_w^{(n)}\big)_{w\in\rightpar_n}
\end{equation*}

\begin{definition}[Bipartite configuration model, $\BC_n(\bfm,\bfk)$]\label{def:bipartite-conf-model}
    Let the sequences $\bfm=\bfm^{(n)}$ and $\bfk=\bfk^{(n)}$ of non-negative integers satisfy
    \begin{equation}\label{eq:matching-condition}
        H_n := \sum_{x\in\leftpar_n} m_x^{(n)} = \sum_{w\in\rightpar_n} k_w^{(n)}.
    \end{equation}
    We assign to every individual $x\in \leftpar_n$, $m_{x}^{(n)}$ many membership `half-edges’ and to every community $w\in\rightpar_n$, $k_{w}^{(n)}$ many community ‘half-edges’, respectively. We denote the set of all membership half-edges by $\cH^{\mem}$ and the set of all community half-edges by $\cH^{\comm}$.
    
    We define the \emph{bipartite configuration model} $\BC_n(\bfm,\bfk)$ as a distribution on (multi-) graphs constructed as follows:
    We perform a \textit{uniform random bipartite matching} of these half-edges, i.e., we only match membership half-edges with community half-edges. For each matched pair of half-edges, we associate an edge between the respective individual and community  in $\BC_n(\bfm,\bfk)$. Finally, we write $e_{\BC}(x,w)$ for the number of edges between an individual $x\in \leftpar_n$ and a community $w\in \rightpar_n$.
\end{definition}
Note that by assuming that \eqref{eq:matching-condition} is satisfied we ensure that every membership half-edge can be paired with a community half-edge. 

In the second step, we construct the actual random intersection graph $\RIG$ on the set of individuals via a procedure we call \textit{community-projection}.  First we define for every community $w\in \rightpar_n$, the set
\begin{equation*}
    \cK^{(n)}_w:=\{x\in \leftpar_n: e_{\BC}(x,w)\geq 1\}.
\end{equation*}
The idea is that every half-edge of $w$ represents a (possibly fictitious) person in a community and all the members of the community should connect to each other, giving a complete graph of size $k_w^{(n)}$. If each half-edge of $w$ is matched to different individuals $x_1, \dots, x_{k_w^{(n)}}$, then this simply translates to superimposing a complete graph on this set of vertices in the random intersection graph. We think of these edges as coming from community $w$, another community $w'$ may superimpose further edges on the same pair of vertices. However, as the same individual can have more than one half-edge matched to the community $w$, an individual can take multiple of the fictious community roles, and the complete graph may become smaller. We would like to define the procedure to creating the random intersection graph in a way that the total number of edges due to community $w$ is deterministic and equals $\binom{k_w^{(n)}}{2}$. This is done as follows.   
\begin{enumerate}
    \item We add an edge between $x,y\in \leftpar_n$ in the new graph $\RIG$ for every possible half-edge-labelled cherry between $x,y$ via any community $w\in\rightpar_n$, i.e.,
    \begin{equation}\label{eq:edge-rule-rig}
        e_{\RIG}(x,y):=\sum_{w\in \rightpar_n}e_{\BC}(x,w)e_{\BC}(y,w).
    \end{equation}
    \item Furthermore, we add a self-loop $x\in \leftpar_n$ for every possible half-edge-labelled cycle of length $2$ of the form $(x,w,x)$ for every $w\in \rightpar$, i.e.,
    \begin{equation}\label{eq:loop-rule-rig}
        e_{\RIG}(x,x):=\sum_{w\in \rightpar}\binom{e_{\BC}(x,w)}{2}.
    \end{equation}
\end{enumerate}
This `community-projection' procedure results in a multi-graph where the vertices that have a half-edge matched to community $w$ form a (possibly partly contracted) complete graph, and the number of edges associated to community $w$ in RIG is exactly ${k_w^{(n)} \choose 2}$. Note that the membership vectors and the community-size vectors $\bfm := (m_v)_{v\in\leftpar_n}$, $\bfk := (k_w)_{w\in\rightpar_n}$ together with the random matching fully determine the distribution of the resulting (multi)-graph.

\begin{definition}[Random intersection graph, $\mathrm{RIG}_n(\bfm, \bfk)$]\label{def:random-intersection-graph}
    Let $\bfm=\bfm^{(n)} := (m_v)_{v\in\leftpar_n}$, $\bfk=\bfk^{(n)} := (k_w)_{w\in\rightpar_n}$ satisfy \eqref{eq:matching-condition} and let $\BC_{n}(\bfm,\bfk)$ be the corresponding bipartite configuration model as defined in Definition~\ref{def:bipartite-conf-model}. Then the \textit{random intersection graph} $\mathrm{RIG}_n(\bfm, \bfk)$ is defined as the random multi-graph obtained by the community projection of $\BC_{n}(\bfm,\bfk)$. 
\end{definition}
Throughout this section, we make use of the following description. Let $U_n \sim \mathrm{Unif}[\leftpar_n]$ and $A_n \sim \mathrm{Unif}[\rightpar_n]$ denote uniformly chosen individuals and community vertices respectively, and define
\begin{equation*}
    M_n := m_{U_n},
    \qquad K_n := k_{A_n}. 
\end{equation*}
Then we define the probability mass functions as
\begin{equation}\label{eq:empirical-probability-M}
    p^{(n)}_\ell := \frac{|\{x\in\leftpar_n:\, m^{(n)}_x =\ell\}|}{|\leftpar_n|} \quad \text{ and } \quad q^{(n)}_\ell :=  \frac{|\{w\in\rightpar_n:\, k^{(n)}_w = \ell\}|}{|\rightpar_n|}
\end{equation}
for each $\ell\in\Z^+$. We collect the probability mass functions in the (infinite-dimensional) probability vectors $\bitp^{(n)}  = ( p^{(n)} _\ell)_{\ell\in\Z^+}$, $\bitq^{(n)}  = ( q^{(n)} _\ell)_{\ell\in\Z^+}$. 

To be able to connect graphs across different values of $n$, we shall work under the following assumption throughout the paper.

\begin{assumption}
\label{asmp:convergence} Letting~$M_n \sim \mathbf p^{(n)}$ and~$K_n \sim \mathbf q^{(n)}$ for each~$n$, we assume for some $\eta>1$
		\begin{equation*}
			M_n \xrightarrow[\mathrm{d}]{n \to \infty} M,
            \qquad \mathbb E[M_n^{\eta}] \xrightarrow{n \to \infty} \mathbb E[M^{\eta}]\in \R_+,
		\end{equation*}
		and similarly,
		\begin{equation*}
			K_n \xrightarrow[\mathrm{d}]{n \to \infty} K,
            \qquad \mathbb E[K_n^\eta] \xrightarrow{n \to \infty} \mathbb E[K^\eta]\in \R_+.
		\end{equation*}
        for some $\eta>1$.
\end{assumption}
Due to $h_n=\E[M_n] n = \E[K_n] a_n$, note that these assumptions already imply that
\begin{equation*}
    \gamma_n:=\frac{|\rightpar_n|}{|\leftpar_n|} = \frac{\E[M_n]}{\E[K_n]} \to \frac{\E[M]}{\E[K]} =: \gamma \in \R_+ \quad \text{ as } \quad n\to \infty.
\end{equation*}

Under these assumptions, the local weak limit of the random intersection graphs is known, which we define in the next section. Intuitively, the local weak limit describes the distribution of the local graph neighbourhood of a uniformly chosen vertex in the large network limit (as $n$ tends infinity).
For local weak limit specific to graphs constructed using half-edges,  we refer the reader for a compact introduction to \cite[Section 3]{van2021random}.

\subsection{Branching community trees}\label{sec:Definition-community-tree}
In this section we describe the local weak limit of the random intersection graph. Just like with~$\mathrm{RIG}$, we construct it in two steps. First we define a branching process with two different offspring distributions, one applied in even and one in odd generations. This corresponds to the local weak limit of the bipartite configuration model in Definition \ref{def:bipartite-conf-model}. Then we describe the equivalent notion of the `community projection' for this tree, obtaining the branching community tree.

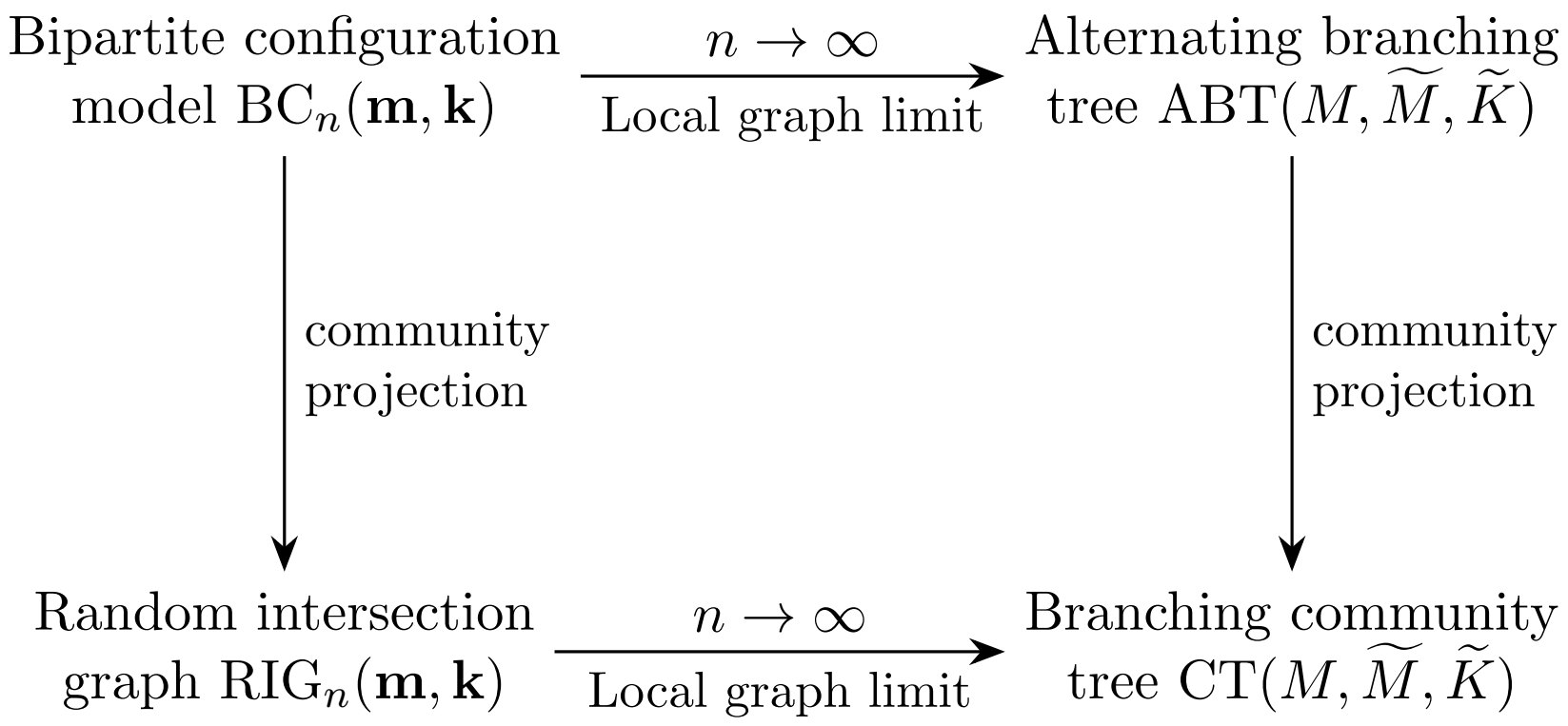
\begin{figure}[ht]
\centering
\begin{tikzpicture}[
  obj/.style={align=center},
  arr/.style={-{Stealth[length=2.5mm]}, semithick},
  lab/.style={font=\small, align=left}
]
  \node[obj] (B)   at (0, 0)  {Bipartite configuration\\ model $\mathrm{BC}_n(\mathbf m,\mathbf k)$};
  \node[obj] (ABT) at (7, 0)  {Alternating branching\\ tree $\mathrm{ABT}(M,\widetilde{M},\widetilde{K})$};
  \node[obj] (RIG) at (0,-4)  {Random intersection\\ graph $\mathrm{RIG}_n(\mathbf m, \mathbf k)$};
  \node[obj] (CT)  at (7,-4)  {Branching community\\ tree  $\mathrm{CT}(M,\widetilde{M},\widetilde{K})$};

  \draw[arr] (B)   -- node[above] {$n\to\infty$}
                      node[below, font=\small] {Local graph limit} (ABT);
  \draw[arr] (RIG) -- node[above] {$n\to\infty$}
                      node[below, font=\small] {Local graph limit} (CT);

  \draw[arr] (B)   -- node[right, lab] {community\\ projection} (RIG);
  \draw[arr] (ABT) -- node[right, lab] {community\\ projection} (CT);
  \end{tikzpicture}
  \caption{Diagram showing the local graph limits as~$n \to \infty$ of the bipartite configuration model~$\mathrm{BC}_n(\mathbf{m},\mathbf{k})$ and the random intersection graph~$\mathrm{RIG}_n(\mathbf{m},\mathbf{k})$. Here,~$M$ and~$K$ are random variables with the limiting distributions of Assumption~\ref{asmp:convergence}, and~$\widetilde{M}$ and~$\widetilde{K}$ are their downshifted size-biased versions, see Definition~\ref{def:size-biased-notation} below.}
\end{figure}

\begin{figure}[!t]
  \centering
  \includegraphics[width=\textwidth]{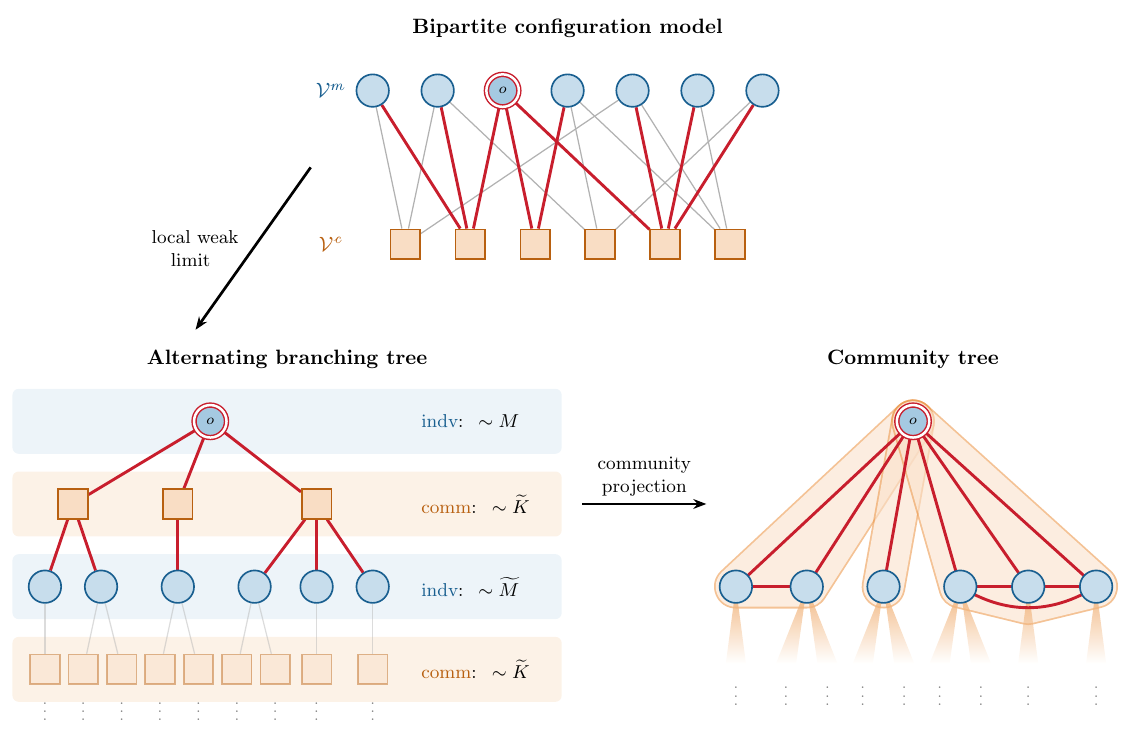}
  \caption{On top, we see the creation of the bipartite configuration model: the top row contains vertices corresponding to individuals and the bottom row contains  vertices corresponding to communities. The breadth first search exploration of the vertex $o$ marked in red is indicated by the red edges. Below on the left we see the corresponding local weak limit, which is an alternating branching tree, and below on the right we see the corresponding branching community tree which results from the community projection.}
  \label{fig:BCLWL}
\end{figure}

\begin{definition}[Alternating branching trees with special root]\label{def_alternating}
    Let $D$, $X$ and $Y$ be three distributions with values in $\IN_0$. Consider the tree where the root $\root$ has degree distributed as $D$, and every other vertex in an even generation has independent offspring distributed as $X$, while every vertex in odd generations has independent offspring distributed as $Y$. We denote  random tree obtained in this way by $\mathrm{ABT}(D, X,Y)$.
\end{definition}
This definition is almost enough to identify the local weak limit (in probability) of the bipartite configuration model under Assumption \ref{asmp:convergence}. We still need to define the size-biased versions and size-biased down-shifted versions of discrete random variables.
 \begin{definition}[Size-biasing]\label{def:size-biased-notation} Let $X$ be a $\IN_0$ valued random variable with $\IE[X]\in(0,\infty)$. Then we denote by $X^{\star}$ the size-biased version and by $\widetilde{X}$ the downshifted size-biased version, with respective mass functions
\begin{equation}\label{eq:size-biasing-def}
    \IP(X^{\star}=k)=\frac{k\IP(X=k)}{\IE[X]}  \quad \text{ and } \quad \IP(\widetilde{X}=k):=\IP(X^{\star}=k-1) \quad \forall\, k\geq 0
\end{equation}
\end{definition}
Recall the limit variables $M$ and $K$ from Assumption \ref{asmp:convergence}.
Then \cite[Theorem 2.14]{van2021random} establishes that under Assumption \ref{asmp:convergence}, with a uniformly chosen root vertex $U_n$ from $\leftpar_n$, for the bipartite configuration model it holds that 
\begin{equation*}
    (\BC_n(\bfm, \bfk), U_n)\xrightarrow{\IP-\mathrm{loc}} (\mathrm{ABT}(M, \widetilde M, \widetilde K),\root) \quad \text{ as } \quad n\to\infty.
\end{equation*}
This is visualised in Figure~\ref{fig:BCLWL}.

From here it is easy to guess the local weak limit of the random intersection graph as well. Namely, the community projection, intuitively, replaces each community vertex by a complete graph on the individuals belonging to that community (see Figure~\ref{fig:BCLWL}). When we carry out this operation on the periodic-offspring branching tree $\mathrm{ABT}(M, \widetilde M, \widetilde K)$, we obtain an object that is not a tree but that still preserves much of a tree-like structure. (In particular, it has tree-width $1$).

For a more precise definition, and to be able to also refer to vertices of the obtained graph, we first define branching community trees using given membership- and community size arrays, then later randomize these choices. 
Let $\cU$ be the Ulam-Harris label set $\cU=\{\root\}\cup \cup_{\ell\in \IN_+} (\IN_+)^\ell$, and let us consider two arrays $\IM:=(m_{x})_{x\in\cU}$ and $\IK:=(k_{x,j})_{x \in \cU, j\in \IN}$ of numbers with values in $\IN_0$. 
One can think of these deterministic arrays as a \emph{given realisation} of the random memberships and community sizes. 

Given $\IM, \IK$, let the set of vertices in generation $0$ be defined as $\cZ_{0}=\{\root\}$ and define the set of vertices in generation $\ell+1$ for every $\ell\geq 0$ as  
\begin{equation*}
\cZ_{\ell+1}=\bigg\{xn: x\in \cZ_{\ell}: n\leq \sum_{j=1}^{m_x}k_{x,j}\bigg\}.
\end{equation*}
For example, the root is a member of $m_\root$ many communities, the $j$th one having $k_{\root,j}$ many additional members, so the the size of the first generation $\mathcal Z_1$ is $\sum_{j=1}^{m_\root} k_{o,j}$. When we use independent random variables $M_\root$ for $m_\root$ and independent arrays of independent variables distributed as $\widetilde M$ in $\IM$ and as $\widetilde K$ in $\IK$, we obtain the vertex set of every even generation of the periodic-offspring branching tree $\mathrm{ABT}(M, \widetilde M, \widetilde K)$.
Next we explain the community projection step. 
We introduce notation for the set of vertices in the $i$-th community of $x\in \cZ_\ell$, for $i\leq m_x$, by
\begin{equation*}
    \cK_i(x):=\bigg\{xn: \sum_{j=1}^{i-1}k_{x,j} < n\leq  \sum_{j=1}^{i}k_{x,j}\bigg\}\cup \{x\}.
\end{equation*}
Note that $|\cK_{i}(x)|=k_{x,i}+1$. Finally, the vertex set consists of the union of all generations, and the edge set consists of superimposing a complete graph on each present community. We obtain the community tree $\CT:=\CT(\IM,\IK)$ associated to the arrays $\IM, \IK$ defined as the graph with the following edge and vertex set.
\begin{equation*}
\cV(\CT):=\bigg\{x\in \cU: x\in \bigcup_{\ell\geq 0} \cZ_{\ell}\bigg\}, \quad \quad
\cE(\CT):=\bigcup_{x \in \mathcal V}\bigcup_{i\le m_x}\big\{\{y,z\}\subset \cV(\CT): y,z\in \cK_i(x)\big\}.
\end{equation*} 

We are ready to define the random community tree. To be able to state the local limit of random intersection graphs smoothly, we specify the distribution of the membership-number of the root separately.
\begin{definition}[Branching community trees (CT)]\label{def:branching-community-tree}
      Let $D, X, Y$ be three random variables with values in $\IN_0$. Let $m_0$ be distributed as  $D$ independently of everything else, and let  $\IM\setminus\{m_0\}$ and $\IK$ be independent arrays of iid random variables with respective underlying distributions $X$ and $Y$. 
    Then we denote by $\CT(D, X,Y):=\CT(\IM, \IK)$ the \emph{branching  community tree}. When the laws of $D$ and $X$ agree, we write $\CT(X,Y)$. Furthermore, we call $\IP(X=\cdot)$ the \emph{membership distribution} and  $\IP(Y=\cdot)$ the \emph{community distribution} of the community tree.
\end{definition}

Finally, we cite the following result from \cite{van2021random},about the local weak limit of random intersection graphs of the previous section. Recall the size-biased down-shifted distribution from Definition \ref{def:size-biased-notation}.
\begin{theorem}[Local weak convergence of the RIG]
    Consider $\RIG_n = \RIG_n(\bfm,\bfk)$ under Assumption~\ref{asmp:convergence}. Then, with $U_n\sim \mathrm{Unif}[\leftpar_n]$ as $n\to \infty$.
    \begin{equation*}
        (\RIG_n,U_n)\xrightarrow{\IP-\mathrm{loc}} (\CT(M, \widetilde M, \widetilde K),\root) \quad \text{ as } \quad n\to\infty.
    \end{equation*}
\end{theorem}

\section{Main Results}\label{sec:results}
\subsection{Contact process on branching community trees: the critical values.}
 To be able to understand the behavior of the contact process on a finite  but large random intersection graph, first we shall understand its behavior on its local weak limit, the branching community tree.
Our first result allows for general membership and community distributions, as well as an arbitrary (but not identically $0$) membership number for the root. As such, we suppress the notation for the root. 
Consider the contact process $\bfxi$ on a  community graph. Then, we can identity the \textit{critical infection rate for global survival}
\begin{equation*}
	\lambda_{1}:=\inf \big\{\lambda\geq 0:\  \IE[\theta(\one_{\{\root\}},\lambda;\CT(X, Y))]>0\big\},
\end{equation*}
and the critical value for \textit{local survival} by 
\begin{equation*}
	\lambda_{2}:=\inf \big\{\lambda\geq 0:\  \IE[\phi(\one_{\{\root\}},\lambda;\CT(X,Y)]>0\big\}.
\end{equation*}
Clearly $\lambda_{1}\le \lambda_{2}$ even more generally, for the contact process on all infinite graphs $G$. Our first main result provides a criterion when the process is supercritical for all $\lambda>0$, i.e., it has a positive chance to survive for arbitrarily small infection rates $\lambda>0$. 

\begin{theorem}\label{thm:Survival-Criterion-Community-Tree}
Let $\offmem, \offcomm$ be two random variables with values in $\IN$.
    Let $\bfxi$ be a contact process on the community tree $\CT(\offmem,\offcomm)$. Suppose that $\mu=\IE[\offmem]\IE[\offcomm]>1$ and either
    \begin{equation}\label{eq:conditions-both}
       \E[\exp(c\offmem)]=\infty \quad \text{ or } \quad \E\big[\exp\big(c\offcomm\log(\offcomm)\big)\big]=\infty
    \end{equation}
    for all $c>0$. Then $\lambda_{2}=\lambda_{1}=0$.
\end{theorem}
In particular, if $K$ follows a Geometric distribution with arbitrary mean then this condition is satisfied and we obtain $\lambda_1=0$.
The proof can be found in Section \ref{sec:lwl-survival}. We can give an almost matching counterpart to this theorem. I.e., the existence of a subcritical phase for small infection rates $\lambda>0$ under `almost' the negation of the conditions in Theorem \ref{thm:Survival-Criterion-Community-Tree}.
\begin{theorem}\label{thm:Existence-of-Subcritical-Phase}
Let $\offmem, \offcomm$ be two random variables with values in $\IN$.
    Let $\bfxi$ be a contact process on the community tree $\CT(\offmem,\offcomm)$. 
    Suppose that there exists  $c>0$ such that 
    \begin{equation}\label{eq:conditions-extinction-possible}
       \E[\exp(c\offmem)]<\infty \quad \text{ and } \quad \E[\exp(\offcomm \log(c \offcomm))]<\infty
    \end{equation}
    then $\lambda_{1}>0$.
\end{theorem}
In particular, if $K$ follows a Poisson distribution with arbitrary mean $\mu>0$, then we have~$\mathbb E[\exp(K\log(cK))]<\infty$ for $c<1/(e\mu)$; provided that for some~$c$ within this range we also have~$\mathbb E[\exp(cM)] < \infty$, we obtain $\lambda_1>0$.
We also mention that the condition~$\mathbb E[\exp(K\log(cK))]< \infty$ is equivalent to the other two conditions
\begin{equation}\label{eq:extinction-K}
\exists c>0: \E[(1+ c K)^K]<\infty, \qquad 
\exists C<\infty \text{ s.t. } \E[e^{\offcomm\log (\offcomm)-C \offcomm}]<\infty.
\end{equation}
Note that this is just slightly stronger than the negation of the criterion $\E[e^{cK\log K}]=\infty$ for all $c>0$ in Theorem \ref{thm:Survival-Criterion-Community-Tree}. Here it matters that the coefficient of $K\log K$ in the exponent is~$1$. The equivalence of the two formulations follows from taking expectation in
\begin{equation}
(cK)^K\le (1+c K)^K=(cK)^K \Big(1+\frac{1}{c K}\Big)^K\le e^{1/c} (cK)^K.
\end{equation}
We mention that a class of distributions that is in the `gap' left open by Theorems \ref{thm:Survival-Criterion-Community-Tree} and~\ref{thm:Existence-of-Subcritical-Phase}, that is, where it is unsettled whether $\lambda_1=0$ or $>0$,  can be constructed as follows. Take~$\alpha\in(0,1)$ and define for $k\in \IN$
\begin{equation}
\P(K=k)=k^{-\alpha k}/Z_\alpha, \qquad Z_\alpha:=\sum_{k=1}^\infty k^{-\alpha k}<\infty.
\end{equation}
Then $\E[e^{c K\log K}]<\infty$ for $c\in(0,\alpha)$, and also $\E[e^{K\log K-C K}]=\infty$ for all $C>0$.

The proof of Theorem \ref{thm:Existence-of-Subcritical-Phase} can be found in Section \ref{sec:lwl-extinction}. 

For $d$-regular infinite trees, \cite{pemantle1992contact, stacey1996existence} shows that $0<\lambda_{1} < \lambda_{2}<\infty$. It is generally not known even for branching process trees (Bienaym\'e-Galton-Watson trees) whether $\lambda_{1}<\lambda_{2}$. 
For the branching community tree, 
we show that the same phenomenon prevails for deterministic membership numbers and community sizes, i.e., when $m_x\equiv m$ and $k_{x,j}\equiv k$, then the two critical values differ for certain values of $(m,k)$.
\begin{example}\label{example:local-global-difference}
Let $\CT(\IM,\IK)$ be a branching community tree with $m_{\root}=m+1$, and all other $m_x\equiv m$ and $k_{x,j}\equiv k$. Then 
\begin{equation}
\lambda_{1}\le \frac{1}{mk-1}, \qquad\qquad  \lambda_{2} \ge \frac{1}{k-1+2\sqrt{mk}}.
\end{equation}
In particular, $\lambda_{1}<\lambda_{2}$ whenever $(m,k)$ satisfies $mk-1>k-1+2\sqrt{mk}$.
\end{example}
The upper bound for $\lambda_1$ is direct consequence of Proposition~\ref{prop:Coupling-with-BGWtree}, which provides a coupling with a contact process on a $mk$-regular tree and then the bound follows by \cite[Theorem 4.1(a)]{liggett1999stochastic}. The lower bound for $\lambda_2$ can be shown analogously to \cite[Theorem 4.1(b)]{liggett1999stochastic} with minor modifications.

\subsection{Fast extinction vs long persistence on random intersection graphs}
In this section, we present our main results concerning the extinction time of a contact process on random intersection graphs. 
On a finite graph, the contact process eventually becomes extinct, however, tthis can take exponentially long in the number of vertices of the underlying graph. We call this phenomenon \emph{long survival}. 
Our next results concern this phenomenon for random intersection graphs. Recall Definition \ref{def:size-biased-notation} for size-biased variables.
When the underlying graph is random, we write $\IP_\lambda(\cdot)$ for the joint law of the contact process and the graph. For a graph $G_n$, let then $T^{G_n}_{\mathrm{ext}}$ denote the extinction time of the contact process $\bfxi=\bfxi^{(n)}$ on $G_n$ started from the all-infected initial state, i.e., from  $\bfxi_0=G_n$. 
\begin{theorem}\label{thm:survival-finite-graph}
Consider the random intersection graph $G_n:=\mathrm{RIG}_n(\bfm^{(n)}, \bfk^{(n)})$ satisfying Assumption \ref{asmp:convergence} with limiting distributions $\offmem$ and $\offcomm$ satisfying $\E[\widetilde M]\E[\widetilde K]>1$, and  assume that either $\E[e^{c\widetilde M}]=\infty$ or $\E[e^{c\widetilde K\log \widetilde K}]=\infty$ (or both) hold for all $c>0$. 
Then for any $\lambda > 0$ there exists a constant $c=c(\lambda)>0$ such that
	\begin{equation*}
			\IP_{\lambda}(T^{G_n}_{\mathrm{ext}} \geq e^{c n}) \to 1.
		\end{equation*}
\end{theorem}
\begin{remark}\normalfont
In Theorem \ref{thm:survival-finite-graph} one could equivalently use the original distributions $M, K$ in the criteria $\E[e^{cM}]=\infty$ $\E[e^{cK\log K}]=\infty$ in place of the size-biased downshifted distributions $\widetilde M$ and $\widetilde K$. We kept the size-biased formulation as it reflects the membership and community-size distributions in the associated local-weak limit branching community tree.
\end{remark}
The counterpart of the previous theorem is the following. It requires only that the model has a local weak limit, i.e., convergence of the first moment of $M_n$ and $K_n$, as in \cite{van2021random}.

\begin{theorem}\label{thm:rigc-extinction}
Consider the random intersection graph $G_n:=\mathrm{RIG}_n(\bfm^{(n)}, \bfk^{(n)})$ satisfying Assumptions~\ref{asmp:convergence} with limiting distributions $\offmem$ and $\offcomm$. Then for any $\lambda > 0$ there exists a constant $C=C(\lambda)>0$ such that
\begin{equation*}
	\IP_{\lambda}(T^{G_n}_{\mathrm{ext}} \geq e^{C n}) \to 0 \quad \text{ as } \quad  n\to \infty.
\end{equation*}
\end{theorem}
We mention that the condition of the $\eta$th moment convergence in Assumption \ref{asmp:convergence} can be relaxed to assuming the convergence of $\E[K_n \log (1+K_n)]$ to $\E[K\log (1+K)]$.
The at most exponential extinction time is the consequence of a more general result on extinction times that we state; we believe it is of interest in its own right.
\begin{theorem}\label{thm:general-exponential-extinction}
Let $(G_n)_{n\ge 1}$ be a sequence of possibly random graphs with $|V(G_n)|=n$. Let $D_n$ be the distribution of the degree of a uniformly chosen vertex from  $G_n$. 
Assume that there exists $\delta\in (0,1)$ such that 
\begin{equation}\label{eq:condition-gen-exp-ext}
\limsup_{n\to\infty} \E[\log(1+\delta D_n)]<\infty
\end{equation}
holds. Then for any $\lambda > 0$ there exists a constant $C=C(\lambda)>0$ such that
\begin{equation*}
	\IP_{\lambda}(T^{G_n}_{\mathrm{ext}} \geq e^{C n}) \to 0 \quad \text{ as } \quad  n\to \infty.
\end{equation*}
\end{theorem}
\begin{remark}\normalfont
    This result is applicable for a broad class of (finite) random graphs, and shows that the contact process persists for at most an exponentially long time in the graph size. Beyond the RIG, the theorem also applies to the configuration model, the Chung Lu and Norros-Reitu graph and the preferential attachment graph, and spatial variants such as the geometric inhomogeneous random graphs, hyperbolic random graphs, and the weight-dependent random connection model under the necessary regularity conditions. We refer the reader for these models \cite{van2017random}.
\end{remark}

\paragraph{Discussion on Theorem \ref{thm:rigc-extinction}} The conditions for Theorem \ref{thm:rigc-extinction} only requite that the limiting distributions $M$ and $K$ have $\eta$-th moment for some $\eta>1$. Exponential extinction time is surprising when $\E[K^2]=\infty$, as in this case the number of edges of the graph is superlinear. Indeed, the number of edges is $\sum_{w\in \rightpar_n} {k_w \choose 2}/n$  tends to infinity whenever $\E[K_n^2]\to \infty$. For instance, for setting $\IP(K=k)\asymp k^{-\tau_c}$ with $\tau_c\in(2,3)$ gives communities in the graph that are of order $n^{1/(\tau_c-2)}$ containing of order $n^{2/(\tau_c-2)}\gg n$ many edges. Moreover, two communities of order $n^x$ and $n^y$ typically share $n^{x+y-1}$ many common members, so these communities are closely tied to each other. Yet, we see that extinction is governed on a strict exponential scale in the number of vertices.  

The proofs can be found in Section \ref{sec:long-survival-on-rig}.

\section{Survival on branching community trees}\label{sec:lwl-survival}
In this section we prove Theorem \ref{thm:Survival-Criterion-Community-Tree} about $\lambda_{c,2}=\lambda_{c,1}=0$ under the assumption that either the moment generation function of $M$ is infinite on~$(0,\infty)$ or the moment generating function of $K\log K$ is infinite on~$(0,\infty)$.  We separate the two cases, and first we show survival for any small $\lambda>0$ if the criterion for the membership distribution $M$ is satisfied.
\subsection{When large-degree nodes are responsible for survival}
Recall Definition \ref{def:branching-community-tree}.
\begin{proposition}\label{prop:Coupling-with-BGWtree}
    Let $\CT=\CT(M,K)$ be a community tree and let $\cT$ be a BGW-tree with offspring distribution $D:=\sum_{i=1}^{M}K_i$, where $(\offcomm_i)_{i\geq 1}$ is a sequence of independent copies of $\offcomm$. There exists a coupling such that 
    \begin{equation*}
        \cV(\cT)=\cV(\CT) \quad \text{ and } \quad \cE(\cT)\subset\cE(\CT).
    \end{equation*}
    This implies in particular that for all $\lambda > 0$,
    \begin{equation}\label{eq:Domination-of-survival-probability}
	       \theta(\xi,\lambda;\CT)\geq \theta(\xi, \lambda; \cT)\quad \text{ and } \quad 
           \phi(\xi,\lambda;\CT)\geq \phi(\xi, \lambda; \cT).
    \end{equation}
\end{proposition}
\begin{proof}
    It suffices to provide the claimed coupling of $\CT$ and $\cT$, since then the second claim is a direct consequence. For that, recall the definition of the community tree $\CT$ from Definition \ref{def:branching-community-tree} in Section~\ref{sec:Definition-community-tree}. We define the random tree $\cT$ by setting
    \begin{equation*}
        \cV(\cT):=\cV(\CT) \quad \text{ and }\quad
        \cE(\cT):=\big\{\{x,y\}\in \cE(\CT): \exists n\in \IN_0 \text{ s.t. } y=xn \big\}\subset \cE(\CT),
    \end{equation*}
    i.e.\ we only keep edges between different generations and drop all within-generation edges. Now, by definition it is clear that $\cT$ is a subgraph of $\CT$. Furthermore, it is easy to see that $\cT$ is a BGW-tree, where the number of offspring of a vertex $x\in \cV$ is given by the random variable $D_x:=\sum_{i=1}^{M_{x}}K_{i,x}$, which concludes the proof.
\end{proof}
The next result is a direct consequence of this fact and \cite[Theorem~1.4]{huang2020contact}.
\begin{theorem}\label{thm:Survival-via-stars}
     Let $(\bfxi_t)_{t\geq0}$ be contact process with infection rate $\lambda>0$ on $\CT(M,K)$ with $\bfxi_{0}(o)=1$. If $\mu:=\E[\offmem]\E[\offcomm]>1$ and  
    \begin{equation}\label{eq:Condition-for-survival-via-stars}
	   \IE[\exp(cM)]=\infty \qquad \text{for all }c > 0,
    \end{equation}
then $\bfxi$ survives locally with positive probability for every $\lambda>0$, i.e., $\lambda_{c,2}=0$.
\end{theorem}
\begin{proof}
    By Proposition~\ref{prop:Coupling-with-BGWtree} we know that we can couple a BGW-tree $\cT$ with offspring distribution $D:=\sum_{i=1}^{M}\offcomm_i$, where $(\offcomm_i)_{i\geq 1}$ is a sequence of independent copies of $\offcomm$, with the community tree $\CT$ such that $\cT\subset \CT$ almost surely. 
    Since $\mu>1$ we know that $K>0$ with positive probability, which implies that for any $c>0$ there exists an $\varepsilon>0$ such that $\IE[e^{cK}]>1+\varepsilon$. Thus, we see that
    \begin{equation*}
        \IE[e^{cD}]=\IE\bigg[\prod_{i=1}^{M} \IE[e^{cK_{i}}]\bigg]\geq \IE[(1+\varepsilon)^{M}]=\infty
    \end{equation*}
    for every $c>0$. Furthermore, $\IE[D]=\IE[\offmem]\IE[\offcomm]>1$. Now, \cite[Theorem~1.4]{huang2020contact} implies that the contact process defined on $\cT$ survives locally for every $\lambda>0$. Since $\cT\subset \CT$ almost surely, the same holds true for the contact process on $\CT$, and thus $\lambda_{c,2}=0$.
\end{proof}
\subsection{When large communities are responsible for survival}
We continue with the necessary lemmas to prove Theorem \ref{thm:Survival-Criterion-Community-Tree}. The previous theorem showed that the blow-up of the moment-generating function of the membership distribution makes survival possible for any $\lambda>0$. On a heuristic level, survival in this case comes from high-degree vertices and the `star-graph' (i.e., the center high-degree vertex together with its neighbors) formed around them, as the proof is a reduction to the embedded branching process tree. 
We move on to the novel part of Theorem \ref{thm:Survival-Criterion-Community-Tree}, when survival for arbitrary $\lambda>0$ is due to the communities.  
We will show that the infection survives for every $\lambda>0$ solely due to the existence of large communities, when the moment generating function of $K\log K$ blows up on $\mathbb R_+$. 
In this case, instead of star-graphs, the local substructures that sustain the infections are complete graphs. We show that sufficiently large complete graphs, once ignited, are able to sustain the infection \emph{much longer} than a star-graph with the same number of vertices. This allows us to search for the next sufficiently large community in a much wider distance-range in the branching community tree. This gained distance results makes the blow-up of the $K\log K$ moment sufficient for survival. 

Our proof strategy is thus as follows: 
\begin{itemize}
\item[1.] We first study the contact process in a single fixed but large community. We identify the minimally sufficient community size $k_\lambda$ so that all three of the following happen with probability close to $1$. 
\begin{itemize}
\item[1a.] The contact process \emph{ignites} the community, i.e., upon reaching it at a single vertex, it also infects more than a quarter of the community with probability close to $1$. We call such a community in such a state \emph{infested}. 
\item[1b.] Then, we show that  the infection \emph{persists}, i.e., the community stays infested for a long time, of order $\exp(\Theta(k_\lambda\log (\lambda k_\lambda))$.
\item[1c.]
We then show that an infested community can \emph{transmit} the infestation to another community of size at least $k_\lambda$ whenever that community is within distance $r_\lambda$. 
\end{itemize}
\item[2.] We then show that a (one-ended) infinite chain of communities of size $k_\lambda$, with neighboring communities at graph distance at most $r_\lambda$ from each other, survives locally with probability close to $1$ once the infection is started from a single vertex.
\item[3.] Based on the community size distribution $K$ and on $\lambda>0$, we identify $k_\lambda$ and $r_\lambda$, and we \emph{embed} the previous `community line graph' structure in the branching community tree. 
\end{itemize}
We start with the details for step 1.  Let us fix some notation. 
Let us write $\compgraph_k$ for the complete graph on $k$ vertices, where we set  $V(\compgraph_k)=\{\root\}\cup\{1, \dots, k-1\}$ and $E(\compgraph_k)=\{\{u,v\}\subset V(\compgraph_k): u\neq v\}$.

The first lemma ignites the community. It provides a lower bound on the probability to generate enough infected individuals in a complete graph $\compgraph_k$ if we start from a single infected vertex.
 \begin{lemma}\label{lem: from 1 to many}
    Let $(\bfxi_t)_{t\geq0}$ be contact process with infection rate $0<\lambda<1$ on a complete graph $\compgraph_k$ with $k\in \IN$ satisfying $\lambda k >2$. Assume that $|\bfxi_0|\geq 1$ and define
    \begin{equation}
        \sigma_{k/2}:=\inf\Big\{t\geq 0: |\bfxi_t|\geq \frac{k}{2}\Big\}\quad \text{ and }\quad \sigma_{\emptyset}:=\inf\{t\geq 0: \bfxi_t= \emptyset\},
    \end{equation}
    then
    \begin{equation*}
        \P\big(\sigma_{k/2}\leq \sigma_{\emptyset}\big)\geq 1-\frac{2}{\lambda k}.
    \end{equation*}
\end{lemma}
\begin{proof}
    For a contact process $(\bfxi_t)_{t\geq0}$ which is defined on a complete graph $\compgraph_k$ it is easy to see that the process $(|\bfxi_t|)_{t\geq 0}$ is again a Markov process. In fact $|\bfxi_t|$ is a continuous-time birth-and-death chain that increases by $1$ with rate $\lambda|\bfxi_t|(k-|\bfxi_t |)$ and decreases by $1$ with rate $|\bfxi_t|$. Thus, it is clear that the associated jump chain is given by a simple random walk $(U_n)_{n\geq 0}$ which gets trapped in $0$ and has transitions 
    \begin{align*}
        U_{n+1}= U_n+1 & \text{ with prob. } \tfrac{\lambda(k-U_n)}{\lambda(k-U_n)+1} \text{ and }\\
        U_{n+1}= U_n-1 & \text{ with prob. } \tfrac{1}{\lambda(k-U_n)+1}.
    \end{align*}
    By assumption, the chain is started from $|\bfxi_{0}|\ge 1$. Thus, as long as $U_n\leq \frac{k}{2}$ we get 
    \begin{equation*}
        \frac{\lambda(k-U_n)}{\lambda(k-U_n)+1}\geq \frac{\lambda (k-\tfrac{k}{2})}{\lambda (k-\tfrac{ k }{2})+1}=\frac{\lambda k}{\lambda k + 2}.
    \end{equation*}
    Now let $W_n$ be a (biased) random walk on $\{0, 1, \dots, \lceil k/2\rceil\}$ with $W_0= 1$ and transitions
    \begin{align*}
        W_{n+1}= W_n+1 & \text{ with prob. } p_W=\tfrac{\lambda k}{\lambda k + 2} \text{ and }\\
        W_{n+1}= W_n-1 & \text{ with prob. } 1-p_W=\tfrac{2}{\lambda k + 2},
    \end{align*}
    and absorbed at $0$ and reflected at $\lceil k/2\rceil$.
    Observe that $W_n$ is stochastically dominated from above by $U_n$. Furthermore, define 
    \begin{equation*}
        \overline{\sigma}_W:=\inf\big\{n\geq 0:W_n\geq  \tfrac{k}{2}\big\} \quad \text{ and } \quad 
        \underline{\sigma}_W:=\inf\big\{n>0:W_n=0\big\}
    \end{equation*}
    Note that $p_{W}\geq 1-p_{W}$, and thus the stochastic domination (using biased Gambler's ruin) gives that
    \begin{equation*}
        \IP(\sigma_{k/2}\leq \sigma_0) \geq 1- \IP(\overline{\sigma}_W\geq \underline{\sigma}_W)=1-\frac{1-p_W}{p_W}= 1-\frac{2}{\lambda k}.
    \end{equation*}
\end{proof}
The next lemma carries out Step 1b in the proof outline and shows that a complete graph, once ignited, stays infested for a long time.
\begin{lemma}\label{lem:Keep the infection alive}
    Let $(\bfxi_t)_{t\geq0}$ be contact process with infection rate $0<\lambda<1$ on a complete graph $\compgraph_k$ with $k\in \IN$ satisfying $\lambda k >640$.
 Let $S_k:=\exp\big(\tfrac{k}{16}\log\big(\lambda k\big)\big)$. Then it holds that
    \begin{equation*}
         \P\big( \inf_{0\leq s\leq S_k} |\bfxi_s|\geq \tfrac{k}{4} \,\big{|}\, |\bfxi_0|\geq \tfrac{k}{2}\big)\geq 1-\tfrac{1}{S_k}.
    \end{equation*}
\end{lemma}
\begin{proof}
    See \cite[Lemma~2.1]{van2018super}
\end{proof}
The next lemma combines Steps 1a and 1b and shows that a complete graph does not go extinct for a long time. 
\begin{lemma}\label{lem:Dying out in a comp graph}
         Let $(\bfxi_t)_{t\geq0}$ be contact process with infection rate $0<\lambda<1$ on a complete graph $\compgraph_k$ with $k\in \IN$ satisfying $\lambda k >640$. Let $S_k:=\exp\big(\tfrac{k}{16}\log\big(\lambda k\big)\big)$. Then it holds that
    \begin{equation*}
         \P\big( \bfxi_{S_k}=\emptyset \,\big{|}\, |\bfxi_0|=1\big)\leq \frac{2}{\lambda k}+\frac{1}{S_k}.
    \end{equation*}
\end{lemma}
\begin{proof}
    We see that
    \begin{equation*}
        \P\big( \bfxi_{S_k}=\emptyset \,\big|\,|\bfxi_0|=1\big)\leq \IP(\sigma_{k/2}>\sigma_0)+\IP(\sigma_{k/2}<\sigma_0,\bfxi_{S_k}=\emptyset).
    \end{equation*}
    By Lemma~\ref{lem: from 1 to many} we can bound the first term from above by $\frac{2}{\lambda k}$. For the second term we can use the strong Markov property together with Lemma~\ref{lem:Keep the infection alive} which implies
    \begin{equation*}
       \IP(\sigma_{k/2}<\sigma_0,\bfxi_{S_k}=\emptyset)\leq \IP\big(\bfxi_{S_k}=\emptyset\,|\,|\bfxi_0|=\tfrac{k}{2})\leq \tfrac{1}{S_k}.
    \end{equation*}
    This yields the claim.
\end{proof}

After obtaining the necessary results for a contact process on a complete graph of size $k$ the next step is to study its behaviour on a line consisting of complete graphs $\compgraph_k$, all a certain distance $r$ apart.
We first define a community line graph, which is a chain of complete graphs separated by paths of length $r$. 
\begin{figure}[t]
    \centering

    \begin{tikzpicture}[
        x=1.15cm,
        y=1.15cm,
        vertex/.style={circle,draw,fill=black,inner sep=1.4pt},
        lab/.style={font=\scriptsize}
    ]

\foreach \i in {0,...,12}{
    \node[vertex] (v\i) at (\i,0) {};
}

\foreach \i in {0,...,11}{
    \pgfmathtruncatemacro{\j}{\i+1}
    \draw (v\i) -- (v\j);
}

\node[lab,below=2pt] at (v0) {$o=v_0$};

\foreach \i in {1,...,12}{
    \node[lab,below=2pt] at (v\i) {$v_{\i}$};
}

\draw (v12) -- ++(0.65,0);
\node[right] at ($(v12)+(0.72,0)$) {$\cdots$};

\node[vertex] (a0) at ($(v0)+(-0.55,0.95)$) {};
\node[vertex] (b0) at ($(v0)+(0.00,1.90)$) {};
\node[vertex] (c0) at ($(v1)+(0.00,1.90)$) {};
\node[vertex] (d0) at ($(v1)+(0.55,0.95)$) {};

\draw
    (v0)--(a0)
    (v0)--(b0)
    (v0)--(c0)
    (v0)--(d0)
    (v1)--(a0)
    (v1)--(b0)
    (v1)--(c0)
    (v1)--(d0)
    (a0)--(b0)
    (a0)--(c0)
    (a0)--(d0)
    (b0)--(c0)
    (b0)--(d0)
    (c0)--(d0);


\node[vertex] (a1) at ($(v5)+(-0.55,0.95)$) {};
\node[vertex] (b1) at ($(v5)+(0.00,1.90)$) {};
\node[vertex] (c1) at ($(v6)+(0.00,1.90)$) {};
\node[vertex] (d1) at ($(v6)+(0.55,0.95)$) {};

\draw
    (v5)--(a1)
    (v5)--(b1)
    (v5)--(c1)
    (v5)--(d1)
    (v6)--(a1)
    (v6)--(b1)
    (v6)--(c1)
    (v6)--(d1)
    (a1)--(b1)
    (a1)--(c1)
    (a1)--(d1)
    (b1)--(c1)
    (b1)--(d1)
    (c1)--(d1);

\node[vertex] (a2) at ($(v10)+(-0.55,0.95)$) {};
\node[vertex] (b2) at ($(v10)+(0.00,1.90)$) {};
\node[vertex] (c2) at ($(v11)+(0.00,1.90)$) {};
\node[vertex] (d2) at ($(v11)+(0.55,0.95)$) {};

\draw
    (v10)--(a2)
    (v10)--(b2)
    (v10)--(c2)
    (v10)--(d2)
    (v11)--(a2)
    (v11)--(b2)
    (v11)--(c2)
    (v11)--(d2)
    (a2)--(b2)
    (a2)--(c2)
    (a2)--(d2)
    (b2)--(c2)
    (b2)--(d2)
    (c2)--(d2);     

    \end{tikzpicture}

    \caption{The community line graph $\compline(4,6)$.}
    \label{fig:community-line}
\end{figure}
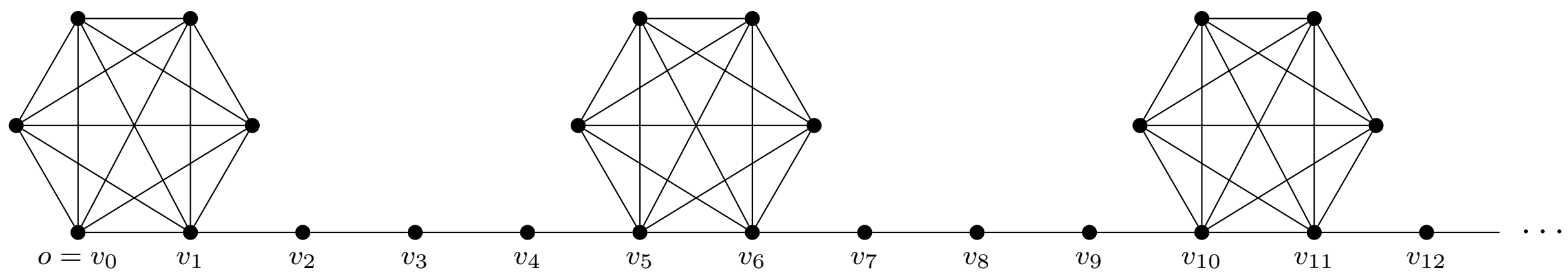
\begin{definition}[Community line graph]\label{def:line-of-complete-graphs}
    Let $k,r\in \IN$. We denote by $\compgraph_{k}^{(i)}=\big(V_k^{(i)},E_k^{(i)}\big)$ a complete graph of size $k$ for every $i\in \IN_0$ and assume that $V^{(i)}_k\cap V^{(i')}_k=\emptyset$ for every $i\neq i'$. Additionally, let $\cL=\{v_k:k\in \IN_0\}$ be a set of distinct vertices with $v_0=o$. We denote by 
\begin{equation}\label{eq:line-segements}
        \cL_i(r):=\{v_{k}: i(r+2)\leq k< (i+1)(r+2)\},
    \end{equation}
    the line segments of vertices that connect two complete graphs, where $i\in \IN_0$. We further assume that $V^{(i)}_k\cap\cL_{i}(r)=\{v_{i(r+2)}\}$  and $V^{(i+1)}_k\cap\cL_{i}=\{v_{(i+1)(r+2)}\}$ for every $i\in \IN_0$. Then, we define
    \begin{align*}
        \cV_{r,k}&:=\bigcup_{i\in \IN_0} (V_k^{(i)}\cup \cL_i(r)),\\
        \cE_{r,k}&:=\bigcup_{i\in \IN_0} E^{(i)}_k\cup \big\{\{v_i,v_{i+1}\}: i\in \IN_0 \big\}.
    \end{align*}
    Finally, we set $\compline(r,k):=(\cV_{r,k},\cE_{r,k})$, which is a line of complete graphs that are all of size $k$ and $r$ is the distance between each of them. Furthermore, we define 
\begin{equation}\label{eq:block-segments}
         \cB^{(i)}(r,k):= V_k^{(i)}\cup \cL_i(r)\cup V_k^{(i+1)}
    \end{equation}
    for every $i\geq 0$ and call them block segments.
\end{definition}
The next result gives a lower bound on the probability to transmit the infection between two vertices which are a distance $r$ apart in a graph.
\begin{lemma}\label{lem:probapath}
     Let $\bfxi$ be a contact process with infection rate $\lambda\in(0,1)$ on some arbitrary graph $G=(V,E)$ with $o\in V$ and $\bfxi_0=\one_{o}$. We assume that there exist a vertex $y\in V$ with $d(y,o)=r\geq 1$. Then, it holds that 
	\begin{equation*}
	        \P(\exists s\in[0,2r]: \bfxi_{s}(y)=1)\geq\frac{1}{2}\bigg(\frac{\lambda}{\lambda+1}\bigg)^{r}.
	\end{equation*}
\end{lemma}
\begin{proof}
    Lemma \cite[Lemma 3.1]{huang2020contact} yields that
    \begin{equation*}
        P(\exists s\in[0,2r]: \bfxi_{s}(y)=1)\geq(1-e^{-\gamma r})\bigg(\frac{\lambda}{\lambda+1}\bigg)^{r},
    \end{equation*}
    where $\gamma=I(2)>1$ and $I$ is the rate function of an exponential distribution. We see that $(1-e^{-1})>\frac{1}{2}$. Thus, the claim follows. A more detailed proof for this type of result can be found in \cite[Lemma 5.1]{cardona2021fitness}.
\end{proof}
The next result is exactly Step 1c, and is crucial to show survival. It states that if one community contains sufficiently many infected individuals, then with a high probability the infection will reach the next community.
Recall that $\cB^{(i)}(r,k)$ is the union of two complete graphs of size $k$, connected by a path of length $r$. 
\begin{lemma}\label{lem:transmit-whp}
    Let $(\bfxi_t)_{t\geq 0}$ be a contact process with infection rate $\lambda\in (0,1)$ on $\cB^{(i)}(r,k)$, with  $r,k\in \IN$ satisfying $\lambda k>640$. 
    Assume that $|\bfxi_0(u)\cap \cK^{(i)}_k|\geq 1$. Let $S_k=\exp\big(\tfrac{k}{16}\log\big(\lambda k\big)\big)$. Then
\begin{equation}\label{eq:error-with-rm}
        \begin{aligned}
             \P\Big(\exists s\leq S_k: \, |\bfxi_s(u)\cap \cK^{(i+1)}_k|\geq 1\Big)\geq 1-\bigg(1- \frac{1}{2}\bigg(\frac{\lambda}{2}\bigg)^{r+1}\bigg)^{ \big\lfloor\frac{S_k}{2(r+1)}\big\rfloor}-R(\lambda,k)
        \end{aligned}
    \end{equation}
    where $R(\lambda,k):=\frac{2}{\lambda k}+\frac{1}{S_k}$.
\end{lemma}
\begin{proof}
    We define the events
    \begin{equation*}
        A:=\Big\{\inf_{0\leq t\leq S_k} |\bfxi_t\cap \compgraph^{(i)}_k|\geq 1 \Big\} \quad \text{ and }\quad 
        B:=\Big\{\exists s\leq S_k: |\bfxi_s\cap \compgraph^{(i+1)}_k|\geq 1\Big\}.
    \end{equation*}
    Appealing to Lemma~\ref{lem:Dying out in a comp graph} we see that $\IP(A^{\complement})\leq  R(\lambda,k)$. We use this to conclude that
    \begin{equation*}
         \P(B^{\complement})\leq \IP(B^{\complement}\cap A)+\IP(A^{\complement})\leq \IP(B^{\complement}\cap A)+R(\lambda,k).
    \end{equation*}
    To verify the claim it remains to show that
    \begin{equation*}
        \IP(B^{\complement}\cap A)\leq \bigg(1- \frac{1}{2}\bigg(\frac{\lambda}{\lambda+1}\bigg)^{r+1}\bigg)^{ \big\lfloor\frac{S_k}{2(r+1)}\big\rfloor}.
    \end{equation*}
    We decompose time into blocks of length $r+1$ and define for every $j\in \IN$,
    \begin{equation*}
        B_{j}:=\{\exists s\in [2(r+1)(j-1),2(r+1) j] : \,\, |\bfxi_s(u)\cap \cK^{(i+1)}_k|\geq 1\}.
    \end{equation*}
    We define $L:=\big\lfloor\tfrac{S_k}{2r}\big\rfloor$ and we see that $B\supseteq \bigcup_{j=1}^{L} B_{j}$. Furthermore, we set for every $k\in \IN_0$ 
    \begin{equation*}
        A_k:=\{|\bfxi_{(2r+1)k}\cap \compgraph^{(i)}_m|\geq 1\}.
    \end{equation*}
    Using the strong Markov-property and Lemma~\ref{lem:probapath} on the path of length $r+1$ between any infected vertex in $\compgraph_k^{(i)}$ and the other community $\compgraph_k^{(i+1)}$ yields
    \begin{equation*}
        \IP(B_k\mid A_k)\geq \frac{1}{2}\bigg(\frac{\lambda}{\lambda+1}\bigg)^{r+1}.
    \end{equation*}
    Finally, we observe that $A\subseteq A_j$ for all $j\leq L$ yields that $\IP(B_j^{\complement } \cap A)\le \IP(B_j^{\complement } \cap A_j) \le \IP(B_j^{\complement } \mid A_j)$. Observe now that that the events $B_j\mid A_j$ can be stochastically dominated by independent events where we set an arbitrary single vertex in $\cK_k^{(i)}$ infected at time $2(r+1)(j-1)$, and so we may apply Lemma \ref{lem:probapath} giving
        \begin{equation*}
        \IP(B^{\complement}\cap A)\leq \prod_{j=1}^{L} \IP(B_j^{\complement}\cap A) \leq \bigg(1- \frac{1}{2}\bigg(\frac{\lambda}{\lambda+1}\bigg)^{r+1}\bigg)^{ \big\lfloor\tfrac{S_j}{2(r+1)}\big\rfloor}.
    \end{equation*}
   Since $\lambda<1$, this yields the claim.
\end{proof}
The next lemma describes the relation between the community size $k$ and the distance $r$ between communities in the line graph so that transmission is possible with a given $\lambda\in(0,1)$.   
\begin{lemma}\label{lem:passing-the-infection-on}
For every~$\varepsilon  > 0$ there exists~$L > 640$ such that the following holds. Let~$\lambda \in (0,1)$ and~$k,r \in \mathbb N$ satisfying
\begin{equation}\label{eq:Minimal-degree-of-a-star}
k \ge \frac{L}{\lambda} \quad\text{ and }\quad \log\bigg(\frac{2}{\lambda}\bigg)^{-1}\frac{k}{48}\log(\lambda k)\geq  r.
\end{equation}
Let~$(\bfxi_t)_{t \ge 0}$ be a contact process with infection rate $\lambda\in (0,1)$ on $\cB^{(i)}(k,r)$ for some $i\geq 0$. Assume that $|\bfxi_0\cap \compgraph^{(i)}_k|\geq 1$. Let $S_k=\exp\big(\tfrac{k}{16}\log\big(\lambda k\big)\big)$.
Then, we have
\begin{equation}\label{eq:Passing-on-infection-with-prob}
            \P\big(\exists s\leq S_k: |\bfxi_s(u)\cap \compgraph^{(i+1)}_k|\ge 1\big)
            \geq 1-\varepsilon.
        \end{equation}
\end{lemma}
\begin{proof}
    Let $\varepsilon>0$ be given. As $L>640$, we know that $\lambda k>640$ for every $\lambda\in (0,1)$.
    Thus, can we can apply Lemma~\ref{lem:transmit-whp} to obtain that 
    \begin{equation}\label{eq:transmit-whp}
        \begin{aligned}
             \P\Big(\exists s\leq S_k: \, |\bfxi_s(u)\cap \cK^{(i+1)}_k|\geq 1\Big)\geq 1-\bigg(1- \frac{1}{2}\bigg(\frac{\lambda}{2}\bigg)^{r+1}\bigg)^{ \big\lfloor\frac{S_k}{2(r+1)}\big\rfloor}-R(\lambda,k)
        \end{aligned}
    \end{equation}
    First note that $R(\lambda,k)=2/(k\lambda) + 1/S_k\leq \frac{3}{L}$   for all $k\geq k_{\lambda}(L)$. Therefore, there exists an $L>640$ such that $R(\lambda,k)<\frac{\varepsilon}{2}$ for all $k\geq k_{\lambda}(L)$ with $L\geq L_0$. Note that the choice of $L_0$ does not depend on $\lambda$. It remains to show that there exists a  $L_1>L_0$ such that,
\begin{equation}\label{eq:Desired-Convergences}
         \bigg(1- \frac{1}{2}\bigg(\frac{\lambda}{2}\bigg)^{r+1}\bigg)^{\big\lfloor \frac{S_k}{2(r+1)}\big\rfloor}< \frac{\varepsilon}{2}.
    \end{equation}
    for all $k>k_{\lambda}(L)$ with $L>L_1$.  This then implies \eqref{eq:Passing-on-infection-with-prob} as an immediate consequence.
    
    Using $1-x\geq e^{-x}$, $\lfloor x\rfloor>\tfrac{x}{2}$ for $x\geq 2$ and $\lambda<1$ we see that
    \begin{equation*}
        \bigg(1- \frac{1}{2}\bigg(\frac{\lambda}{2}\bigg)^{r+1}\bigg)^{ \big\lfloor\frac{S_k}{2(r+1)}\big\rfloor}\leq  \exp\bigg(-\frac{1}{2}\bigg(\frac{\lambda}{2}\bigg)^{r+1}\frac{S_k}{4(r+1)}\bigg).
    \end{equation*}
    Now, we consider the exponent of the right-hand side. To be precise, we apply the function $\log(\cdot)$ to this exponent and see that
    \begin{equation*}
        \log\big(\big(\tfrac{\lambda}{2}\big)^{r+1}\tfrac{S_k}{8(r+1)}\big)= \tfrac{k}{16}\log(\lambda k)-(r+1)\log\big(\tfrac{2}{\lambda}\big)-\log(r+1)-\log(8).
    \end{equation*}
    Next, we use that $2r\geq r+1$ and \eqref{eq:Minimal-degree-of-a-star}, which yields
    \begin{equation*}
         \log\big(\big(\tfrac{\lambda}{2}\big)^{r+1}\tfrac{S_k}{8(r+1)}\big)\geq  \tfrac{k}{16}\log(L)-\log(r+1)-\log(8)
    \end{equation*}
    for all $k$ that satisfy \eqref{eq:Minimal-degree-of-a-star}. We see that the right hand side converges to $\infty$ as $L\to \infty$. Because, the functions $\log(\cdot)$ and $\exp(-\cdot)$ are monotone, this implies that there exists a $L_1>L_0$ such that \eqref{eq:Desired-Convergences} holds. Note that the right hand side does not depend on $\lambda$ anymore, and thus the choices of $L_1$ does not depend on $\lambda$.
\end{proof}
Now that we have shown that the a block of two communities transmits the infection to each other with probability at least $1-\varepsilon$, we are ready to show the local survival of the community line graph.
\begin{proposition}\label{thm:Survival-Line-Of-Communities}
     For every $\delta>0$ there exists an $L_\delta>640$ large enough such that the following holds. Let $r\in \IN$ and $(\bfxi_t)_{t\geq 0}$ be a contact process with infection rate $\lambda\in (0,1)$  on $\compline(r,k)$, where $k, r\in \IN$ satisfies \eqref{eq:Minimal-degree-of-a-star} with $L=L_\delta$. Then,
     \begin{equation*}
        \phi\big(\one_{\{\root\}},\lambda; \compline(r,k)\big)>1-\delta.
     \end{equation*}
\end{proposition}
\begin{proof}
    The proof of this result is an adaptation from the proof of \cite[Proposition~4.7]{valesin2024contact}. The core idea is to couple the contact process with a supercritical $1$-dimensional oriented percolation. We only provide a sketch and highlight the differences.
    
    Let $\lambda\in(0,1)$ and recall that $S_k= \exp\big(\tfrac{k}{16}\log\big(\lambda k\big)\big)$,
    and let $k$ satisfy \eqref{eq:Minimal-degree-of-a-star} for the given $\lambda$ and $L_\delta$.
    Furthermore, recall from  Definition~\ref{def:line-of-complete-graphs} the block segments graphs $\cB^{(i)}(r,k))$.
    
    Fix an arbitrary $\varepsilon\in (0,1)$, and set $t_n:=nS_{k}$. We introduce the auxiliary process $\bfzeta^{(i,n)}= (\bfzeta_t^{(i,n)})_{t\in[0,S_k]}$ which is defined on $\{0,1\}^{\cB^{(i)}(r,k)}$ by setting 
    \begin{equation}\label{eq:CP-on-block-segments}
        \bfzeta_{0}^{(i,n)}(x):=\bfxi_{t_n}(x) \qquad \text{ for every }\quad  x\in \cB^{(i)}(r,k).
    \end{equation}
    In words, $\bfzeta$ shall be the time-shifted restriction of the contact process to the block graph $\cB^{(i)}(r,k)$: the process $\bfzeta$ evolves according to the same graphical representation that we use for the contact process $\bfxi$, just restricted to the block segment $\cB^{(i)}(r,k)$. Thus, the two processes are coupled such that $\bfzeta_{t-t_n}^{(i,n)}(x)\leq \bfxi_{t}(x)$ for all $x\in \cB^{(i)}(r,k)$ and $t\in[t_n,t_{n+1}]$.
    Lemma~\ref{lem:Dying out in a comp graph} and Lemma~\ref{lem:passing-the-infection-on} together with the strong Markov property we see that there exists $L$ large enough such that for any $k\geq k_{\lambda}(L)$ holds that if $|\bfzeta_{0}^{(i,n)}\cap \compgraph^{(i)}_k|\geq 1$, then the restricted process reaches the next community within time $S_k$ with probability at least $1-\varepsilon$, that is, we have
\begin{equation}\label{eq:Lower-bound-for-X}
        \P\big(|\bfzeta^{(i,n)}_{S_k} \cap \compgraph^{(i+1)}_k\}|\geq 1\big)\geq 1-\varepsilon.
    \end{equation}
  To define the oriented percolation we introduce the following oriented graph $(\Lambda_0,\vec{E}_0)$, where
    \begin{equation*}
        \begin{aligned}
            \Lambda_0&:=\{(j,n)\in \IN_0^2: j+n\in 2\IN_0\}\\
         \vec{E}_0&:=\big\{\langle (j,n), (j',n+1)\rangle: (j,n)\in \Lambda_0, |j'-j|=1\big\},
        \end{aligned}
    \end{equation*}
    where we use the notation $\langle u, v\rangle$ to indicate a directed edge from $u$ to $v$.
    Now we couple the contact process on $\compline(r,k)$ with an $1$-dependent oriented percolation $\bfX=(\bfX(\vec{e}))_{\vec{e}\in \vec{E}_0}$. Let $\vec{e}=\langle (i,n), (j,n+1)\rangle$ with $|i-j|=1$.
     \begin{itemize}
        \item[(i)] If $|\bfzeta_{0}^{(i,n)} \cap \compgraph^{(i)}_k|\geq 1$, then we set
            \begin{equation*}
                \bfX(\vec{e})=\begin{cases}
                    1 &\text{ if } |\bfzeta^{(i,n)}_{S_k} \cap \compgraph^{(j)}_k|>1\\
                    0 & \text{ otherwise.}
                \end{cases}
            \end{equation*}
        \item[(ii)] On the other hand, if $|\bfzeta_{0}^{(i,n)} \cap \compgraph^{(i)}_k|=0$, then we set $\bfX(\vec{e})=1$.
    \end{itemize}
    Now let $\cG_n$ be the $\sigma$-algebra generated by the graphical representation up to time $t_n$. Then we know by \eqref{eq:Lower-bound-for-X} and using symmetry of the process that conditional on $\cG_n$ the random variables
    \begin{equation*}
        \{\bfX(\langle i,n\rangle, \langle j,n+1\rangle): i,j\in \IN_0,|i-j|=1\}
    \end{equation*}
    are a family of $1$-dependent Bernoulli random variables with success probability lower bounded by $1-\varepsilon$. 
    Furthermore, by definition of the process $\bfX$, for every $n\geq 0$ it follows that
    \begin{equation*}
        \big\{|\bfxi_{t_n} \cap \compgraph^{(0)}_k|\geq 1\}\supset \{\exists \text{ an open path from } (0,0) \text{ to } (0,n) \text{ w.r.t. } \bfX\}.
    \end{equation*}
    To finish the proof we can follow the proof of \cite[Proposition~4.7]{valesin2024contact}, which we sketch only.  
First one shows (e.g. using standard reduction of high-intensity $1$-dependent percolation to independent percolation, see Liggett-Schonman-Stacey \cite{liggett1997domination} that $\bfX$ dominates a independent oriented percolation $\widehat{\bfX}$ with sufficiently high success probability $\IP(\widehat{\bfX}(\vec{e})=1)$. This establishes a coupling with an independent oriented percolation. Since we can choose $1-\varepsilon$ to be arbitrary close to $1$, the corresponding independent oriented percolation $\widehat{\bfX}$ will also be locally supercritical and thus revisit any site infinitely often. With this one can show that contact process does survive locally. 
\end{proof}
The next theorem carries out Step 3, and embeds the community line graph into the community branching tree. It gives the community-driven survival of Theorem \ref{thm:Survival-Criterion-Community-Tree}.
\begin{theorem}\label{thm:Survival-via-Communities}
Let $M, K$ be two random variables with values in $\IN$.
    Let $(\bfxi_t)_{t\geq0}$ be contact process with infection rate $\lambda>0$ on $\CT(M,K)$ and $\bfxi_{0}(o)=1$. Suppose that $\mu=\IE[M]\IE[K]>1$ and
    \begin{equation}\label{eq:condition-survival-in-line-of-communities}
	   \IE[\exp(cK\log(K))]=\infty
    \end{equation}
    for all $c> 0$. Then $\bfxi$ survives locally with positive probability for every $\lambda>0$, i.e., $\lambda_2=0$.
\end{theorem}
\begin{proof}
With the local survival of the contact process on  the community line graph $\compline(r,k)$ at hand in Proposition \ref{thm:Survival-Line-Of-Communities}, our goal is to embed such a graph into the branching community tree $\mathrm{CT}(M,K)$ with positive probability. The general arc of the proof is similar to that of \cite[Theorem~4.1]{valesin2024contact}. As details differ, we provide the proof.  In particular, our goal is to show that for appropriately chosen $k,r\in \IN$, conditional on the the $\root$ having a community of size $k$, with positive probability, $\CT(M,K)$ contains a copy of the graph $\compline(r,k)$ where the first community $\cK^{(0)}$ community contains the root $\root$.

   The assumption \eqref{eq:condition-survival-in-line-of-communities} gives that for every $a,b>0$,
    \begin{equation}\label{eq:Uniform-lower-bound-Comm-Dist}
        \P(\offcomm = k)\geq \exp(- b k \log(a k))
    \end{equation}
    for infinitely many~$k \in \mathbb N$.
    Let $\lambda\in (0,1)$ and for a given $\delta>0$ let $L=L_{\delta}$ be chosen as in Proposition~\ref{thm:Survival-Line-Of-Communities}. Now, we set $b_{\lambda}:=\log\big(\frac{2}{\lambda}\big)^{-1}$ and choose $k\in \IN$ such that 
    \begin{equation}\label{eq:Parameter-choice-1}
        r_k:=\lfloor b_{\lambda} k\log(\lambda k)\rfloor\geq 1
    \end{equation}
    and as $\mu>1$, $M$ is not identically $0$. So we may choose $k$ that satisfies \begin{equation}\label{eq:Parameter-choice-2}
        \E[M]\P(\offcomm = k)\geq \exp\Big(- \tfrac{1 }{4}\log(\mu) b_{\lambda} k \log( \lambda k)\Big).
    \end{equation}
    In case the root has a community of size $k$, say $\cK_1(o),$  we
    choose one arbitrary other vertex (offspring) $x$ from this community and set $\mathrm{CT}(1)$ to be the subgraph consisting of $\cK_1(o)$, and all further offspring of $x$, inheriting generation numbers from the original graph. Set 
    \begin{equation*}
	N(r,k):=| \{ v \text{ in generation } r+1 \text { of } \mathrm{CT}(1)\text{ with a non-parent community of size}\ k \}|.
    \end{equation*}
    By non-parent community we mean communities that only contain  the vertex $v$ and other vertices in generation number one more than that of $v$.
    As the membership numbers and community sizes of vertices in generation $r+1$ are independent of the size of the generation, the expected number of vertices in generation $r+1$ having a size $k$ community is lower bounded by 
    \begin{equation*}
    \begin{aligned}
	       \E[N(r,k) \mid |\compgraph_{1}(\root)|=k, \offmem_{\root}\geq 1] &\geq k \mu^{r}\E[M]\P(K = k)\\
           &\geq k \mu^{r}\exp\big(- \tfrac{1}{4}\log(\mu)b_{\lambda} k \log(\lambda k)\big).
           \end{aligned}
    \end{equation*}
    The factor $\E[M]\P(K=k)$ corresponds to the expectation of the indicator that a vertex has a community of size exactly $k$. 
    Now setting $r=r_k$ in \eqref{eq:Parameter-choice-1} and \eqref{eq:Parameter-choice-2} yields that
    \begin{equation} 
    \begin{aligned}
           \E[N(r,k) \mid |\compgraph_{1}(\root)|=k, \offmem_{\root}\geq 1] &\geq k \exp\Big( r_k\log\mu - \tfrac14 \log(\mu)b_{\lambda}k\log(\lambda k)\Big)\\
           &\ge 
              k\exp\Big( \frac{1}{2}\log(\mu)b_{\lambda}k\log(\lambda k)\Big)>1,
\end{aligned}\label{eq:Auxiliary-Branching-Surives}
    \end{equation}
  since $\mu>1$ and $k\geq L/\lambda$, where we used that $\lfloor x\rfloor\geq \frac{x}{2}$ for $x\geq 1$. Let us now define an auxiliary branching process: this process only counts vertices in every $r+1$st generation that have a non-parent community of size $k$. By the above computation, this branching process is supercritical.

    Thus, we have shown that for any $\lambda>0$, based on $\lambda$ and the mass function $\P(K=k)$,we can find a pair  $k$ and  $r=r_k$ such that if that conditional on $\root$ having a community of size $m$, the community tree $\CT$ contains a copy of the graph $\compline(r,k)$ where one community contains the root $\root$ with positive probability. By the choice above, $(k,r_k)$ satisfies \eqref{eq:Minimal-degree-of-a-star}, and thus Proposition \ref{thm:Survival-Line-Of-Communities} is applicable. 
    Applying now  Proposition~\ref{thm:Survival-Line-Of-Communities} gives that for every $\varepsilon>0$ there exists an $L>0$ large enough that the contact process survives locally on the graph $\compline(k,r_k)$ with probability $(1-\varepsilon)$ for all $\lambda\in(0,1)$. As this is embedded in the original graph as a subgraph,  by monotonicity, the contact process survies on the whole graph also with positive probability.
    
    The event that the root $\root$ has a community of size $k$ is lower bounded by $\P(M\ge 1)\P(K=k)$, which is strictly positive.  It follows that $\phi(\1_{\{o\}},\lambda,\CT)>0$ for all $\lambda>0$.
\end{proof}

Finally, we can conclude our main result.
\begin{proof}[Proof of Theorem~\ref{thm:Survival-Criterion-Community-Tree}]
    This is a direct consequence of Theorems~\ref{thm:Survival-via-stars} and \ref{thm:Survival-via-Communities}.
\end{proof}
\section{Existence of a subcritical phase on branching community trees}\label{sec:lwl-extinction}
In this section we prove Theorem~\ref{thm:Existence-of-Subcritical-Phase}, the existence of a subcritical phase for the contact process on a branching community tree $\CT(M,K)$ under the conditions that the moment generating functions of $M$ and $K^2$ are both finite in an open interval around $0$, see condition \eqref{eq:conditions-extinction-possible}. To show this, we stochastically dominate the contact process by an inhomogeneous-rate reinforced contact process on a tree, and study extinction of the latter. A reinforced contact process is a modification of the contact process, where vertices are only allowed to heal at rate $1$ when their whole subtree is healthy. 

This process satisfies nice recursive properties and allows for analysis of more general trees with additional edges within the children of a given vertex, that we call augmented trees.  
\subsection{Augmented and (deterministic) community trees}
Recall that we denote by $\cU$ the Ulam-Harris label set
\begin{equation*}
    \cU:=\{x=\root x_1\cdots x_n : n\in \IN,\ x_i\in \IN, \ \forall \  i\leq n\}.
\end{equation*}
Let $T=(V,E_1)$ be a tree with $V\subset \cU$ and recall some typical notions for vertices in a tree. We call $\root$ the \emph{root}, to be precise we consider a rooted tree $(T,\root)$. A vertex $x\in V$ with $d(\root,x)=k$ is in the $k$-th generation $k$ and vertex $y\in V$ such that $\{x,y\}\in E_1$ and $d(\root,x)< d(\root,y)$ is a \emph{child} of $x$ and $x$ is the parent of $y$. Furthermore, $z\in V$ is a \emph{descendant} of $x$ if $d(o,x)< d(o,z)$ and a sequence $v_0,\dots,v_n$, with $n\geq 2$ such that $v_0=x$, $v_n=z$, $\{v_{i},v_{i+1}\}\in E_1$ and $d(\root,v_{i})<d(\root,v_{i+1})$ for all $i\leq n-1$ exists.
\begin{definition}[Augmented tree]
	An (rooted) \emph{augmented tree} is a pair~$(\at,o)$, where~$\at$ is a graph of the form $\at=(V,E_1 \cup E_2)$,  where $E_2$ is a set of edges between vertices that have the same parent in~$(V,E_1)$.
\end{definition}
 Note that any realisation of a branching community tree $\mathrm{CT}(M,K)$ is a specific augmented tree, where edges are added between vertices that have a same parent and \emph{belong to the same community}.
\begin{remark}\normalfont
    For augmented trees as well as branching community trees, we will use the tree terminologies parent, child and descendant that is implied by the embedded tree $(V,E_1)$ rooted at~$o$.
\end{remark}
\subsection{Extinction via reinforced contact process}
For an augmented tree, we will use tree terminology (`parent', `child', `descendant' etc.), always in reference to the tree~$(V,E_1)$ rooted at~$o$.

Given an augmented tree~$((V,E_1\cup E_2),o)$ and a vertex~$v \in V$, we let~$G_{v}$ be the subgraph induced by the set of vertices~$\{v\} \cup \{u \in V: u \text{ is a descendant of }v\}$. Rooting this at~$v$, we obtain a new augmented tree.

The \emph{reinforced contact process} on an augmented tree is a contact process with the dynamics modified so that vertices are only allowed to recover from the infection when all their descendants are healthy.

\begin{proposition}\label{prop_augmented}
Let~$(G,o)$ be an augmented tree and~$\lambda > 0$. Let~$\tau$ be the extinction time of the reinforced contact process with rate~$\lambda$ on~$(G,o)$, started from only~$o$ infected, and for each~$v \sim o$, let~$\tau_v$ be the extinction time of the reinforced contact process with rate~$\lambda$ on~$G_v$, started from only~$v$ infected. Then,
\begin{equation*}
\mathbb E[\tau] \le \prod_{v: v \sim o} \left(1 + \lambda \mathrm{comm}_*(v)\cdot \mathbb E \left[ \tau_v \right]\right),
\end{equation*}
where~$\mathrm{comm}_*(v):=|\{p(v) \} \cup \{u \in V: u \sim v \text{ and }u \sim p(v)\}|$, where $p(v)$ is the ancestor of $v$ in the branching community tree.
\end{proposition} Observe that $\mathrm{comm}_\star(v)=k-1$ if $k$ is the size of the community that $v$ and its parent shares.
\begin{proof}
Let~$\tau'$ be the extinction time of a reinforced contact process on~$G$, started from only the root infected, and with modified dynamics  so that:
\begin{itemize}
\item the root infects each neighbour~$v$ with rate~$\lambda \mathrm{comm}_*(v)$, rather than~$\lambda$;
\item neighbours of the root cannot infect each other.
\end{itemize}
Infections upwards in the tree are irrelevant as vertices with an infected descendant are infected by definition. 
It is straightforward to check, using the graphical construction and a coupling, that~$\mathbb E[\tau] \le \mathbb E[\tau']$. Next, we claim that
\begin{equation*}
\mathbb E[\tau'] \le \prod_{v: v \sim o} \left(1 + \lambda \mathrm{comm}_*(v)\cdot \mathbb E \left[ \tau_v \right]\right).
\end{equation*}
This is proved in the same way as~\cite[Lemma 4.12]{valesin2024contact}, which is itself a reproduction of the method in~\cite{bhamidi2021survival}. There are only two differences worth mentioning in our current setting. First, the infection rate from the root may vary from neighbour to neighbour; this has no impact on the proof. Second, our graph~$G$ is not a tree. However, in the modified dynamics described above, the edges between neighbours of~$o$ are irrelevant and can be removed, making the dynamics in the subtrees rooted at neighbours of~$o$ independent in periods when~$o$ is infected; this is all that matters for the proof.
\end{proof}
\begin{remark}\normalfont
The identity above says the following: the original contact process on an augmented tree can be stochastically dominated from above by a \emph{degree-dependent reinforced contact process} on the non-augmented tree. In this process, a vertex does not heal until its whole subtree is healthy, and each infected vertex $v$ sends the infection to its descendant $u$ at rate $\lambda \mathrm{comm}_*(u)$. In the reinforced contact process  infections upwards in the tree are irrelevant as the parent of any infected vertex is per definition infected. 
\end{remark}
\begin{proof}[Proof of Theorem~\ref{thm:Existence-of-Subcritical-Phase}]
The arch of the proof takes ideas from  ~\cite[Section 4.3]{valesin2024contact}, and also from~\cite{bhamidi2021survival} to analyse the degree dependent reinforced contact process. 
For~$L \in \mathbb N_0$, let~$\mathcal T_L$ denote the graph obtained from~$\CT$ by excluding all vertices at distance larger than~$L$ from the root. Let~$f(L)$ denote the expectation of the extinction time of the reinforced contact process on~$\mathcal T_L$ started from only the root infected (integrating over the randomness both of the random graph and of the contact process).
Note that~$f(0)=1$. We will prove that~$f(L)\le 2$ for all~$L$ when~$\lambda$ is small enough. For an inductive argument, fix~$\lambda$ small (it will be chosen during the proof) and assume that~$f(\ell)\le 2$ has already been established for all~$\ell \le L$.
For a deterministic graph~$G$ with distinguished vertex~$o$, let~$h(G,o)$ denote the expected extinction time of the reinforced contact process on~$G$ with rate~$\lambda$ started with only~$o$ infected. 
 Recall from Section \ref{sec:Definition-community-tree} that $K_{o,i}$ is the number of additional vertices beyond the root  $o$ in the $i$th community of the root $o$. As then $\mathrm{comm_*(v)}=(K_{o,i}+1)-1=K_{o,i}$ for each non-root member $v\in \mathcal K_{o,i}$, we obtain by Proposition~\ref{prop_augmented} that
\[
h(\mathcal T_{L+1},o) \le \prod_{i=1}^{M_o} \prod_{v \in \cK_{o, i}\setminus \{o\}} (1 + \lambda \cdot K_{o,i} \cdot h(\mathcal T_{L,v},v)),
\]
where~$\mathcal T_{L,v}$ is the sub-graph of~$\mathcal T_{L+1}$ induced by~$v$ and its descendants. As $h(\mathcal{T}_{L,v},v)$ is (conditionally) independent across $v$ in the above double product, taking conditional expectation gives
\begin{equation*}
\mathbb E\left[h(\mathcal T_{L+1},o) \mid M_o, (K_{o,i})_i \right] \le \prod_{i=1}^{M_o} (1 + \lambda \cdot K_{o,i} \cdot f(L))^{K_{o,i}}.
\end{equation*}
By Definition \ref{def:branching-community-tree}, each $K_{o,i}$ is independent and has the same distribution as $K$. Thus taking conditional expectation given $M_o$,   
\begin{equation}
\label{eq_inter_M}
\mathbb E\left[h(\mathcal T_{L+1},o) \mid M_o \right] \le  \mathbb E\left[(1 + \lambda \cdot K \cdot f(L))^{K} \right]^{M_o}.
\end{equation}
We now fix~$c > 0$ such that~$\mathbb E[(1+c(K+1))^K]<\infty$ and bound
\[
\mathbb E\left[(1 + \lambda \cdot K \cdot f(L))^{K} \right]
\le \mathbb E\left[(1 + 2\lambda K)^K\right],
\]
using the induction hypothesis that $f(L)\le 2$. If~$\lambda \le c/2$, the integrand on the right-hand side is bounded by~$(1+cK)^K$. Hence, applying dominated convergence on the collection of random variables $X(\lambda):=(1+2\lambda K)^{K}$, for any $c'>0$ we can find that for all $\lambda$ sufficiently small we have
\begin{equation}
\label{eq_inter_for_K}
\mathbb E\left[(1 + 2\lambda K)^K\right]
\le e^{c'}.
\end{equation}
By the assumption that~$M$ has a finite exponential moment, there exists~$c' > 0$ such that~$\mathbb E[e^{c' M}] < 2$. Substituting this back in~\eqref{eq_inter_M} and integrating with respect to the law of $\mathcal T_{L+1}$,
\[
f(L+1) = \mathbb E[h(\mathcal T_{L+1},o)] \le \mathbb E[e^{c' M}] < 2.
\]
This completes the proof of the induction. 
To conclude, we note that the extinction time of the reinforced contact process on~$\CT$ started from the root infected is equal to~$\lim_{L \to \infty}f(L)$, by monotone convergence. Since this is finite when~$\lambda$ is small, the extinction time itself is almost surely finite. Since the reinforced contact process dominates the contact process. The proof is complete by noting that finiteness of $\E[(1+c K)^K]$ is equivalent to that of $\E[e^{K \log (cK)}]$ by elementary bounds.  
\end{proof}
Observe that the branching community tree is not a maximal augmented tree in the sense that not every child of a single vertex is connected to each other, only the ones in the same community. One may wonder about condition for the existence of a subcritical phase in maximal augmented trees. In fact our theorem above provides the answer by setting the distribution $M\equiv 1$ in Definition \ref{def:branching-community-tree}. 
\begin{example}
Consider the maximal augmented tree arising from a  Bienaym\'e-Galton-Watson tree with offspring distribution $X$, and maximal augmentation, i.e., when every child of a vertex $v$ forms a complete graph.  
As a corollary of the above theorem  we obtain that when $\E[(1+cX)^X]<\infty$ for some $c>0$, then $\lambda_1>0$.
\end{example}
\section{Exponentially long persistence in the network size}\label{sec:long-survival-on-rig}

In this section we consider a contact process $\bfxi=\bfxi^{(n)}$ on a random intersection graph $G_n:=\mathrm{RIG}_n(\bfm^{(n)}, \bfk^{(n)})$ satisfying Assumption~\ref{asmp:convergence} with limiting distributions $\offmem$ and $\offcomm$. Furthermore, recall the that $T^{n}_{\mathrm{ext}}=T^{G_n}_{\mathrm{ext}}$ is the extinction time with initially every individual infected, i.e.\ $\bfxi_0=G_n$. The aim of this section is to show that for every $\lambda$, $T^{n}_{\mathrm{ext}}$ is of order $\exp(\Theta(n))$ as $n\to \infty$ with high probability. 

\subsection{Upper bound via limited individual vertex contributions}

 We first begin with the proof of the upper bound showing that the extinction happens on an exponential scale. The key idea is to introduce events that imply no contribution for vertices, that work on time-scales relevant to the degree of the vertex under consideration.  In this section, $\deg(x)$ denotes the degree of vertex $x$ in a graph $G_n$.
\begin{proof}[Proof of Theorem \ref{thm:general-exponential-extinction}]
 We set out to prove that  there exists $c>0$ such that
\begin{equation}\label{eq:Bound-dying-out-after-1}
	\IP(T_{\text{ext}}^{n}<1)\geq e^{- cn}.
\end{equation}
From this inequality the claim follows immediately from the strong Markov property combined with the fact that the process is monotone with respect to the set of infected vertices, and in $T_{\text{ext}}^{n}$ we use the all-infected state as starting state. Indeed, suppose \eqref{eq:Bound-dying-out-after-1} is true. Let $X_n\sim \text{Bin}(\lfloor e^{2cn}\rfloor, e^{- cn})$, then using the strong Markov property to restart after every time step $t$ to $t+1$ and the fact that disjoint parts of the graphical representation are independent allows us to conclude
\begin{equation*}
	\IP(T_{\text{ext}}^{n}> e^{2cn})\leq \IP(X_n=0)\to 0
\end{equation*}
as $n\to \infty$, since $\lfloor e^{2cn}\rfloor e^{- cn}\to \infty$.

Thus, it remains to to show \eqref{eq:Bound-dying-out-after-1}. Recall that $|V(G_n)|=n$. 
First we work conditionally on a realisation of $G_n$. We denote by $\deg(x)$ the degree of vertex $x$ in this realisation.
Recall from Section \ref{sec:graphical-rep} that $(\mathcal R_x)_{x\in V(G_n)}$ and $(\mathcal I_e)_{e\in \vec{E}(G_n)}$ are the Poisson processes governing the recovery and infection times of the contact process graphical representation.  
We set $T^{\text{rec}}_x:=\inf\{t\geq 0: t\in \cR_{x}\}$ and $T^{\text{inf}}_{(x,y)}:=\inf\{t\geq 0: t\in \cI_{(x,y)}\}$ to be the first recovery of $x$ and the first infection times out of $x$ towards $y$, respectively. We compute and write
\begin{equation}\label{eq:Lx}
	L_x:=\mathbb E[\min(T_x^{\mathrm{rec}}, \min_{y\sim x} T_{(x,y)}^{\mathrm{inf}})\mid \deg(x)]=\frac{1}{1+\lambda \deg(x)}\leq 1,
\end{equation}
the expected time it takes for either a recovery event or an out-going infection event from $x$ to a neighbour $y$ to be triggered. We denote by
\begin{equation*}
	A_x:=\{T^{\text{rec}}_x\le L_x\}\cap \bigcap_{y\in \cN_x}\{T^{\text{inf}}_{(x,y)}>L_x\},
\end{equation*}
and note that on $A_x$, $x$ recovers before it can send an infection to a neighbour. As the event $A_x$ captures that vertex $x$ itself does not contribute to the spread of infection, we may observe that if this event happens for all $x\in G_n$,  then infection dies out within a unit time step. Thus, the intersection over all $x$ is a subset $\{T_{\text{ext}}^{n}<1\}$, and thus in particular
\begin{equation*}
	\IP_\lambda(T_{\text{ext}}^{n}<1|G_n)\geq \IP\bigg(\bigcap_{x\in G_n}A_x\mid G_n\bigg).
\end{equation*}
We make the following important observation:  we only use outgoing edges in the definition of the event, implying that $(A_x)_{x\in \leftpar_n}$ is an independent collection of events. Therefore, we have that
\begin{equation*}
	\IP\bigg(\bigcap_{x\in G_n}A_x\mid G_n\bigg)=\prod_{x\in G_n} \IP (A_x\mid G_n)=
	\prod_{x\in G_n}\big(1-e^{-L_x}\big)e^{-\lambda\deg(x) L_x}.
\end{equation*}
As $L_x<1$, we may use that $1-e^{-L_x}\geq L_x/2$. Further, using the formula for $L_x$ in \eqref{eq:Lx}, we have that $\frac{\lambda\deg(x)}{1+\lambda\deg(x)}\leq 1$, so  $\exp(-\lambda \deg(x)L_x)\ge \exp(-1)$ for each $x\in G_n$.  We arrive at
\begin{equation*}
    \begin{aligned}
	   \IP\bigg(\bigcap_{x\in G_n}A_x\mid G_n\bigg)&    \geq  \exp\bigg(- n+\sum_{x\in G_n}\log\bigg(\frac{1}{2(1+\lambda \deg(x)}\bigg)\bigg)\\
	   &=\exp\bigg(-  n-\sum_{x\in G_n}\log\big(1+\lambda \deg(x)\big)\bigg)  .
    \end{aligned}
\end{equation*}
We now take expectation over the law of $G_n$. As the function $e^{-y}$ is convex, we may use Jensen's inequality.
\begin{equation}\label{eq:intermediated-upper-bound}
     \begin{aligned}
	   \IP\bigg(\bigcap_{x\in G_n}A_x\bigg)& \ge  e^{-n} \IE\bigg[  \exp\bigg(-\sum_{x\in G_n}\log\big(1+\lambda \deg(x)\big)\bigg)\bigg]\\
	   &=e^{-n}\exp\bigg(-n\E\bigg[\frac1n\sum_{x\in G_n}\log\big(1+\lambda \deg(x)\big)\bigg]\bigg) \\
       &= e^{-n} \exp\left(-n\E[\log(1+\lambda D_n]\right).
    \end{aligned}
\end{equation}
By assumption and subadditivity of $x\mapsto\log(1+\lambda x)$ for $x\geq 0$, there exists a constant $\Delta_{\lambda}>0$ such that $\E[\log(1+\lambda D_n)]<\Delta_{\lambda}$ for all $n$ sufficiently large. By applying this to \eqref{eq:intermediated-upper-bound} we obtain
\begin{equation*}
     \begin{aligned}
	   \IP\bigg(\bigcap_{x\in G_n}A_x\bigg)& \ge \exp\big(-(1+\Delta_{\lambda})n\big),
    \end{aligned}
\end{equation*}
which proves the claim.
\end{proof}
After proving the more general Theorem~\ref{thm:general-exponential-extinction}, Theorem~\ref{thm:rigc-extinction} follows as corollary.
\begin{proof}[Proof of Theorem~\ref{thm:rigc-extinction}] 
Recall that we assumed that $M_n\stackrel{d}{\longrightarrow}M$ and $K_n\stackrel{d}{\longrightarrow}K$ with means converging as well. Let $D_n$ be the degree of a uniformly chosen vertex in the random itersection graph.

We shall show not only the condition of Theorem \ref{thm:general-exponential-extinction}, but the stronger statement that
$\E[\log(1+ c D_n)]$ converges to $\E[\log (1+ c D)]<\infty$, where $D=\sum_{j=1}^{M} \widetilde K^{(j)}$ a random sum of iid entries. Note that $D_n$ itself cannot be written as a sum of independent entries. We now describe its distribution exactly.
In a random intersection graph, the degree of a uniformly chosen vertex is determined by the random matching in the bipartite representation of the model in Section \ref{sec:rig-def}. Given the community size list $\mathbf k=(k_w)_{w\in [a_n]}$, create a list $h_n=\sum_{w\in[a_n]}k_w$ items, where $k_w$ many items are labelled with belonging to $w$. Let us sample from this list sequentially \emph{without replacement}, obtaining $(\bar K_n^{(1)}, \bar K_n^{(1)}, \dots, \bar K_n^{(h_n)})$, where $\bar K_n^{(j)}$ is the size $-1$ of the community that the $j$th sampled item belongs to. As the sampling is governed by a random permutation, the random variables $(\bar K_n^{(1)}, \bar K_n^{(1)}, \dots, \bar K_n^{(h_n)})$ are \emph{exchangeable}. 
Let now $U_n$ be a discrete uniform variable in $[n]$, and $m_{U_n}=M_n$ be the membership number of a uniformly chosen vertex $U_n$. Then the degree of vertex $U_n$ in the corresponding random intersection graph is exactly $D_n\ {\buildrel d \over =}\ \sum_{j=1}^{M_n}\bar K_{n}^{(j)}$. This is due to the rules of drawing multiple edges and loops in \eqref{eq:edge-rule-rig} and \eqref{eq:loop-rule-rig}, as if we match $\ell$ half-edges of $U_n$ to community $w$, then the contribution of that community to the degree is exactly $\ell (k_w-1)$.   
We now show that $\log(1+ c D_n)$ is uniformly integrable. 
We can use the sub-additivity of the logarithm, followed by Jensen's inequality $(n^{-1}\sum_{i=1}^m x_n)^{1+\varepsilon}\le n^{-1}\sum_{i=1}^m x_n^{1+\varepsilon} $ to obtain
\begin{equation}
    \begin{aligned}
    \log^{1+\varepsilon}(1+ cD_n)&=
    \Big(\log (1+ c \sum_{j=1}^{M_n}\bar K_{n}^{(j)})\Big)^{1+\varepsilon}\\ &\le \Big(\sum_{j=1}^{M_n}\log (1+ c\bar K_{n}^{(j)}  )\Big)^{1+\varepsilon} 
        \le  M_n^{\varepsilon} \sum_{j=1}^{M_n} \log^{1+\varepsilon}(1+ c \bar K_{n}^{(j)} )
    \end{aligned}
       \end{equation}
       We take expectation, apply law of total probability, then we can use exchangeability by noting that the marginals of $\bar K_n^{(\cdot)}$ do not depend on $M_n$. This gives
       \begin{equation}
       \begin{aligned}
      \E&\big[\log^{1+\varepsilon}(1+ cD_n)\big]\le \sum_{m=1}^{\infty}\P(M_n=m)m^{\varepsilon}\E\Big[\sum_{j=1}^m \log^{1+\varepsilon}(1+ c \bar K_n^{(j)})\Big]\\
&=\E[M_n^{1+\varepsilon}] \E[\log^{1+\varepsilon}(1+ c \bar K_n^{(1)})].
    \end{aligned}
\end{equation}
 For $\varepsilon<\eta$, we have $\lim_{n\to \infty}\E[M_n^{1+\varepsilon}]=\E[M^{1+\varepsilon}]$, hence the first factor on the right hand-side is uniformly bounded in $n$. To bound the second factor, we use that the first half-edge is uniform among $h_n=\E[K_n] a_n$ many half-edges. 
\begin{equation}
 \begin{aligned}
        \E\Big[\log^{1+\varepsilon} (1+ c\bar K_{n}^{(1)}  )\Big] &= \frac{1}{a_n\E[K_n]}  \sum_{w\in [a_n]} k_w \log^{1+\varepsilon} (1+ c(k_w-1)) \\&= 
        \frac{\E[K_n \log^{1+\varepsilon}(1+ c K_n)]}{\E[K_n]}.
    \end{aligned}
\end{equation}
As the function $g(k):=k \log^{1+\varepsilon}(1+ ck)$ satisfies $g(k)=o(k^{1+\eta})$, and $\lim_{n\to \infty}\E[K_n^{1+\eta}]=\E[K^{1+\eta}]<\infty$ per assumption, we have that $(g(K_n))_{n\ge 1}$ is uniformly integrable. Also $g(K_n)\ {\buildrel d \over\to}\  g(K)$, so the  the right hand side converges to  $\E[K \log^{1+\varepsilon}(1+ cK)]/\E[K]<\infty$. This gives the uniform integrability of $\log(1+c D_n)$. Using the representation $D_n=\sum_{j=1}^{M_n}\bar K_{n}^{(j)})$, we see that $D_n$ converges in distribution to $D:=\sum_{i=1}^{M}\widetilde K^{(i)}$, where $M$ is the distributional limit of $M_n$  and $\widetilde K^{(i)}$ are iid from down-shifted size-biased version of $K$, the distributional limit of $K_n$. This gives that 
$\E[\log(1+ c D_n)]$ converges to $\E[\log (1+ c D)]$, and establishes the condition of the Theorem \ref{thm:general-exponential-extinction}. 
\end{proof}

\section{Persistence lower bound via embedding linearly many communities}\label{sec:persistence-in-rig}
In this section we prove that the contact process persists exponentially long in the finite graph. We carry out the embedding of beads as described in the introduction, see Figure \ref{fig:bead-embedding}. Instead of working directly with the random intersection graph, we work with the bipartite configuration model representation that generated it, see Section \ref{sec:rig-def}. Once the matching in the bipartite representation is fixed, the random intersection graph is fully determined, as there is a deterministic map from the bipartite representation to the final graph. We make use of this by searching for the relevant `beads' in the bipartite representation. This has the advantage that large communities and vertices with high membership number can be `symmetrically represented', that is, our algorithm works in both cases when we search for high degree nodes vs when we search for large enough communities to sustain the infection. 
\subsection{Finding an embedded network in a bipartite configuration model}

As explained above, we will consider the bipartite configuration model as given in Definition~\ref{def:bipartite-conf-model}, except that we will not think of vertices as ``members'' and ``communities'' as we did there. Rather, we will simply split vertices into ``bottom vertices'' and ``top vertices''. The reason for this change is that we will prove a meta-result (expressed in the ``bottom-top'' setting), and then  apply this result once with members in the bottom and communities in the top, and later with these reversed, depending on which exponential moment condition in Theorem \ref{thm:survival-finite-graph} is satisfied.

The letter~$y$ will be associated to the bottom, and the letter~$z$ to the top.
The set of bottom vertices is~$\mathcal V^{\mathscr{y}}_n$, and the set of top vertices is~$\mathcal V^{\mathscr{z}}_n$. We express the degrees as
\[
\mathbf{y}=\mathbf{y}^{(n)}=\left( y^{(n)}_x\right)_{x \in \mathcal V^{\mathscr{y}}_n} \quad \text{and} \quad \mathbf{z}=\mathbf{z}^{(n)}=\left( z^{(n)}_x\right)_{x \in \mathcal V^{\mathscr{z}}_n},\quad n \in \mathbb N.
\]
We assume that
\begin{equation}
    \label{eq_divergence_Vs}
\frac{|\mathcal V^{\mathscr{y}}_n|}{n} \text{ and } \frac{|\mathcal V^{\mathscr{z}}_n|}{n} \text{ converge to positive constants as~$n \to \infty$}
\end{equation}
and
\[\sum_{x \in \mathcal V^{\mathscr{y}}_n} y^{(n)}_x =\sum_{x \in \mathcal V^{\mathscr{z}}_n} z^{(n)}_x.
\]
We define the probability vectors~$\mathbf{p}^{(\mathscr{y},n)}=(p_\ell^{(\mathscr{y},n)})_{\ell \in \mathbb Z^+}$ and~$\mathbf{p}^{(\mathscr{z},n)}=(p_\ell^{(\mathscr{z},n)})_{\ell \in \mathbb Z^+}$ as
\begin{equation}\label{eq_empirical_yz}
    p_\ell^{(\mathscr{y},n)}:=\frac{|\{x \in \mathcal V^{\mathscr{y}}_n: y^{(n)}_x=\ell\}|}{|\mathcal V^{\mathscr{y}}_n|},\qquad p_\ell^{(\mathscr{z},n)}:=\frac{|\{x \in \mathcal V^{\mathscr{z}}_n: z^{(n)}_x=\ell\}|}{|\mathcal V^{\mathscr{z}}_n|}.
\end{equation}

For convenience, we restate Assumption~\ref{asmp:convergence} with the new notation:
\begin{assumption}
\label{asmp:convergence_YZ} Letting~$Y_n \sim \mathbf p^{(\mathscr{y},n)}$ and~$Z_n \sim \mathbf p^{(\mathscr z,n)}$ for each~$n$, we assume that there exist random variables~$Y$ and~$Z$ and~$\eta>1$ such that
		\begin{equation*}
			Y_n \xrightarrow[\mathrm{d}]{n \to \infty} Y,
            \qquad \mathbb E[Y_n^{\eta}] \xrightarrow{n \to \infty} \mathbb E[Y^{\eta}]\in \R_+,
		\end{equation*}
		and similarly,
		\begin{equation*}
			Z_n \xrightarrow[\mathrm{d}]{n \to \infty} Z,
            \qquad \mathbb E[Z_n^\eta] \xrightarrow{n \to \infty} \mathbb E[Z^\eta]\in \R_+.
		\end{equation*}
\end{assumption}

Recall the definition of alternating branching trees (Definition~\ref{def_alternating}).   Let~$\mathcal T$ be such a tree, and let~$d \in \mathbb N$. A \emph{breadth-first search of~$\mathcal T$ for vertices of degree~$d$ in even generations} is a procedure in which we progressively reveal the number of children of vertices in~$\mathcal T$ in a breadth-first order (that is, exhausting a generation before moving on to the next generation), with the caveat that whenever we reveal the number of children of a non-root vertex in an even generation, and this number is exactly~$d-1$, we do not inspect its descendants in the exploration (so that the whole tree descending from it is left unexplored).

We fix natural numbers~$\mathsf d$ and~$\mathsf r$ and assume the following.

\begin{assumption}\label{assumption_new}
	We assume that~$\{(\mathbf y^{(n)}, \mathbf{z}^{(n)})_n,\mathsf d, \mathsf r\}$ satisfy  the following assumptions:
    \begin{itemize}
        \item[$\mathrm{(a)}$] Assumption~\ref{asmp:convergence_YZ} is satisfied by~$\{(\mathbf y^{(n)}, \mathbf z^{(n)})\}$, with limiting random variables~$Y$ and~$Z$.
        \item[$\mathrm{(b)}$] We have
\begin{equation}\label{eq_extra_asmp_candidate}
\mathsf d\cdot \mathbb P(Y=\mathsf d) < \mathbb E[Y]/16.
        \end{equation}
\item[$\mathrm{(c)}$] Let $\widetilde Y$ and $\widetilde Z$ be the down-shifted size-biased version of $Y$ and $Z$ given in \eqref{eq:size-biasing-def}. 
    
    Let~$\mathcal T$ be distributed as~$\mathrm{ABT}(\mathsf d,\widetilde{Y},\widetilde{Z})$ {(that is, 
the root has deterministically $\mathsf d$ children, $\widetilde Y$ is the offspring distribution in even generations and $\widetilde Z$ is the offspring distribution in odd generations.)
    }
		
        Let~$\mathcal A$ be the event that in a breadth-first search of~$\mathcal T$ for vertices of degree~$\mathsf d$ in even generations, there are at least two instances when the search succeeds within distance~$\mathsf r$ of the root (in other words, the exploration reveals at least two non-root
        vertices in even generations within distance~$\mathsf r$ of the root with degree equal to~$\mathsf d$).
        
        We assume that~$\mathbb P(\mathcal A)>3/4$.
        \end{itemize}
\end{assumption}

\begin{definition}\label{def_star_line}
    A connected tree~$G=(V,E)$ with a root vertex~$o$ is called a \emph{star-line tree} (with degree parameter~$\mathsf d$ and distance parameter~$\mathsf r$) if there exist~$V' \subset V$ with
    \[o \in V'\quad \text{and}\quad \deg(x)=\mathsf d \text{ for all } x \in V'\]
    and a \emph{parent function}~$\mathsf p: V' \setminus\{o\} \to V'$ with~$\mathsf p(x) \neq x$ for all~$x \in V'$, and such that the following holds:
    \begin{itemize}
        \item the graph~$G':=(V',E')$, where~$E':=\{\{\mathsf p(x),x\}: x \in V' \setminus \{o\}\}$, is a connected tree with degrees bounded by~$3$; we call $\{\mathsf{p}(x),x\}$ a \emph{renormalised edge};
        \item letting~$[\mathsf p(x),x]$ denote the set of vertices in the geodesic in~$G$ from~$\mathsf p(x)$ to~$x$ for each~$x \in V' \setminus \{o\}$, we have
        \[|[\mathsf p(x),x]| \le \mathsf {r+1} \quad \text{and} \quad V' \cap [\mathsf p(x),x] = \{\mathsf p(x),x\};\]
        \item letting~$\mathcal N_x$ denote the set of neighbours of~$x$ in~$G$, we call $\{x\}\cup \mathcal N_x$ a \emph{star} and 
        \[
        V = \mathcal N_o \cup  \bigcup_{x \in V' \setminus \{o\}} (\mathcal N_x \cup [\mathsf p(x),x]).
        \]
    \end{itemize}
\end{definition}
{A star-line tree, in words, contains the `center' vertices forming $V'$, which are arranged in a `renormalised' tree structure of maximal degree $3$ using the parent-child relationships (the renormalised edges). Each parent-child relationship corresponds to a path between the center vertices of length at most $r$ edges, that may or may not overlap between different children of the same parent vertex. Additionally, each center vertex has surrounding neighbors so that its total degree is exactly $d$. On Figure \ref{fig:bead-embedding}, a similar object is shown, where the center vertices are replaced by communities insted of `stars'. }
Our main goal in this section is to prove the following.
\begin{proposition}\label{prop:suitable-finite-subgraph}
	
    Assume that~$\{(\mathbf y^{(n)}, \mathbf{z}^{(n)})_n,\mathsf d, \mathsf r\}$ satisfy Assumption~\ref{assumption_new}. Then, there exists~$\delta_0 > 0$ such that with probability tending to~$1$ as~$n \to \infty$, the bipartite configuration model~$\mathrm{BC}_n(\mathbf y, \mathbf z)$ has a subgraph which is a star-line graph (with degree parameter~$\mathsf d$ and distance parameter~$\mathsf r$) in which all stars are bottom vertices,  containing~$\delta_0 n$ stars.
\end{proposition}

The proof of this proposition is done via an exploration algorithm of the random bipartite graph, which we develop and analyze in Sections~\ref{ss_prelim_algo}--\ref{ss_algo_analysis}.

\subsubsection{Proof of Theorem~\ref{thm:survival-finite-graph} with~${\displaystyle \mathbb E[e^{c\widetilde{K}\log \widetilde{K}}]=\infty}$ for all~$c$}

\begin{lemma}\label{lem_asmp_yonly}
     Let~$(\mathbf y^{(n)})_n$ and~$(\mathbf z^{(n)})_n$ be degree sequences satisfying Assumption~\ref{asmp:convergence_YZ}, with limiting distributions~$Y$ and~$Z$ satisfying~$\mathbb E[\widetilde{Y}] \mathbb E[\widetilde{Z}] > 1$ and~$\mathbb E[e^{c\widetilde{Y}\log \widetilde Y}]=\infty$ for all~$c > 0$. Then, for all~$\lambda \in (0,1)$, there exists~$\mathsf d_\star = \mathsf d_\star(\lambda)$ such that for all~$\mathsf d \ge \mathsf d_\star$, letting
     \begin{equation}\label{eq:r-of-d-def}
     \mathsf r =\mathsf r(\mathsf d,\lambda):=\lfloor \log(2/\lambda)^{-1}\mathsf d \log(\lambda \mathsf d) \rfloor,
     \end{equation}
     we have that $\{(\mathbf y^{(n)}, \mathbf{z}^{(n)})_n,\mathsf d, \mathsf r\}$ satisfies Assumption~\ref{assumption_new}.
    \end{lemma}
\begin{proof}
Part~$\mathrm{(a)}$ of Assumption~\ref{assumption_new} is satisfied by assumption. We will prove that there is a~$\lambda$-dependent infinite set of values~$\mathsf d$ and corresponding choice $\mathsf{r}$ in \eqref{eq:r-of-d-def} such that part~$\mathrm{(c)}$ is satisfied by~{$\mathrm{ABT}(\mathsf{d}, \widetilde Y, \widetilde Z)$ with distance $\mathsf{r}$}.
Then, by taking one such~$\mathsf d$ to be sufficiently large, we can also ensure that part~$\mathrm{(b)}$ is satisfied. 

We start by proving an auxiliary result. Note that in the following claim the switch between odd and even generation distribution is intentional. 
\begin{claim}\label{cl_seneta}
For~$\mathcal T' \sim \mathrm{ABT}(1,\widetilde{Z}, \widetilde{Y})$, we have
\[
\mathbb P(|\mathcal T'|=\infty)>0
\]
and with $\mu:=\mathbb E[\widetilde{Y}]\mathbb E[\widetilde{Z}]$ we have
\begin{equation}\label{eq_liminf_exponent}
\liminf_{m \to \infty} |\{v \in \mathcal T': d(o,v)=m\}|^{1/m} \ge \sqrt{\mu} \text{ almost surely on } \{|\mathcal T'|=\infty\}.
\end{equation}
\end{claim}
\begin{proof}
The process~$W=(W_m)_{m \ge 0}$ given by~$W_m:=|\{v \in \mathcal T': d(o,v) = 2m\}|$ is a branching process with mean offspring~$\mu$, so the first statement follows from our assumption that~$\mu > 1$. 

Next, the Seneta-Heyde theorem implies that there exists a sequence $(c_m)_{m\geq 0}$ with $\tfrac{c_{m+1}}{c_{m}}\to \mu$ and $c_{m}\to \infty$ as $m\to \infty$ and a random variable $\mathcal W$ such that
    \begin{equation*}
        \frac{W_m}{c_m}\to \mathcal W \text{ a.s. } \quad \text{ and }\quad  \{\mathcal W>0\}=\{W \text{ survives}\}.
    \end{equation*}
Fix~$\mu_\star < \mu$. By the ratio assumption on $c_m$, it follows that $\frac{c_m}{(\mu_\star)^m}\to \infty$ as~$m \to \infty$. Thus, we obtain
    \begin{equation*}
        \begin{aligned}
            &\IP(W_m\geq (\mu_\star)^m \text{ for all }m \ge m_0)\\
            &\geq \IP(\{W_m\geq (\mu_\star)^m\text{ for all }m \ge m_0\} \cap \{W \text{ survives}\})\xrightarrow{m_0 \to \infty} \mathbb P(W \text{ survives}).
        \end{aligned} 
    \end{equation*}
    This proves that 
    \[
    \liminf_{m \to \infty} |\{v \in \mathcal T': d(o,v)=2m\}|^{1/m} \ge \mu \text{ almost surely on } \{|\mathcal T'|=\infty\}.
    \]
To prove that
\[
    \liminf_{m \to \infty} |\{v \in \mathcal T': d(o,v)=2m+1\}|^{1/m} \ge \mu \text{ almost surely on } \{|\mathcal T'|=\infty\},
    \]
we repeat the argument for the process~$(|\{v\in \mathcal T': d(o,v)=2m+1\}|)_{m \ge 0}$.
\end{proof}

Fix~$\mathsf d$ to be chosen later, and~$\mathsf r:=\lfloor \log(2/\lambda)^{-1}\mathsf d \log(\lambda \mathsf d) \rfloor$ as in the statement of the lemma, and also let~$\mathsf r'$ be the largest even number smaller than or equal to~$\mathsf r$.
Let~$\mathcal T \sim \mathrm{ABT}(\mathsf d, \widetilde{Y},\widetilde{Z})$, and let~$\mathcal A$ be the event described in part~$\mathrm{(c)}$ of Assumption~\ref{assumption_new}. Define
\begin{align*}
    &X_{\mathsf d} := |\{u \in \mathcal T: d(o,u)=\mathsf r'\}|, \qquad X_{\mathsf d}':=|\{u \in \mathcal T: d(o,u)=\mathsf r',\;\deg(u)=\mathsf d\}|.
\end{align*}
Then,
\begin{equation*}
    \mathcal A \supseteq \{X_{\mathsf d}' \ge 2\}.
\end{equation*}
For any~$k \in \mathbb N$ we then have the lower bound: 
\begin{equation}\label{eq_even_branching}
\mathbb P(X_{\mathsf d}' \ge 2) \ge {\mathbb P(X_{\mathsf d}' \ge 2 \mid X_{\mathsf d} = k)} \cdot {\mathbb P(X_{\mathsf d} \ge k)}.
\end{equation}
Fix~$\mu_\star \in (1,\mu)$. We now prove that taking~$\mathsf d$ large enough and
\begin{equation}\label{eq:k-d-def-temp}
k=k(\mathsf{d}):=\mu_\star^{\mathsf r'/2},
\end{equation}
the right-hand side of~\eqref{eq_even_branching} is larger than~$7/8$.

We start giving a lower bound on~$\mathbb P(X_{\mathsf d} \ge k)$. For each children~$v$ of the root of~$\mathcal T$ in $\mathrm{ABT}(1,\widetilde{Y}, \widetilde{Z})$, let~$\mathcal T^{(v)}$ be the subtree including~$v$ and all its descendants, rooted at~$v$. Note that these trees are independent, all with the distribution of~$\mathrm{ABT}(1,\widetilde{Z},\widetilde{Y})$, (with even and odd generation distribution switched). 
Let~$\theta:=\mathbb P(|\mathcal T^{(v)}|=\infty)$.
Using Claim~\ref{cl_seneta} and the fact that~$\mu_\star < \mu$, by making~$\mathsf r'$ large enough (by increasing~$\mathsf d$), we have
\[
\mathbb P(|\{u \in \mathcal T^{(v)}: d(v,u) = \mathsf r'-1\}| \ge k) > \theta/2.
\]
Then,
\begin{equation}\label{eq_binomial_theta}\begin{aligned}
    \mathbb P(X_{\mathsf d} \ge k) &\ge \mathbb P(\exists v: v \sim o, \; |\{u \in \mathcal T^{(v)}: d(v,u)=\mathsf r'-1\}| \ge k)\\
    &\ge \mathbb P(\mathrm{Bin}(\mathsf d, \theta/2) \ge 1)  = 1-(1-\theta/2)^{\mathsf d}.
\end{aligned}\end{equation}
Next, writing~$\mathsf q(\mathsf d):=\mathbb P(\widetilde{Y}=\mathsf d)$, we bound
\begin{equation}\label{eq_mathsfq}
    \mathbb P(X_{\mathsf d}'\ge 2 \mid X_{\mathsf d}=k) = \mathbb P(\mathrm{Bin}(\lfloor k \rfloor, \mathsf q(\mathsf d))  \ge 2).
    \end{equation}
Note that by \eqref{eq:Uniform-lower-bound-Comm-Dist}, the assumption that $\IE\big[e^{c\widetilde{Y}\log(\widetilde{Y})}\big]=\infty$ for all $c>0$ implies that
    \begin{equation*}
       \mathsf q(d)= \IP(\widetilde{Y} = d)\geq \exp(- \tfrac{1}{4}\log(\mu_{\star})\log(2/\lambda)^{-1} d \log(\lambda d)) \text{ for infinitely many }d \in \mathbb N.
    \end{equation*}
    We bound now the expectation of the Binomial distribution in question.  By the choice of $k(\mathsf d)$ in \eqref{eq:k-d-def-temp} and the choice of $r'$ being the largest even integer at most $\lfloor \log(2/\lambda)^{-1}\mathsf d \log(\lambda \mathsf d) \rfloor$, there exist infinitely many choices of $\mathsf d$ with $\lambda \mathsf d\geq 640$ that satisfy 
    \begin{align*}
        \lfloor k(\mathsf d) \rfloor \cdot \mathsf q(\mathsf d) =\lfloor (\mu_{\star})^{\mathsf r'/2} \rfloor \cdot \mathsf q&\geq \exp(\tfrac{1}{4}\log(\mu_{\star})\log(2/\lambda)^{-1}\mathsf d\log(\lambda \mathsf d))\\
        &\geq \exp(\tfrac{1}{4}\log(\mu_{\star})\log(2/\lambda)^{-1}\mathsf d\log(640)).
    \end{align*}
    The right-hand side tends to $\infty$ as $\mathsf d\to \infty$ at least linearly fast. Hence, by taking~$\mathsf d$ among the values that satisfy the above inequality, and sufficiently large, we can make the expectation of~$\mathrm{Bin}(\lfloor k(\mathsf d) \rfloor,\mathsf q(\mathsf d))$ as large as desired; in particular, we can make the right-hand side of~\eqref{eq_mathsfq} as close to 1 as desired. Together with~\eqref{eq_binomial_theta}, this proves that the right-hand side of~\eqref{eq_even_branching} can be made above~$7/8$, as required.
\end{proof}

\begin{proof}[Proof of Theorem~\ref{thm:survival-finite-graph} under the assumption~${\displaystyle \mathbb E[e^{c\widetilde K \log \widetilde K}]=\infty}$ for all $c>0$]

It suffices to treat~$\lambda \in (0,1)$. For larger~$\lambda$, the statement readily follows by monotonicity.

We apply our bottom-top setup for the bipartite configuration model with communities at the bottom and members at the top; in other words, we let~$(\mathbf y, \mathbf z)=(\mathbf k, \mathbf m)$. Since Assumption~\ref{asmp:convergence} (or equivalently, Assumption~\ref{asmp:convergence_YZ}) is in place, Lemma~\ref{lem_asmp_yonly} provides that there exists~$\mathsf d_{\star} > 0$ so that for all~$\mathsf d \ge \mathsf d_\star$, Assumption~\ref{assumption_new} is satisfied by
\begin{equation*}
    \{(\bfk^{(n)}, \bfm^{(n)})_n,\mathsf d, \mathsf r \}, \qquad \text{where } \mathsf r := \lfloor \log(2/\lambda)^{-1}\mathsf d \log(\lambda \mathsf d) \rfloor.
\end{equation*}
Thus, appealing to Proposition~\ref{prop:suitable-finite-subgraph}, with high probability the bipartite configuration model $\mathrm{BC}_n(\mathbf m, \mathbf k)$ contains a subgraph that is a star-line graph with degree parameter~$\mathsf d$ and distance parameter~$\mathsf r$, containing $\delta_0 n$ stars, all of which are community vertices. We denote by~$G=(V,E)$ this subgraph.

For this star-line graph, we will employ the terminology and notation of Definition~\ref{def_star_line}. In particular,~$V'$ denotes the vertices which are centres of stars (all of which are community vertices in the current setting), there is a root~$o \in V'$, and we will refer to the parent function~$\mathsf p$ and the line segments~$\mathsf p(w),w]$ for (non-root) vertices~$w \in V'$. Note that along this line segment, vertices alternate between~$\rightpar_n$ and~$\leftpar_n$.

Recall that the graph~$\mathrm{RIG}_n(\mathbf m,\mathbf k)$, whose vertex set is~$\mathcal V^{\mathscr m}_n$, is obtained from~$\mathrm{BC}_n(\mathbf{m},\mathbf{k})$ via the community projection. In particular, recall that if~$w$ is a community vertex of~$\mathrm{BC}_n(\mathbf m, \mathbf k)$, then we denote by~$\mathcal K^{(n)}_w$ the complete subgraph of~$\mathrm{RIG}_n(\mathbf m, \mathbf k)$ whose set of vertices is the set of~$v \in \mathcal V^{\mathscr m}_n$ that are neighbours of~$w$ in~$\mathrm{BC}_n(\mathbf m,\mathbf k)$. (We will abuse notation and refer to~$\mathcal K^{(n)}_w$ as the set of vertices of this complete graph).

We now want to describe a subgraph~$\mathscr G=(\mathscr V,\mathscr{E})$ of~$\mathrm{RIG}_n(\mathbf m,\mathbf k)$ that is obtained from~$G$ via this projection. We set
\[
\mathscr V:= \mathcal K^{(n)}_o \cup \bigcup_{w \in V'\setminus \{o\}}  (\mathcal K^{(n)}_w \cup ([\mathsf p(w),w] \cap \mathcal V^{\mathscr m}_n)). 
\]
In words,~$\mathscr V$ includes:
\begin{itemize}
    \item the communities of~$\mathrm{RIG}_n(\mathbf m,\mathbf k)$ corresponding to  the community vertices of~$\mathrm{BC}_n(\mathbf m,\mathbf k)$ that belong to~$V'$;
    \item the vertices in the line segments~$[\mathsf p(w),w]$ that belong to~$\mathcal V^{\mathscr m}_n$ (and are thus also vertices of~$\mathrm{RIG}_n(\mathbf m,\mathbf k)$).
\end{itemize}
The set of edges~$\mathscr E$ is defined simply as
\[
\mathscr E:=\{\{x,y\}: \{x,y\} \text{ is an edge of }\mathrm{RIG}_n(\mathbf m,\mathbf k),\; x,y \in \mathscr V\}.
\]
We see~$\mathscr G=(\mathscr V,\mathscr E) \subset \mathrm{RIG}_n(\mathbf m, \mathbf k)$ as a sort of image of~$G=(V,E) \subset \mathrm{BC}_n(\mathbf m,\mathbf k)$ under the community projection, this image being governed by the set of community-stars~$V'$ and the line segments connecting these community-stars.

Recall that in Definition~\ref{def_star_line}, we have defined the (renormalised) edge set~$E':=\{\{\mathsf p(x),x\}: x \in V' \setminus \{o\}\}$. {To accommodate for discrete time}, we now introduce the oriented graph $(\Lambda_0,\vec{E}_0)$, where
    \begin{equation*}
        \begin{aligned}
            \Lambda_0:=\{(w,n)\in V'\times\IN_0\}\quad \text{ and } \quad 
            \vec{E}_0:=\big\{\langle (w,n), (z,n+1)\rangle:  \{w,z\}\in E'\big\}.
        \end{aligned}
    \end{equation*}
    We now outline how the contact process on~$\mathscr G$ induces bond percolation on this oriented graph, similarly to the proof of Proposition~\ref{thm:Survival-Line-Of-Communities}.

    We define~$t_n:=n\exp(\frac{\mathsf d}{16}\log(\lambda \mathsf d))$ for~$n \ge 0$. Analogously as in \eqref{eq:CP-on-block-segments}, for~$\{w,z\} \in E'$, we define~$\bfzeta^{((w,z),n)}$,  the contact process restricted to $\cK^{(n)}_w\cup\cK^{(n)}_z\cup ([w,z]\cap \leftpar_n)$ in the interval $[t_n,t_{n+1}]$, started at time~$t_n$ from the configuration~${\bm \xi}_{t_n}$ restricted to this subgraph.

Now we define (again, as in the  proof of Proposition~\ref{thm:Survival-Line-Of-Communities}) a $1$-dependent oriented percolation $\bfX=(\bfX(\vec{e}))_{\vec{e}\in \vec{E}_0}$, with respect to the auxiliary processes $\bfzeta^{((w,z),n)}$. By choice $(\mathsf d,\mathsf r)$ satisfy \eqref{eq:Minimal-degree-of-a-star}, and thus appealing to Lemma~\ref{lem:passing-the-infection-on} we can choose $\mathsf d$ large enough such that  $\IP(\bfX(\vec{e})=1)$ is sufficiently close to $1$.

Using that $G'=(V',E')$ has degrees bounded by $3$ we see that $\bfX$ dominates an independent oriented percolation $\widehat{\bfX}$ with sufficiently high success probability $\IP(\widehat{\bfX}(\vec{e})=1)$, so that finally by \cite[Proposition 5.2]{mountford2016exponetial} the claim follows.
\end{proof}

\subsubsection{Proof of Theorem~\ref{thm:survival-finite-graph} with~${\displaystyle \mathbb E [e^{c\widetilde{M}}]=\infty}$ for all~$c>0$}

\begin{lemma}\label{lem_asmp_ylogy}
    Let~$(\mathbf y^{(n)})_n$ and~$(\mathbf z^{(n)})_n$ be degree sequences satisfying Assumption~\ref{asmp:convergence_YZ}, with limiting distributions~$Y$ and~$Z$ satisfying~$\mathbb E[\widetilde{Y}] \mathbb E[\widetilde{Z}] > 1$ and~$\mathbb E[e^{c\widetilde{Y}}]=\infty$ for all~$c > 0$.  Then, for all~$\lambda > 0$, there exists~$\mathsf d_\star = \mathsf d_\star(\lambda)$ such that for all~$\mathsf d \ge \mathsf d_\star$, letting
    \[\mathsf r = \mathsf r(\mathsf d,\lambda):= \lfloor C\log(\tfrac{1}{\lambda})^{-1} \lambda^2 \mathsf d \rfloor\}\] (where~$C$ is a constant that does not depend on~$\lambda$ or~$\mathsf d$), we have that $\{(\mathbf y^{(n)}, \mathbf{z}^{(n)})_n,\mathsf d, \mathsf r\}$ satisfies Assumption~\ref{assumption_new}. 
    \end{lemma}
The proof is very similar to that of Lemma~\ref{thm:survival-finite-graph}; we leave the verification to the reader.

\begin{proof}[Proof of Theorem~\ref{thm:survival-finite-graph} under the assumption~${\displaystyle \mathbb E[e^{c\widetilde M}]=\infty}$ for all $c$]
The proof under assumption~${\displaystyle \mathbb E[e^{c\widetilde M}]=\infty}$ for all $c$ can be conducted in an analogous manner as in the prior case, and there are only a few differences; most importantly that now stars (high degree vertices and their direct neighbors) carry the infection. Thus, we denote by $\stargraph^{(n)}_x=\stargraph_x$ a star graph with $x\in \RIG_n(\mathsf m,\mathsf k)$ at its center, this means $V(\stargraph_x)=\cN_x\cup\{x\}$ and $E(\stargraph_x)=\{\{x,y\}: y\in \cN_x \}$. Thus, this time we apply our bottom-top setup for the bipartite configuration model with members at the bottom and communities at the top, i.e.\ $(\mathbf y, \mathbf z)=(\mathbf m, \mathbf k)$. We again use Lemma~\ref{lem_asmp_ylogy} to see that 
    \begin{equation*}
        \{( \bfm^{(n)}, \bfk^{(n)})_n,\mathsf d, \mathsf r \}, \qquad \text{where } \mathsf r := \lfloor C\log(\tfrac{1}{\lambda})^{-1} \lambda^2 \mathsf d \rfloor.
    \end{equation*}
    satisfies Assumption~\ref{assumption_new}. 

    The coupling construction with the oriented percolation can be done analogously, with two main differences. First the block segments we restrict the process to are $\cS^{(n)}_x\cup\cS^{(n)}_y\cup ([x,y]\cap \leftpar_n)$, i.e.\ two connected stars instead of two connected communities. Secondly, we define~$t_n:=n\exp(c\lambda^2\mathsf d)$ for some $c>0$ and~$n \ge 0$, which is the typical time a star preserves the infection. To guarantee that infections persist for a time length $t_n$ with sufficiently high probability we use \cite[Lemma 4.10]{valesin2024contact} and choose $\mathsf d$ large enough, where the constant $c$ is choose as in the cited lemma. The rest of the proof can be done analogously.    
\end{proof}

\subsubsection{Preliminaries for the exploration algorithm}\label{ss_prelim_algo}

We now turn to proving Proposition~\ref{prop:suitable-finite-subgraph}, the embedding of the star line graph into the bipartite configuration model. 
For the rest of this section, we fix degree sequences~$(\mathbf y^{(n)})_n$ and~$(\mathbf z^{(n)})_n$ (for bottom and top, respectively) and numbers~$\mathsf d$,~$\mathsf r$ such that~$\{(\mathbf y^{(n)}, \mathbf z^{(n)})_n, \mathsf d, \mathsf r\}$ satisfy Assumption~\ref{assumption_new}.
The next lemma ensures that a branching processes that will arise from explorations of the graph can be coupled to the limiting alternating branching tree, which already satisfies the conditions due to Assumption~\ref{assumption_new}.
\begin{lemma}\label{lem_new_bgw} 
	There exists a positive integer~$\mathsf N$ such that the following holds for~$n$ large enough. 

    Let~$Y_n \sim \mathbf p^{(\mathscr y,n)}$ and~$Z_n \sim \mathbf p^{(\mathscr z,n)}$. Let $\widetilde Y_n$ and $\widetilde Z_n$ be the down-shifted size-biased version of $Y_n$ and $Z_n$.
    
    Let~$\mathcal T_n$ be distributed as either~$\mathrm{ABT}(\mathsf d-1,\widetilde{Y}_n,\widetilde{Z}_n)$ or~$\mathrm{ABT}(\mathsf d,\widetilde{Y}_n,\widetilde{Z}_n)$
    
        Let~$\mathcal A'$ be the event that in a breadth-first search of~$\mathcal T_n$ for vertices of degree~$\mathsf d$ in even generations, there are at least two instances when the search succeeds within distance~$\mathsf r$ of the root, \emph{and moreover}, these occur among the first~$\mathsf N$ vertices revealed by the exploration.
        We then have~$\mathbb P(\mathcal T_n \in \mathcal A') \ge 5/8$.
\end{lemma}
Note that the statement of the above lemma differs from Assumption~\ref{assumption_new} in two ways, apart from the minor point that we also allow for degree~{$\mathsf d-1$} at the root. First, the offspring distributions of the alternating branching trees are made from the (size-biased, shifted) $n$-dependent empirical laws, rather than the limiting laws. Second, we modify the event~$\mathcal A$ by adding the restriction that the two successes happen within the first~$\mathsf N$ vertices explored.

\begin{proof}
    We only treat~$\mathcal T_n \sim \mathrm{ABT}(\mathsf d-1,\widetilde{Y}_n,\widetilde{Z}_n)$, since the other case,~$\mathcal T_n \sim \mathrm{ABT}(\mathsf d,\widetilde{Y}_n,\widetilde{Z}_n)$, then follows from a simple monotonicity argument.

    Recall the event~$\cA$ defined in Assumption~\ref{assumption_new}. This is the event that, in a breadth-first search of~$\mathcal T$ for vertices of degree~$\mathsf d$ in even generations, we have two instances of success within distance~$\mathsf r$ of the root. Note that this event only depends on the tree up to depth~$\mathsf r$.

    Let~$\mathcal T \sim \mathrm{ABT}(\mathsf d-1, \widetilde Y, \widetilde Z)$, where~$\widetilde Y$ and~$\widetilde Z$ correspond to the limiting distributions, as in Assumption~\ref{assumption_new}.
Cutting~$\mathcal T_n$ and~$\mathcal T$ after generation~$\mathsf r$, we obtain trees which we denote by~$\mathcal T_n'$ and~$\mathcal T'$, respectively.

    Suppose that~$T$ is a deterministic tree of depth~$\mathsf r$, and denote by~$\zeta_w(T)$ the number of children of~$w\in T$. Then
    \begin{equation}\label{eq:distribution-of-generation-truncated-tree}
        \IP(\cT_n'=T)=\prod_{w\in T}\IP(\zeta_w(\cT_n')=\zeta_w(T)).
    \end{equation}
Assumption~\ref{asmp:convergence_YZ} (which is in force, since it is included in Assumption~\ref{assumption_new}) implies that~$\widetilde Y_n \stackrel{\mathrm{d}}{\longrightarrow}\widetilde Y$ and~$\widetilde Z_n \stackrel{\mathrm{d}}{\longrightarrow} \widetilde Z$ as~$n\to \infty$, and thus the generating functions of both~$\widetilde Y_n$ and~$\widetilde Z_n$ converge uniformly. Together with the fact that, for a given~$T$, the product in~\eqref{eq:distribution-of-generation-truncated-tree} has finitely many factors and~$r$ is fixed, this yields that the total variation distance between the distribution of the truncated trees 
    \begin{equation*}
        d_{TV}(\cT_n', \cT')\to 0 \quad \text{ as } n\to \infty.
    \end{equation*}
    Now let~$\mathsf N\in\IN$ be arbitrary but fixed, to be chosen in the next step. The convergence in total variation  and the definition of the event $\mathcal A'$ in the lemma statement implies that, for every~$\varepsilon>0$, and for~$n$ large enough we have \begin{equation}\label{eq:lower-bound-for-the-prelimit-tree}
    \begin{aligned}
        \IP\big(\cT_n\in \cA'\big)&\geq \IP\big(\{ \cT_n\in \cA\}\cap\{|\cT_n'|<\mathsf N\}\big)\\
        &\geq \IP\big( \{ \cT\in \cA\}\cap\{|\cT'|<\mathsf N\}\big)-\varepsilon\\
        &\geq \IP\big(\cT\in \cA\big)-\IP\big(|\cT'|\geq \mathsf N\big)-\varepsilon,
    \end{aligned}
    \end{equation}
    where we used that~$\cA$ and~$\cA'$ only depend on the tree up to depth~$\mathsf r$. Next, since~$|\cT'|$ is almost surely finite, for every~$\varepsilon>0$ we can choose~$\mathsf N$ large enough that
    \begin{equation*}
        \IP\big(|\cT'|\geq \mathsf N \big)< \varepsilon.
    \end{equation*}
    Finally, choosing~$\varepsilon=\frac{1}{16}$, it follows from the previous inequality,~\eqref{eq:lower-bound-for-the-prelimit-tree} and Assumption~\ref{assumption_new} that
    \begin{equation*}
        \IP\big(\cT_n\in \cA'\big)\geq \frac{5}{8}. \qedhere
    \end{equation*}
\end{proof}

We will also need the following lemma, which will allow us to say that the empirical degree distributions~$\mathbf p^{(\mathscr y,n)}$ and~$\mathbf p^{(\mathscr z,n)}$ are not too impacted when we delete a few vertices from the graph (together with their half-edges).

\begin{lemma}\label{lem:admissible_empirical}
For every~$\epsilon > 0$, there exists~$\kappa > 0$ such that, for~$n$ large enough, we have
\begin{equation*}
    \sum_{j \in A} y^{(n)}_j \le \epsilon n \text{ for every }A \subset \mathcal V^{\mathscr y}_n \text{ with }|A| \le \kappa n 
\end{equation*}
and 
\begin{equation*}
    \sum_{j \in A} z^{(n)}_j \le \epsilon n \text{ for every }A \subset \mathcal V^{\mathscr z}_n \text{ with }|A| \le \kappa n.
\end{equation*}
\end{lemma}

In order to prove this, we will need the following auxiliary result.

\begin{lemma}\label{lem:admissible}
Let~$\mathbf d = (d_1,\ldots, d_\ell)$ be a sequence of positive integers, and let~$q>1$. Then, for every~$\alpha \in (0, 1- \frac{1}{q})$, 
 we have that
\[
\sum_{j \in A} d_j \le \left( \frac{|A|}{\ell}\right)^\alpha \cdot \ell \quad \text{ for all }A \subset [\ell] \text{ such that } \frac{|A|}{\ell} \le  \left(\tfrac{1}{\ell} \|\mathbf d \|_q^q \right)^{-\frac{1}{q-1-\alpha q}}
\]
(where we write~$\|\mathbf d\|_q := (\sum_{j\in[\ell]} d_j^q)^{1/q}$).
\end{lemma}
\begin{proof}
Using H\"older's inequality,
\[\sum_{j\in A}d_j\le \| \mathbf{d}\|_q \cdot |A|^{1-\frac{1}{q}}  = \ell^{-\frac{1}{q}}\cdot \|\mathbf{d}\|_q \cdot \left(\frac{|A|}{\ell}\right)^{\frac{q-1-\alpha q}{q}}\cdot \left(\frac{|A|}{\ell}\right)^\alpha \cdot \ell.\]
The right-hand side is smaller than~$(|A|/\ell)^\alpha \cdot \ell$ if 
\[
\ell^{-\frac{1}{q}} \cdot \|\mathbf{d}\|_q \cdot \left(\frac{|A|}{\ell}\right)^{\frac{q-1-\alpha q}{q}} \le 1,
\]
or equivalently, if~$|A|/\ell \le (\tfrac{1}{\ell}\|\mathbf d\|_q^q)^{-\frac{1}{q-1-\alpha q}}$.
\end{proof}

\begin{proof}[Proof of Lemma~\ref{lem:admissible_empirical}]
Let~$\eta \in (1,\infty)$ be the exponent of Assumption~\ref{asmp:convergence_YZ}, and fix~$\alpha  \in (0,1-\tfrac{1}{\eta})$.
    
    Recall that our sequences~$\mathbf y^{(n)}$ and~$\mathbf z^{(n)}$ are given by
    \[
    \mathbf y^{(n)}=(y^{(n)}_x: x \in \mathcal V^{\mathscr y}_n) \qquad \text{and}\qquad \mathbf z^{(n)}=(z^{(n)}_x: x \in \mathcal V^{\mathscr z}_n).
    \]
    Letting~$Y_n \sim \mathbf p^{(\mathscr y,n)}=(p^{(\mathscr y,n)}_\ell)_\ell$, note that
    \begin{align*}
        \mathbb E[Y_n^\eta] =\sum_\ell \ell^\eta \cdot p^{(\mathscr y,n)}_\ell &= \sum_\ell \ell^\eta \cdot \frac{|\{x \in \mathcal V^{\mathscr y}_n: y^{(n)}_x = \ell\}|}{|\mathcal V^{\mathscr y}_n|}\\
        &= \frac{1}{|\mathcal V^{\mathscr y}_n|} \cdot \sum_\ell \sum_{x \in \mathcal V^{\mathscr y}_n} \ell^\eta \cdot \mathds{1}\{y^{(n)}_x = \ell\} = \frac{1}{|\mathcal V^{\mathscr y}_n|} \cdot  \sum_{x \in \mathcal V^{\mathscr y}_n} (y^{(n)}_x)^\eta = \frac{\|\mathbf y^{(n)} \|^\eta_\eta}{|\mathcal V^{\mathscr y}_n|}.
    \end{align*}
Hence, Lemma~\ref{lem:admissible} implies that
    \begin{equation}\label{eq_with_alpha_y}
        \sum_{x \in A}y^{(n)}_x \le \left( \frac{|A|}{|\mathcal V^{\mathscr y}_n|}\right)^\alpha \cdot |\mathcal V^{\mathscr y}_n| \quad \text{ for all } A \subset \mathcal V^{\mathscr y}_n \text{ such that }\frac{|A|}{|\mathcal V^{\mathscr y}_n|} \le \left(  \mathbb E[Y^\eta_n] \right)^{-\frac{1}{\eta - 1 - \alpha \eta}}.
    \end{equation}
    Similarly, letting~$Z_n \sim \mathbf p^{(\mathscr z,n)}$, we have
    \begin{equation*}
        \sum_{x \in A}z^{(n)}_x \le \left( \frac{|A|}{|\mathcal V^{\mathscr z}_n|}\right)^\alpha \cdot |\mathcal V^{\mathscr z}_n| \quad \text{ for all } A \subset \mathcal V^{\mathscr z}_n \text{ such that }\frac{|A|}{|\mathcal V^{\mathscr z}_n|} \le \left(  \mathbb E[Z^\eta_n] \right)^{-\frac{1}{\eta - 1 - \alpha \eta}}.
    \end{equation*}
We now define
\begin{equation}\label{eq_def_of_kappa}
\kappa := \min \left\{\epsilon^{1/\alpha},\; (2\mathbb E[Y^\eta])^{-\frac{1}{\eta-1-\alpha \eta}},\; (2\mathbb E[Z^\eta])^{-\frac{1}{\eta-1-\alpha \eta}}\right\}.
\end{equation}
By Assumption~\ref{asmp:convergence_YZ}, for~$n$ large enough we have
\begin{equation}\label{eq_twice_etas}
\mathbb E[Y^\eta_n] \le 2 \mathbb E[Y^\eta],\qquad  \mathbb E[Z^\eta_n] \le 2 \mathbb E[Z^\eta].
    \end{equation}
Assume that~$n$ is large enough that this holds. Let~$A \subset \mathcal V^{\mathscr y}_n$ be a set satisfying~$|A| \le \kappa |\mathcal V^{\mathscr y}_n|$. Then, we have
\[\frac{|A|}{|\mathcal V^{\mathscr y}_n|} \le \kappa \stackrel{\eqref{eq_def_of_kappa}}{\le} (2 \mathbb E[Y^\eta])^{-\frac{1}{\eta-1-\alpha \eta}} \stackrel{\eqref{eq_twice_etas}}{\le} (\mathbb E[Y^\eta_n])^{-\frac{1}{\eta-1-\alpha \eta}},\] 
so
\[\sum_{x \in A}y^{(n)}_x \stackrel{\eqref{eq_with_alpha_y}}{\le} \left(\frac{|A|}{|\mathcal V^{\mathscr y}_n|}\right)^\alpha \cdot |\mathcal V^{\mathscr y}_n|\le \kappa^\alpha |\mathcal V^{\mathscr y}_n| \stackrel{\eqref{eq_def_of_kappa}}{\le} \epsilon |\mathcal V^{\mathscr y}_n|.\]
Similarly, if~$A \subset \mathcal V^{\mathscr z}_n$ satisfies~$|A| \le \kappa |\mathcal V^{\mathscr z}_n|$, then~$\sum_{x \in A}z^{(n)}_x \le \epsilon |\mathcal V^{\mathscr z}_n|$.
\end{proof}

\subsubsection{The exploration algorithm}
\label{ss_algo}
Recall that in our current set-up, the set of vertices of~$\mathrm{BC}_n(\mathbf y, \mathbf z)$ is split into bottom vertices (with degree sequence~$\mathbf y^{(n)}$) and top vertices (with degree sequence~$\mathbf z^{(n)}$).

In this section, we present an algorithm to find a subgraph of~$\mathrm{BC}_n(\mathbf y, \mathbf z)$ which is a star-line graph (as in Definition~\ref{def_star_line}) with degree parameter~$\mathsf d$ and distance parameter~$\mathsf r$, with the additional property that all the star centers are bottom vertices.

With this in mind, bottom vertices with degree equal to~$\mathsf d$ will be particularly important; we call them \emph{candidates}.

The algorithm relies on constructing the bipartite configuration model by repeatedly carrying out the following procedure:
\begin{itemize}
	\item \emph{first}, electing a not-yet-matched half-edge attached to a vertex on either side (bottom or top), and
	\item \emph{second}, matching it to a half-edge chosen uniformly at random among the not-yet-matched half-edges on the other side.
\end{itemize}
One iteration of the above procedure (electing and matching) is called one \emph{step} of the algorithm. We only have freedom of choice in the electing part, as the matching happens randomly. Since we are searching for a desirable subgraph, we try to do the electing in a way that works in our favour.

A \emph{bottom-to-top} step is one where the elected half-edge is attached to a bottom vertex, and a \emph{top-to-bottom} step is one where the elected half-edge is attached to a top vertex.

Steps are grouped in \emph{rounds}. Roughly speaking, in a round, we take a candidate vertex (called the \emph{pivot vertex} of the round) and explore its vicinity in a breadth-first fashion, in an attempt to find two new candidates within distance~$\mathsf r$ of the pivot. When the first candidate is found, we do not include its children in the round exploration. There are some failure conditions that can interrupt a round (notably, finding a connection to a previously seen vertex); if none of them occurs and a second candidate is found, then the round is a success and stops. If all goes well, the exploration within a round will be coupled to an alternating branching tree as described in Lemma~\ref{lem_new_bgw}; this will be used to guarantee that the probability of success is not too small.

\paragraph{Input:} The input of the algorithm consists of:
\begin{itemize}
	\item the bipartite degree sequences~$\mathbf y=(\mathbf{y}^{(n)})_n$ and~$(\mathbf{z}^{(n)})_n$ and the integer parameters~$\mathsf d$ and~$\mathsf r$ satisfying Assumption~\ref{assumption_new};
	\item the integer~$\mathsf N$ from Lemma~\ref{lem_new_bgw}, and
	\item a small constant~$\delta > 0$.
\end{itemize}
The choice of~$\delta$ will be made later, when we analyze the algorithm. For now, we impose that
\begin{equation}\label{eq:delta-small-condition}
	\delta < \tfrac{1}{6}\mathbb P(Y=\mathsf d) \cdot \lim_{n \to \infty} \frac{|\mathcal V^{\mathscr y}_n|}{n}
\end{equation}
(recall that the above limit was assumed to exist in~\eqref{eq_divergence_Vs}).
We also assume that~$n$ is large enough that
\[
|\{x \in \mathcal V^{\mathscr y}_n: y_x^{(n)}=\mathsf d\}| \ge \tfrac12 \mathbb P(Y=\mathsf d)\cdot |\mathcal V^{\mathscr y}_n|,
\]
which is possible because~$\mathbf p^{(\mathscr y,n)}$ converges to the distribution of~$Y$. Together with the condition on~$\delta$ in \eqref{eq:delta-small-condition}, this gives
\begin{equation}
\label{eq_number_candidates}
    |\{x \in \mathcal V^{\mathscr y}_n: y_x^{(n)}=\mathsf d\}| \ge 3\delta n.
\end{equation}
This condition guarantees that the algorithm deterministically cannot run out of candidates to explore; we will come back to this later.

\paragraph{Order:} When we are free to pick a vertex or a half-edge among many, we take the smallest one according to the following order. Since the vertex sets~$\mathcal V^{\mathscr y}_n$ and~$\mathcal V^{\mathscr z}_n$ are~$\{1,\ldots,|\mathcal V^{\mathscr y}_n| \}$ and~$\{1,\ldots,|\mathcal V^{\mathscr z}_n| \}$, respectively, we endow them with the usual order of the integers. We endow the set of half-edges attached to each vertex with an arbitrary order. Given two half-edges~$h,h'$ attached to vertices on the same side (top or bottom), we declare that~$h < h'$ either if the vertex of~$h$ is smaller than the vertex of~$h'$, or if they are attached to the same vertex and~$h < h'$ in the aforementioned arbitrary order.

\paragraph{Queues and lists:} For bookkeeping, we maintain a few queues and lists in the algorithm. The first two of them are \emph{global}, meaning that they are not subordinate to a specific round. The third one is \emph{local}, meaning that it is created together with a round (and ceases to exist at the end of the round).
\begin{itemize}
\item The \emph{candidate queue}, denoted by~$\mathbf{QCand}$, contains candidates (bottom vertices with degree $\mathsf{d}$) that are waiting to be pivot vertices of future rounds.
\item The \emph{list of seen half-edges}, denoted by~$\mathbf{Sh}$, contains all half-edges that have been involved in matchings so far (either by being elected or by being chosen at random).
\item The \emph{queue of half-edges of the round}, denoted by~$\mathbf{Qh}$, contains the half-edges that may be elected in an upcoming step within the current round. We note that if a new candidate is found during the round, then we do not include its half-edges in this queue; rather, we only include the candidate in the candidate queue~$\mathbf{QCand}$ (and this is done at the end of the round).
\end{itemize}

We are now ready to present the algorithm.

\begin{algorithm}[H]
\caption*{\textbf{The algorithm}}
\begin{algorithmic}[1]
\State Initialize~$\mathbf{QCand}$ and~$\mathbf{Sh}$ as empty
\State Initialize the round counter~$i$ as 0
\While{$i < \lfloor \delta n \rfloor$}
\State \Comment{Round starts now}
    \If{$\mathbf{QCand}=\emptyset$}
        \State Let~$v$ be the smallest candidate that currently has no half-edge of~$\mathbf{Sh}$ attached to it. Add~$v$ to~$\mathbf{QCand}$: $\mathbf{QCand} \leftarrow \mathbf{QCand} \cup \{v\}$
    \EndIf
    \State Let~$\hat{v}$ be the first vertex in~$\mathbf{QCand}$. Call~$\hat{v}$ the \emph{round pivot}. Remove it from~$\mathbf{QCand}$
    
        \State Initialize~$\mathbf{Qh}$ as the set of unseen half-edges attached to~$\hat{v}$
        \State Initialize the counter~$j$ of steps taken in this round as 0
        \State Initialize the counter~$\mathbf{foundcounter}$ of new candidates found in this round as 0
        \State Initialize the indicator variable~$\mathbf{break}$ as 0; this will indicate when the round must be stopped due to some undesirable condition being fulfilled
        \While{$j < \mathsf N$ $\mathbf{and}$ $\mathbf{foundcounter}<2$ $\mathbf{and}$ $\mathbf{break}=0$}
                \State \Comment{Steps start now}
            \State Let~$h$ be the first half-edge of~$\mathbf{Qh}$. Remove it from~$\mathbf{Qh}$
            \State Choose a half-edge~$h'$ uniformly at random from the set of half-edges that are not in~$\mathbf{Sh}$ on the side opposite to that of~$h$ (that is: if~$h$ is attached to a bottom vertex, then~$h'$ is attached to a top vertex, and vice-versa). Let~$\mathcal X$ be the vertex to which~$h'$ belongs
            \State Match~$h$ and~$h'$, and add both of them to~$\mathbf{Sh}$
            \State Set $j \leftarrow j+1$
            \State Let~$\mathcal N$ be the number of half-edges of~$\mathcal X$ belonging to~$\mathbf{Sh}$
            \State \Comment{If~$\mathcal N > 1$, then the newly made matching created a loop}
            \If{$\mathcal N > 1$}
            \State{Set~$\mathbf{break} \leftarrow 1$}
            \State \Comment{If the loop involved a candidate that was in the candidate queue, we consider that candidate spoiled and remove it from the queue}
            \If{$\mathcal X \in \mathbf{QCand}$}
\State Remove~$\mathcal X$ from~$\mathbf{QCand}$
            \EndIf
                        \ElsIf{$\mathcal N = 1$ $\mathbf{and}$ $\mathcal X$ is not a candidate}
                \State Add to~$\mathbf{Qh}$ all the half-edges of~$\mathcal X$ that are different from~$h'$ (and hence are not in~$\mathbf{Sh}$)
            \ElsIf{$\mathcal N = 1$ $\mathbf{and}$ $\mathcal X$ is a candidate}
                \State Set $\mathbf{foundcounter} \leftarrow \mathbf{foundcounter} + 1$
            \EndIf
        \If{$\mathbf{Qh}=\emptyset$ $\mathbf{or}$, in the graph distance induced by the set of edges found in the present round,~$\mathbf{Qh}$ has at least one half-edge belonging to a vertex at distance more than~$\mathsf r$ from~$\hat v$}
\State Set $\mathbf{break} \leftarrow 1$
            \EndIf
        \EndWhile
    \State\Comment{We now check whether the round was a success; if so, we add the new candidates to the candidate queue}
    \If{$\mathbf{foundcounter}=2$ $\mathbf{and}$ $\mathbf{break}=0$}
\State Add the two candidates that were found in the round to the end of~$\mathbf{QCand}$
        \EndIf
    \State Set $i \leftarrow i+1$
\EndWhile
\end{algorithmic}
\end{algorithm}

When a round starts, if~$\mathbf{QCand}=\emptyset$, the algorithm picks a fresh candidate to be the pivot (that is, a candidate with no half-edge currently in~$\mathbf{Sh}$). Condition~\eqref{eq_number_candidates} ensures that the set of fresh candidates is never empty. Indeed, each round touches half-edges of at most three fresh candidates (namely, the pivot and possibly the two new candidates that are found), and the number of rounds is~$\lfloor \delta n \rfloor$.

A round is considered successful if it ends with~$\mathbf{foundcounter}=2$ and~$\mathbf{break}=0$, meaning that two new candidates have been found (and put in~$\mathbf{QCand}$), and no error has occurred.
There are several ways in which a round can fail (i.e.\ not be a success):
\begin{itemize}

    \item \textbf{failure due to a loop:} a matching is made to a half-edge belonging to a vertex~$\mathcal X$ having a half-edge that had already been involved in a matching in a previous step (either in the present round or in a previous one);
    \item \textbf{failure due to a dead end:} at a moment when two new candidates have not yet been found, the round runs out of half-edges to explore, that is,~$\mathbf{Qh}=\emptyset$;
    \item \textbf{failure due to reaching too far:} at a moment when two new candidates have not yet been found, we reach distance larger than~$\mathsf r$ from the pivot (in the graph distance induced by the set of edges that have been found in the present round);
    \item \textbf{failure due to exhaustion of exploration:} the round performs all~$\mathsf N$ steps, but does not discover two new candidates.
\end{itemize}
We denote by~$F^{\mathrm{loop}}_i$,~$F^{\mathrm{dead}}_i$,~$F^{\mathrm{far}}_i$, and~$F^{\mathrm{exh}}_i$ the events that round~$i$ fails due to each of the errors in the above list, in order. These events are not disjoint (for instance, the same matching could cause a loop and also exhaust the half-edges to explore), but this will not be an issue.

If~$F_i^{\mathrm{loop}}$ occurs, there is the further possibility that the loop involved a candidate that was sitting in~$\mathbf{QCand}$; we call this a \textbf{catastrophe}. We let~$F^\mathrm{cat}_i \subset F^\mathrm{loop}_i$ denote the event that a catastrophe occurs in round~$i$.

\subsubsection{Analysis of the exploration algorithm}
\label{ss_algo_analysis}
We index rounds by~$i \in \{0,\ldots, \lfloor \delta n\rfloor -1 \}$ (the number of rounds is deterministic, equal to~$\lfloor \delta n \rfloor$). Let~$S_i$ be the number of steps in round~$i$. Note that~$S_i \le \mathsf N$. We index steps within the round by~$j \in \{0,\ldots,S_i-1\}$. Hence, a double index of the form~$(i,j)$ refers to the step with index~$j$, counting from 0 at the beginning of round~$i$.

Recalling that in each step we match an \emph{elected} half-edge with a \emph{random} half-edge, we denote the elected half-edge by~$h_{(i,j)}$ and the randomly chosen half-edge by~$h_{(i,j)}'$; we denote by~$\mathcal X_{(i,j)}$ the vertex to which~$h_{(i,j)}'$ is attached. We also define a random variable~$\mathrm{Direction}_{(i,j)}$ taking values in the set~$\{ \mathrm{upwards}, \mathrm{downwards} \}$ as follows:
\[
\mathrm{Direction}_{(i,j)}=\begin{cases}
    \mathrm{upwards}&\text{if } h_{(i,j)}' \text{ and }\mathcal X_{(i,j)} \text{ are on the top side};\\
    \mathrm{downwards}&\text{otherwise.}
\end{cases}
\]

Let~$\mathcal F_{i-}$ be the~$\sigma$-algebra generated by all the matchings made by the algorithm before round~$i$ (with~$\mathcal F_{0-}$ being the trivial~$\sigma$-algebra~$\{\emptyset, \Omega\}$). We will prove the following.
\begin{lemma}\label{lem_errors}
For every~$\epsilon > 0$, if~$\delta$ is small enough (depending on~$\epsilon$, $\mathsf d$ and~$\mathsf r$
) and~$n$ is large enough (depending on~$\epsilon$, $\mathsf d$, $\mathsf r$
 and $\delta$), we have
\begin{equation}\label{eq_loops}
\mathbb P\big( F^\mathrm{loop}_i \mid \mathcal F_{i-}\big) < \epsilon
\end{equation}
and
\begin{equation}
    \label{eq_other_errors}
    \mathbb P \big(F^\mathrm{dead}_i \cup F^\mathrm{far}_i \cup F^\mathrm{exh}_i \cup F^\mathrm{loop}_i \mid \mathcal F_{i-} \big) < \frac{3}{8}+2\epsilon
\end{equation}
for all~$i \in \{0,\ldots, \lfloor \delta n\rfloor -1\}$.
\end{lemma}

As a consequence, we immediately obtain:
\begin{corollary}\label{cor_success}
    If~$\delta$ is small enough (depending on~$\mathsf d$ and~$\mathsf r$
    ) and~$n$ is large enough (depending on~$\mathsf d$, $\mathsf r$
    and~$\delta$), then
    \[
    \mathbb P \big(\text{round $i$ is a success} \mid \mathcal F_{i-}\big) > \frac{9}{16}\quad \text{for all }i \in \{0,\ldots, \lfloor \delta n\rfloor -1\}.
    \]
\end{corollary}
\begin{proof}[Proof of Lemma~\ref{lem_errors},~\eqref{eq_loops}]
A loop is created in the algorithm when we perform a matching to a half-edge belonging to a vertex having at least one half-edge that had been seen previously (that is, that belonged to the list~$\mathbf{Sh}$). We therefore need to bound the sum of the degrees of the vertices with half-edges in~$\mathbf{Sh}$, at an arbitrary point during the algorithm.

In round~$i$, assume that the matching for step~$j$ is about to be performed, and condition on the whole algorithm up to this point. We assume that~$\mathrm{Direction}_{(i,j)}=\mathrm{downwards}$, as the other case is treated similarly.
By Assumption~\ref{asmp:convergence_YZ}, $(\sum_{x \in \mathcal V^{\mathscr y}_n}y_x)/|\mathcal V^{\mathscr y}_n|$ converges to~$\mathbb E[Y]$ as~$n \to \infty$. Take~$\delta < \frac{1}{2\mathsf N} \mathbb E[Y]$ and~$n$ large.
Since~$\mathsf N \delta n$ is an upper bound for the number of steps (matchings) performed by the whole algorithm, 
the total number of half-edges on the bottom side still available for matchings in step~$(i,j)$ is then at least
\[
\sum_{x \in \mathcal V^{\mathscr y}_n} y_x - \mathsf N \delta n \ge \frac12 |\mathcal V^{\mathscr y}_n|\cdot \mathbb E[Y].
\]
as we have assumed $\delta < \frac{1}{2\mathsf N} \mathbb E[Y]$.
Let~$A$ be the set of bottom vertices having at least one half-edge in~$\mathbf{Sh}$ immediately before step~$(i,j)$. We have~$|A| \le \mathsf N \delta n$. Fix~$\epsilon_0 > 0$, to be chosen later. By Lemma~\ref{lem:admissible_empirical}, if~$\delta$ is small enough, then for all large enough~$n$
\[
\sum_{x \in A} y_x \le \epsilon_0  n.
\]
Now, step~$(i,j)$ produces a loop if and only if~$h_{(i,j)}'$ is attached to a vertex of~$A$. This probability is bounded from above by an upper bound on the number of unmatched half-edges in $A$ over a lower bound on the total unmatched half-edges, giving
\[
\mathbb P(\text{Step $(i,j)$ produces a loop)} \le \frac{\epsilon_0 n}{\frac12 |\mathcal V^{\mathscr y}_n|\cdot \mathbb E[Y]}\le \frac{4 \epsilon_0}{\gamma^{\mathscr y}\mathbb E[Y]},
\]
where where~$\gamma^{(\mathscr y)}:= \lim_{n \to \infty} |\mathcal V^{\mathscr y}_n|/n$, and the last inequality follows by using~\eqref{eq_divergence_Vs} for~$n$  large enough.
By taking~$\epsilon_0 < \epsilon \gamma^{\mathscr y} \mathbb E[Y]/(4\mathsf N)$, we obtain that the probability that step~$(i,j)$ produces a loop is smaller than~$\epsilon/\mathsf N$. Since at most~$\mathsf N$ steps are performed in the round, the probability that a loop occurs in the round is smaller than~$\epsilon$.
\end{proof}

\begin{proof}[Proof of Lemma~\ref{lem_errors},~\eqref{eq_other_errors}]
Fix the round~$i \in \{0,\ldots, \lfloor \delta n \rfloor - 1\}$ and condition on~$\mathcal F_{i-}$.

We will use some auxiliary random variables to sample~$h_{(i,j)}'$, in a way that facilitates a coupling with an alternating branching tree.

We let~$\mathscr H^{\mathscr z}$ be the set of half-edges attached to top vertices, and~$\mathscr H^{\mathscr y}$ be the set of half-edges attached to bottom vertices. We let~$(H^{\mathscr z}_{j,\ell}: j \ge 0, \; \ell \ge 0)$ be a family of half-edges chosen independently and uniformly at random from~$\mathscr H^{\mathscr z}$ (with replacement), and independently let~$(H^{\mathscr y}_{j,\ell}:j \ge 0,\; \ell \ge 0)$ be a family of half-edges chosen independently and uniformly at random from~$\mathscr H^{\mathscr y}$ (with replacement). We emphasize that these are sampled from the whole sets of half-edges, with no regard to what the algorithm has done. The index~$j$ in~$H^{\mathscr z}_{j,\ell}$ and in~$H^{\mathscr y}_{j,\ell}$ will play the role of steps in the round, {while the index $\ell$ is an auxiliary index. Here we also allow for~$j \ge S_i$. }

For each~$j \in \{0,\ldots, S_i-1\}$, let~$\mathbf{Sh}_{(i,j-)}$ denote the set of half-edges in the list~$\mathbf{Sh}$ immediately before matching~$j$ is performed. Also let
\[
L_j:=\begin{cases}
    \min\{\ell \ge 0: \; H^{\mathscr z}_{j,\ell} \notin \mathbf{Sh}_{(i,j-)} \}&\text{if } \mathrm{Direction}_{(i,j)} = \mathrm{upwards};\\[.2cm]
    \min\{\ell \ge 0: \; H^{\mathscr y}_{j,\ell} \notin \mathbf{Sh}_{(i,j-)} \}&\text{if } \mathrm{Direction}_{(i,j)} = \mathrm{downwards},
\end{cases}
\]
and then set
\[h'_{(i,j)}=\begin{cases}
H^{\mathscr z}_{j,L_j}&\text{if }\mathrm{Direction}_{(i,j)} = \mathrm{upwards};\\[.2cm]
H^{\mathscr y}_{j,L_j}&\text{if }\mathrm{Direction}_{(i,j)} = \mathrm{downwards}.
\end{cases}\]
In words,~$h'_{(i,j)}$ is the first half-edge of the sequence~$(H^\Box_{j,\ell})_{\ell \ge 0}$ (with~$\Box$ being either `${\mathscr z}$' or `${\mathscr y}$') that is not one of the half-edges already matched by the algorithm before step~$j$. This choice guarantees that~$h'_{(i,j)}$ is uniformly distributed on~$\mathscr H^\Box \setminus \mathbf{Sh}_{(i,j-)}$, as required. {We mention that as the number of seen half-edges is small during the whole algorithm, the typical scenario is that $L_j=0$. It may however occur that we sample $H_{j,0}^{\Box}$ from the seen half-edges, resulting in  $L_j>0$.}

We now define a tree~$\mathcal T_0$ rooted at the pivot of the round.
On the event
\begin{equation}\label{eq_event_exploration}
\{L_j=0 \text{ for all~$j \in \{0,1, \dots S_i-1\}$, and the round does not fail due to a loop}\},
\end{equation}
we let~$\mathcal T_0$ be the tree induced by the set of edges constructed in the round. On the complement of this event, let~$J^*$ be the smallest index for which either~$L_j > 0$ or step~$j\le S_i-1$ causes the round to stop due to a loop. Then, let~$\mathcal T_0$ be the tree formed by the set of edges built in the round before step~$J^*$.

The fact that each matching that produces an edge of~$\mathcal T_0$ uses the first element of the sequence~$(H^\Box_{j,\ell})_{\ell \ge 0}$ implies that~$\mathcal T_0$ is an (incomplete) branching tree in which the offspring is distributed as~$\widetilde{Z}_n$ in odd generations and as~$\widetilde{Y}_n$ in even generations (apart from the root). We say ``incomplete'' because the tree is only built in a breadth-first way for as long as the round runs (or until~$J^*$ is reached). Then, using additional independent randomness, and new dummy vertices that are added apart from those of~$\mathcal V^{\mathscr z}_n \cup \mathcal V^{\mathscr y}_n$, we continue building the branching tree with the same (alternating) offspring distributions (potentially making it infinite); we denote by~$\mathcal T \supset \mathcal T_0$ the tree thus produced.

Recall the definition of the event~$\mathcal A'$ from Lemma~\ref{lem_new_bgw}. This event guarantees that in the alternating branching process $\mathrm{ABT}(d, \widetilde Y_n, \widetilde Z_n)$ the following good event happens: within distance $\mathsf{r}$ of the root we find two other vertices of degree $\mathsf{d}$ within the first $\mathsf{N}$ steps of the exploration. 
On the event defined in~\eqref{eq_event_exploration}, the portion of~$\mathcal T$ that is revealed in the first~$\mathsf N$ steps of the exploration is also contained in~$\mathcal T_0$.

We note that on the event $\mathcal{A}'$ none of the events $F_i^{\mathrm{dead}},F^\mathrm{far}_i,F^\mathrm{exh}_i$ occur as then then the round is a success. Thus, either $J_\star\le S_i-1$ and we cannot perfectly couple to the alternating branching process, or $J_\star$ is not defined, there is also no loop in round $i$, and in that case the tree we build is a sample from the alternating branching process $\mathrm{ABP}(d, \widetilde Y_n, \widetilde Z_n)$ up to seeing the first $S_i-1$ edges (in a breadth first search order). 
We can thus bound using the complement of \eqref{eq_event_exploration} and the complement of the event $\mathcal A'$ to obtain that
\begin{align*}
    \mathbb P \Big(F^\mathrm{dead}_i \cup F^\mathrm{far}_i \cup F^\mathrm{exh}_i \cup F_i^{\mathrm{loop}} \mid \mathcal F_{i-} \Big) &\le  \mathbb P(\exists j \in \{0,\ldots, S_i-1\}: L_j > 0 \mid \mathcal F_{i-})\\
    &\ +
    \mathbb P(F_i^{\mathrm{loop}} \mid \mathcal F_{i-}) + \mathbb P(\mathcal T \notin \mathcal A' \mid \mathcal F_{i-})\\
    &\le \mathbb P(\exists j \in \{0,\ldots, S_i-1\}: L_j > 0 \mid \mathcal F_{i-}) +\epsilon+\frac{3}{8},
\end{align*}
where we note that we used here the already established bound $\P(F_i^{\mathrm{loop}})\le \epsilon$.
In order to bound the probability on the right-hand side, fix~$j \ge 0$ and assume that matching~$j$ is performed in the round; in addition to~$\mathcal F_{i-}$, condition on all the steps performed in round~$i$ before step~$j$. As in each round we only match at most $\mathsf N$ edges, and there are at most $\delta n$ rounds, the (conditional) probability that~$L_j > 0$ is then at most
\[
\frac{|\mathbf{Sh}_{(i,j-)}\cap \mathscr H^{\mathscr y}|}{|\mathscr H^{\mathscr y}| } \le \frac{\mathsf N \delta n }{\tfrac12 |\mathcal V^{\mathscr y}_n| \cdot  \mathbb E[Y]}
\]
as soon as~$n$ is large enough that~$|\mathscr H^{\mathscr y}| \ge \tfrac12 |\mathcal V^{\mathscr y}_n| \cdot \mathbb E[{Y}]$.

By taking~$\delta$ small enough, the right-hand side is smaller than~$\epsilon/\mathsf N$; a union bound over the steps in the round concludes the proof.
\end{proof}

The next lemma controls the number of catastrophes that occur in the algorithm. Recall that we write $F_i^{\mathrm{cat}}$ and we call a round a catastrophe when a loop occurs in round $i$ and this loop is matched to another already found candidate vertex, that is, a vertex in $\mathbf{QCand}$.
\begin{lemma}\label{lem_cats}
    For every~$\epsilon > 0$, if~$\delta$ is small enough (depending on~$\epsilon$ and~$\mathsf d$ and~$\mathsf r$
    ), then
\[
\lim_{n \to \infty} \mathbb P\left(\sum_{i=0}^{\lfloor \delta n \rfloor -1}  \mathds{1}(F^{\mathrm{cat}}_i) > \epsilon \delta n \right) =0.
\]
\end{lemma}
\begin{proof}
Instead of directly handling the indicator~$\mathds{1}(F^{\mathrm{cat}}_i)$ that a catastrophe occurs in round~$i$, we will sum over candidate vertices and check whether a candidate vertex~$x$ is involved in a catastrophe. This can happen only if~$x$ is matched twice before it can act as a pivot: first when~$x$ is found, and then when it is accidentally involved in a loop.

{ As candidate vertices are all bottom vertices,  a loop of the latter type can only occur in a \emph{downward} matching.}
Recall that we regard a matching as \emph{downwards} if the elected half-edge is on the top side and the half-edge chosen at random is on the bottom side.

We will use some auxiliary randomness to govern the downwards matchings of the algorithm, in a way that facilitates a convenient coupling. This is similar in flavour to what was done in the proof of~\eqref{eq_other_errors} above, but with a different indexing.

Let~$\mathscr H^{\mathscr y}$ be the set of half-edges attached to bottom vertices. Let~$\zeta_0, \zeta_1,\ldots$ be half-edges chosen independently and uniformly at random from~$\mathscr H^{\mathscr y}$ (with replacement).

We enumerate the downward steps of the algorithm in the order in which they occur, with the index~$b\in \{0,1,\ldots\}$. Note that we are using a single consecutive indexing here, from the beginning to the end of the algorithm, across all rounds.

Let~$\mathbf{Sh}_{b-}$ be the set of half-edges in the list~$\mathbf{Sh}$ immediately before downward step~$b$ is performed.

We now set~$L_0:=0$ and, letting~$B$ be the index of the last downward step performed when the whole algorithm finishes running, we recursively define
\begin{equation}\label{eq:Lb-new-def}
L_{b}:=\inf\{\ell > L_{b-1}:\; \zeta_\ell \notin \mathbf{Sh}_{b-}\},\quad b \in \{1,\ldots, B\}.
\end{equation}
{Again, typically $L_b=L_{b-1}+1$ as the number of seen half-edges in the whole algorithm is small. However, it may occur that we accidentally re-sample a seen half-edge, in which case $L_b> L_{b-1}+1$}. 
We then adopt the sampling rule that, whenever the downward step with index~$b$ is performed, we let~$\zeta_{L_b}$ be the half-edge chosen to match the elected half-edge for that step (that is, recalling our convention that the elected half-edge is denoted by~$h$ and the randomly chosen half-edge is denoted by~$h'$,~$\zeta_{L_b}$ plays the role of~$h'$). This sampling rule has the required property: whenever a half-edge is chosen, it is taken uniformly at random from the set of half-edges that have not been seen so far. We will not need an analogous sampling procedure for the upwards steps.

 Let~$\mathscr C$ denote the set of all candidate vertices (which are bottom vertices), and with $\zeta_\ell$ a uniform variable from all bottom half-edges $\mathscr H^{\mathscr y}$, define
\begin{equation}\label{eq:zeta-x-ell-def}
\mathcal J_{x,\ell}:=\mathds{1}(\zeta_\ell \text{ is attached to vertex $x$}),\quad x \in \mathscr C,\; \ell \ge 0.
\end{equation}
We note that
\[
\sum_{i=0}^{\lfloor \delta n \rfloor -1} \mathds 1(F^{\mathrm{cat}}_i) \le \sum_{x \in \mathscr C} \mathds{1}\left(\sum_{\ell=0}^{L_B} \mathcal J_{x,\ell} \ge 2 \right),
\]
since each catastrophe can be mapped injectively to a candidate that appears at least twice in~$(\zeta_\ell)_{0\le \ell \le L_B}$. We can then bound
\[
\mathbb P\left(\sum_{i=0}^{\lfloor \delta n \rfloor -1}  \mathds{1}(F^{\mathrm{cat}}_i) > \epsilon \delta n \right) \le \mathbb P \left(\sum_{x \in \mathscr C} \mathds{1}\left( \sum_{\ell=0}^{L_B} \mathcal J_{x,\ell}\ge 2\right) > \epsilon \delta n \right).
\]
The probability on the right-hand side above is easier to handle, but we want to replace the sum over~$\ell$ up to the random $L_B$ with a sum with a deterministic number of terms. We hence bound
\begin{equation}\label{eq_two_probabilities}
\begin{aligned}
\mathbb P \left(\sum_{x \in \mathscr C} \mathds{1}\left( \sum_{\ell=0}^{L_B} \mathcal J_{x,\ell}\ge 2\right) > \epsilon \delta n \right) &\le \mathbb P(L_B > 3 \delta \mathsf N n) \\
&\ + \mathbb P \left(\sum_{x \in \mathscr C} \mathds{1}\left( \sum_{\ell=0}^{\lfloor 3 \delta \mathsf N n\rfloor} \mathcal J_{x,\ell}\ge 2\right) > \epsilon \delta n \right).
\end{aligned}
\end{equation}
Here we have chosen the value~$3 \delta \mathsf N n$ because~$L_B$ is unlikely to exceed it. The intuitive reasoning is as follows. The total number of bottom half-edges that enter~$\mathbf{Sh}$ by the end of the algorithm is at most~$\delta \mathsf N n$; when~$\delta$ is small, this is small compared to the total number of half-edges attached to bottom vertices (which is about a fixed constant times~$n$). Because of this, conditionally on~$\mathbf{Sh}_{b-}$, the distribution of~$L_b-L_{b-1}$ is geometric with a parameter close to 1, so this random variable is likely to equal~1 (we dominate it by a~$\mathrm{Geometric}(1/2)$ for convenience). Hence,~$L_B$ is likely to be smaller than~$3B$. Since we also have~$B<\delta \mathsf N n$, we expect~$L_B$ to be below~$3\delta \mathsf N n$.

We now make the above heuristics precise. We will bound the two probabilities on the right-hand side of~\eqref{eq_two_probabilities} separately.

For each~$b \ge 1$, conditionally on all steps before step~$b$, the distribution of~$L_b-L_{b-1}$ is~$\mathrm{Geometric}\big(\frac{|\mathscr H^{\mathscr y}\setminus \mathbf{Sh}_{b-}|}{|\mathscr H^{\mathscr y}|}\big)$. Since~$|\mathbf{Sh}_{b-} \cap \mathscr H^{\mathscr y}| \le \mathsf N \delta n$, by taking~$\delta$ small we can ensure that~$|\mathbf{Sh}_{b-} \cap \mathscr H^{\mathscr y}|\le \frac{|\mathscr H^{\mathscr y}|}{2}$, so that the distribution of~$L_b-L_{b-1}$ is stochastically dominated by $\mathrm{Geometric}(1/2)$ for all~$b\le B$. This shows that
\[
L_B = \sum_{b=1}^B (L_b-L_{b-1}) \preceq \sum_{a=1}^{\lfloor \delta \mathsf N n\rfloor} T_a,\quad \text{where } (T_a)_a \text{ are i.i.d.} \sim \mathrm{Geometric}(1/2).
\]
In particular, by the law of large numbers,
\[
\mathbb P(L_B > 3\delta \mathsf N n) \xrightarrow{n \to \infty} 0.
\]
We now turn to the second probability on the right-hand side of~\eqref{eq_two_probabilities}. We claim that
\begin{equation}\label{eq_claim_dom}
\mathbb E \left[ \mathcal J_{x,\ell} \mid \big(\mathcal J_{x',\ell'}: (x',\ell') \neq (x,\ell)\big)\right] \le \frac{8\mathsf d}{\mathbb E[Y]\cdot \gamma^{\mathscr y} n} \quad \text{ almost surely, for all }(x,\ell),
\end{equation}
where~$\gamma^{\mathscr y}:=\lim_{n \to \infty} |\mathcal V^{\mathscr y}_n|/n$.
To prove this, first note that by the independence of~$\zeta_\ell$ and~$\zeta_{\ell'}$ for~$\ell\neq\ell'$, and the definition of~$\mathcal J_{x,\ell}$ in~\eqref{eq:zeta-x-ell-def}, we have
\begin{equation*}
\mathbb E \left[ \mathcal J_{x,\ell} \mid (\mathcal J_{x',\ell'}: (x',\ell') \neq (x,\ell)) \right]  =
\mathbb E \left[ \mathcal J_{x,\ell} \mid (\mathcal J_{x',\ell}: x' \neq x) \right].
\end{equation*}
Again by~\eqref{eq:zeta-x-ell-def}, for any fixed~$\ell\ge 0$ it holds that~$\sum_{y\in \mathscr C} \mathcal J_{y,\ell} \in \{0,1\}$. This implies that on the event~$\{\sum_{x':x' \neq x}\mathcal J_{x',\ell} = 1\}$, the half-edge~$\zeta_\ell$ belongs to a vertex other than~$x$, so we have~$\mathbb E \left[ \mathcal J_{x,\ell} \mid (\mathcal J_{x',\ell}: x' \neq x) \right]=0$. On the event~$\{\sum_{x':x' \neq x}\mathcal J_{x',\ell} = 0\}$, meaning that the half-edge~$\zeta_\ell$ does not belong to any of the other candidate vertices, we use the fact that there are numerous non-candidate vertices in the graph. This gives that the conditional expectation is
\begin{align}\nonumber
\frac{\mathbb P(\zeta_\ell \text{ is a half-edge of $x$})}{\mathbb P(\zeta_\ell \text{ is not a half-edge of $x'$, for all $x' \in \mathscr C \setminus \{x\}$})} &= \frac{\mathsf d/|\mathscr H^{\mathscr y}|}{(|\mathscr H^{\mathscr y}|-\mathsf d (|\mathscr C|-1))/|\mathscr H^{\mathscr y}|} \\
&= \frac{\mathsf d}{|\mathscr H^{\mathscr y}|-\mathsf d(|\mathscr C|-1)}.\label{eq_aux_dlambda}
\end{align}
Using Assumption~\ref{asmp:convergence_YZ}, we can take~$n$ large enough that
\[\frac{|\mathscr H^{\mathscr y}|}{|\mathcal V^{\mathscr y}_n|} \ge \frac{\mathbb E[Y]}{2}\qquad \text{and} \qquad \frac{|\mathscr C|}{|\mathcal V^{\mathscr y}_n|} \le 2\mathbb P(Y=\mathsf d).\]
Moreover, by~\eqref{eq_divergence_Vs}, we can assume~$n$ is large enough that
\[
\frac{\gamma^{\mathscr y}}{2} \le \frac{|\mathcal V^{\mathscr y}_n|}{n} \le 2\gamma^{\mathscr y}.
\]
This then gives
\[
|\mathscr H^{\mathscr y}| - \mathsf d (|\mathscr C|-1) \ge |\mathscr H^{\mathscr y}| - \mathsf d |\mathscr C| \ge \left(\tfrac{1}{4}\mathbb E[Y] - 4\mathsf d \mathbb P(Y=\mathsf d)\right) \cdot \gamma^{\mathscr y} n \stackrel{\eqref{eq_extra_asmp_candidate}}{\ge} \tfrac18 \mathbb E[Y] \cdot \gamma^{\mathscr y}n.
\]
This shows that the expression in~\eqref{eq_aux_dlambda} is at most~$8\mathsf d/(\mathbb E[Y] \gamma^{\mathscr y}n)$, completing the proof of~\eqref{eq_claim_dom}.

Now, with~\eqref{eq_claim_dom} at hand, it is standard to see that the family of random variables~$\{\mathcal J_{x,\ell}: x \in \mathscr C,\; \ell \ge 0 \}$ is stochastically dominated by a family~$\{\mathcal J_{x,\ell}': x \in \mathscr C,\; \ell \ge 0 \}$ of independent Bernoulli($8\mathsf d/(\mathbb E[Y]\gamma^{\mathscr y} n)$) random variables.
In particular,
\begin{equation}\label{eq_almost_Zs}
\mathbb P \left(\sum_{x \in \mathscr C} \mathds{1}\left( \sum_{\ell=0}^{\lfloor 3 \mathsf N\delta n\rfloor} \mathcal J_{x,\ell}\ge 2\right) > \epsilon \delta n \right) \le
\mathbb P \left(\sum_{x \in \mathscr C} \mathds{1}\left( \sum_{\ell=0}^{\lfloor 3 \mathsf N\delta n\rfloor} \mathcal J_{x,\ell}'\ge 2\right) > \epsilon \delta n \right).
\end{equation}
For each~$x$, using the bound~$\mathbb P(\mathrm{Bin}(N,p)\ge 2) \le N^2p^2$, we have
\[
\mathbb P \left( \sum_{\ell=0}^{\lfloor 3 \mathsf N\delta n\rfloor} \mathcal J_{x,\ell}'\ge 2\right) \le (\lfloor 3 \mathsf N \delta n \rfloor)^2 \cdot \left(\frac{8\mathsf d}{\mathbb E[Y]\gamma^{\mathscr y}n} \right)^2 \le 576\left(\frac{ \mathsf N \delta \mathsf d}{\mathbb E[Y] \gamma^{\mathscr y}}\right)^2 =: \mathsf p.
\]
Hence, the probability on the right-hand side of~\eqref{eq_almost_Zs} is at most
\begin{equation}\label{eq:chernoff-later}
\mathbb P(\mathrm{Bin}(|\mathscr C|,\mathsf p) > \epsilon \delta n ).
\end{equation}
Now, using the fact that~${|\mathscr C|}/{n} \xrightarrow{n \to \infty} \gamma^{\mathscr y} \cdot \mathbb P(Y=\mathsf d) =:p_Y(\mathsf d)$, for all~$n$ sufficiently large we can bound the above by
\[
\mathbb P(\mathrm{Bin}(\lceil 2 p_Y( \mathsf d ) n \rceil ,\mathsf p) > \epsilon \delta n ).
\]
The expectation of this binomial is
\[
\left \lceil 2 p_Y(\mathsf d) n \right \rceil \cdot  576\left(\frac{ \mathsf N \delta \mathsf d}{\mathbb E[Y] \gamma^{\mathscr y}}\right)^2.
\]
Since~$\delta^2 \ll \delta$ when~$\delta$ is small, we can now choose~$\delta$ small enough that the above is smaller than~$\epsilon \delta n/2$ for all~$n$ sufficiently large. A Chernoff bound for the binomial random variable in~\eqref{eq:chernoff-later} concludes the proof.
\end{proof}

\begin{proposition}\label{prop_many_candidates}
If~$\delta$ is small enough (depending on~$\mathsf d$ and~$\mathsf r$
), then with probability tending to~1 as~$n \to \infty$, there are at least~$\delta n/32$ candidates in~$\mathbf{QCand}$ when the algorithm ends.
\end{proposition}
\begin{proof}
We assume that~$\delta$ is small enough, as required by Corollary~\ref{cor_success} (we will need to reduce it further later).

Let~$X_i$ denote the number of elements of~$\mathbf{QCand}$ at the end of round~$i$. Write~$F_i:=F_i^{\mathrm{loop}} \cup F_i^{\mathrm{dead}} \cup F_i^{\mathrm{far}} \cup F_i^{\mathrm{exh}}$ and let~$\Delta X_i := X_{i}-X_{i-1}$.

Note that:
\begin{itemize}
\item if~$X_{i-1}=0$, the round starts by adding a new element to~$\mathbf{QCand}$; if~$X_{i-1}>0$, this is not done;
\item if the round is a success, two new elements are added to~$\mathbf{QCand}$;
\item the round pivot is always removed from~$\mathbf{QCand}$ (regardless of whether the round is a success or a failure);
\item if the round ends in a catastrophe, we lose one element of~$\mathbf{QCand}$ (apart from the pivot), namely, the candidate that was spoiled due to a loop.
\end{itemize}
This gives
\begin{equation}\label{eq_Delta}\begin{aligned}
\Delta X_i &= \mathds{1}(X_{i-1}=0) + 2 \cdot \mathds{1}(F^\complement_i) -1 - \mathds{1}(F^\mathrm{cat}_i)\\
&\ge  2 \cdot \mathds{1}(F^\complement_i)-1 - \mathds{1}(F^\mathrm{cat}_i).
\end{aligned}
\end{equation}
By Corollary~\ref{cor_success}, the random variables~$(\mathds{1}(F^\complement_i))_i$ stochastically dominate independent and identically distributed random variables~$(W_i)_i$ with~$W_i \sim \mathrm{Ber}(9/16)$. Noting that
\[\mathbb E[2W_i-1 ] = 2 \cdot \frac{9}{16}-1 = \frac{1}{8},\]
we obtain
\begin{equation}\label{eq_complements}
\mathbb P\left(\sum_{i=1}^{\lfloor \delta n\rfloor-1} ( 2 \cdot \mathds{1}(F_i^\complement)-1) > \frac{\delta n}{16}\right) \ge \mathbb P\left(\sum_{i=1}^{\lfloor \delta n \rfloor -1}(2W_i-1) > \frac{\delta n}{16} \right) \xrightarrow{n \to \infty} 1.
\end{equation}
Next, Lemma~\ref{lem_cats} (reducing~$\delta$ if necessary, as required by that lemma) gives
 \begin{equation}\label{eq_with_cats}
     \mathbb P\left(\sum_{i=1}^{\lfloor \delta n\rfloor-1}  \mathds{1}(F_i^\mathrm{cat}) \le \frac{\delta n}{32}\right) \xrightarrow{n \to \infty} 1.
 \end{equation}
 Now, using~$X_{\lfloor \delta n\rfloor -1} \ge \sum_{i=1}^{\lfloor \delta n\rfloor -1} \Delta X_i$ together with~\eqref{eq_Delta} and the convergences~\eqref{eq_complements} and~\eqref{eq_with_cats} completes the proof.
 \end{proof}

  \begin{proof}[Proof of Proposition~\ref{prop:suitable-finite-subgraph}]
Fix~$\delta$ small enough, as required by Proposition~\ref{prop_many_candidates}.
For each~$i \in \{0,\ldots,\lfloor \delta n\rfloor -1\}$, let~$v_i$ denote the pivot of round~$i$. Let~$\phi(i)$ denote the index of the round in which~$v_i$ entered the queue~$\mathbf{QCand}$. There are two possibilities:
\begin{itemize}
    \item if~$\mathbf{QCand}$ was empty when round~$i$ started, then~$\phi(i)=i$;
    \item otherwise,~$\phi(i)<i$.
\end{itemize}

Let~$\mathscr I$ denote the set of~$i \in \{0,\ldots,\lfloor \delta n\rfloor -1\}$ for which round~$i$ was a success. For each~$i \in \mathscr I$, we let~$\mathcal S_i$ denote the set of edges connecting~$v_i$ to its neighbours (note that the algorithm has necessarily revealed all these edges, since round~$i$ was a success). We also define a set of edges~$\mathcal L_i$ as follows. If~$\phi(i)=i$, we set~$\mathcal L_i=\emptyset$. If~$\phi(i)<i$, note that the round with index~$\phi(i)$ was also necessarily a success, since it was the round in which~$v_i$ entered~$\mathbf{QCand}$. Consider the graph induced by the set of edges revealed in round~$\phi(i)$; this graph is necessarily a tree. Then, let~$\mathcal L_i$ be the set of edges in the geodesic from~$v_{\phi(i)}$ to~$v_i$ in this graph.

As in the proof of Proposition~\ref{prop_many_candidates}, let~$X_i$ be the number of elements of~$\mathbf{QCand}$ at the end of round~$i$, and work on the event that~$X_{\lfloor \delta n \rfloor -1} \ge \delta n/32$, which has probability tending to~1 by that proposition.
Let~$i^*:=\max\{i \in \{0,\ldots,\lfloor \delta n \rfloor - 1\}: X_i = 0\}$, with the convention that~$i^*:=-1$ if this set is empty. Note that~$i^* < \lfloor \delta n \rfloor -1$. Round~$i^*+1$ starts with~$\mathbf{QCand}$ empty, so~$\phi(i^*+1)=i^*+1$; moreover, by maximality of~$i^*$, we have~$X_{i^*+1}>0$, so round~$i^*+1$ is a success. For all~$i \in \{i^*+2,\ldots, \lfloor \delta n \rfloor -1\}$, round~$i$ starts with~$\mathbf{QCand}$ non-empty, so~$\phi(i)<i$ (and in fact~$\phi(i) > i^*$).

We now define our star-line subgraph as the subgraph induced by the set of edges
\[
\bigcup_{i \in \mathscr I,\, i > i^*} (\mathcal S_i \cup \mathcal L_i);
\]
the set of vertices~$V'$ playing the role of star centers is~$\{v_i: i \in \mathscr I,\; i > i^*\}$. It is readily checked that this subgraph satisfies the properties in Definition~\ref{def_star_line}.

The number of elements of~$V'$ is equal to the number of successful rounds among~$i^*+1,\ldots,\lfloor \delta n \rfloor -1$. By~\eqref{eq_Delta}, round~$i^*+1$ increases~$X$ from~$0$ to~$2$, and for~$i \ge i^*+2$, round~$i$ increases~$X$ by~1 if it is a success and does not increase it otherwise. Since~$X_{\lfloor \delta n \rfloor -1} \ge \delta n/32$, we obtain~$|V'| \ge \delta n /32 - 1$.
 \end{proof}

\textbf{Acknowledgements:}  MS was supported by the German Research Foundation (DFG) Project ID: 531542160 within the priority programme SPP 2265.

\textbf{AI usage statement:} The authors used language models (OpenAI's ChatGPT and Anthropic's Claude) to aid in finding related literature, and to generate the TikZ code for all the figures. No mathematical content and no part of the text was done with the aid of AI.

\printbibliography

\end{document}